\documentclass[12pt]{article}

\usepackage{amssymb, cite}
\usepackage{graphicx}
\usepackage{graphicx}
\usepackage{amssymb, cite}
\usepackage{setspace}
\usepackage[hmargin=1.15in,vmargin=1.15in]{geometry} 

\newcommand{\be}{\begin{equation}}
\newcommand{\ee}{\end{equation}}
\newcommand{\bqn}{\begin{eqnarray}}

\newcommand{\eqn}{\end{eqnarray}}

\newcommand{\bd}{\begin{description}}
\newcommand{\ed}{\end{description}}

\newtheorem{stat}{}[section]

\def\bs{\begin{stat}}
\def\es{\end{stat}}
\def\ben{\begin{enumerate}}
\def\een{\end{enumerate}}

\def\bp{\noindent{\bf Proof}  \ }
\newcommand{\ep}{\hfill $\square$}
\def\epcl{\hfill $\diamondsuit $}

\makeatletter
\@addtoreset{equation}{section}

\makeatother

\begin{document}

\begin{center}
{\bf OPERATIONS ON GRAPHS
\\[1ex]
INCREASING SOME GRAPH  PARAMETERS}
\\[3ex]
{\large {\bf Alexander Kelmans}}
\\[2ex]
{\bf University of Puerto Rico, San Juan, Puerto Rico}
\\[0.5ex]
{\bf Rutgers University, New Brunswick, New Jersey}
\\[2ex]
\end{center}

\begin{abstract}

In this partly expository paper we  discuss and describe
some of our old and 
recent
results on  partial orders  on the set  
${\cal G}_n^m$  of graphs with $n$ vertices and $m$ edges and some operations on graphs within 
${\cal G}_n^m$ that are monotone with respect to these partial orders. 
The partial orders under consideration include those related with some Laplacian characteristics of graphs 
as well as with some
probabilistic characteristics of graphs with randomly deleted edges.
Section \ref{Preliminaries} contains  some basic facts  on the Laplacian polynomial of a graph. 
Section \ref{operations}  
describes various graph operation and  their properties. 
In  
Section \ref{partialOrders}  we introduce some partial orders $\succeq $ on 
${\cal G}_n^m$ related
 with the graph Laplacian  and  the graph reliability ({\em Laplacian posets} and {\em reliability posets}). 
Section \ref{Lposets} contains some 
old and recent results on the 
$\succeq $-monotonicity of some graph  operations with respect to Laplacian posets. 
Section \ref{Rposet} includes some 
old and recent results on the 
$\succeq $-monotonicity of some graph operations
with respect to reliability posets and some open problems. 
In Section \ref{other-parameters}  we consider some other parameters of graphs 
and  establish some results on $\succeq $-monotonicity of our graph  operations with respect to the  linear orders  $\succeq $ on ${\cal G}_n^m$ related with these parameters. The list of these parameters includes
the numbers of Hamiltonian cycles or paths and the numbers of forests of special type.
Section \ref{weighted-graphs} contains   some  generalizations of the described results to weighted graphs.

\end{abstract}

\section{Introduction}
\label{introduction}

\indent

 All notions and facts on graphs, that are  
used but not described here, can be found in 
\cite{BM,D}.
Let $\bar{{\cal G}}_n^m$ (${\cal G}_n^m$) denote the set of  graphs (respectively,  simple graphs) with 
$n$ vertices and $m$ edges.  Replacing in the above notations ${\cal G}$ by ${\cal C}$ results in the notations of the corresponding sets of connected graphs.
If $G, F \in {\cal G}_n^m$, we also say simply that $G$ and $F$ are {\em of the same size}.
\\[1ex]
\indent
In a series of papers (see, for example,
\cite{
Kmxtr,
Kproceed,
Krandgr,
Kextr,
Knotuniform,
KoperProb,
Kcmpr,
Krandgraph,
KoperTrees,
Kprob-network,
Kasympt})
we considered various aspects of the following combinatorial optimization problems related with the synthesis of reliable networks.
Let $G$ be a graph. Suppose that every edge of $G$ has probability $p$ to exist and that all the edge events are mutually independent. Let $R(p,G)$ denote the probability that the random graph $(G,p)$ is connected. We call $R(p,G)$ the {\em reliability function} (or just the {\em reliability}) of $G$.
The problem ${\cal R}_{max}$ is 
to find a most reliable  graph $M(p)$ in $\bar{{\cal G}}_n^m$, i.e. such that   
\\[0.7ex]
\indent
$~~~~~~~~~~~~~~~~~~~~~~~
R(p,M(p)) = \max \{R(p, G): G \in \bar{{\cal G}}_n^m\}.$
\\[0.7ex]
The corresponding dual problem ${\cal R}_{max}^*$ is to find a graph 
$M^*(p)\in \bar{{\cal G}}_n$ such that 
\\[0.7ex]
\indent
$~~~~~~~~~~~~~~~~~~~~~~~
e(M^*(p)) =
 \min \{e(G): G \in \bar{{\cal G}}_n,R(p,G) \ge \alpha (p)\}$.
\\[0.7ex]
It is also interesting to consider the problem ${\cal R}_{min}$ of finding a least reliable graph $L(p)$ in  $\bar{{\cal G}}_n^m$.
If $G$ is not connected, then $R(p,G) = 0$. For that reason, the non-trivial problem    ${\cal R}_{min}$  is
to find a graph $L(p) \in \bar{{\cal C}}_n^m$ such that
\\[0.7ex]
\indent
$~~~~~~~~~~~~~~~~~~~~~~~
R(p,L(p)) = \min \{R(p, G): G \in \bar{{\cal C}}_n^m\}$.
\\[0.7ex]
\indent
Let $t(G)$ denote the number of spanning trees of $G$.
For $p$ close to zero,  the problem ${\cal R}_{max}$ is 
equivalent to problems ${\cal T}_{max}$ of finding 
 a graph $M \in \bar{{\cal G}}_n^m$ such that
 \\[0.7ex]
\indent
$~~~~~~~~~~~~~~~~~~~~~~~~
t(M) = \max \{t(G): G \in \bar{{\cal G}}_n^m\}$
 \\[0.7ex]
 and ${\cal R}_{min}$ is equivalent to
the problem ${\cal T}_{min}$ of finding
a  graph $L$ in  $\bar{{\cal G}}_n^m$ such that
\\[0.7ex]
\indent
$~~~~~~~~~~~~~~~~~~~~~~~~
t(L) = \min \{t(G): G \in \bar{{\cal C}}_n^m\}$, 
where $m \ge n-1$.
 \\[0.7ex]
\indent
The corresponding dual problem ${\cal T}^*_{max}$ is to find a graph $M^*\in \bar{{\cal G}}_n$ such that 
\\[0.7ex]
\indent
$~~~~~~~~~~~~~~~~~~~~~~~~
e(M^*) = \min \{e(G): G \in \bar{{\cal G}}_n,
t(G) \ge \alpha \}$.
\\[0.7ex]
\indent
Although in general these problems are probably $NP$-hard, it turns out that they 
can be solved in some non-trivial particular cases and, in addition,  their analysis leeds to some interesting mathematical results,  ideas, and questions.

In particular, it is not hard to solve problems ${\cal R}_{max}$ and ${\cal R}_{min}$ for graphs of relatively small cyclomatic number or corank (see, for example, 
{\bf \ref{maxR(G),cycl(G)<4}} below).
In \cite{Kprob-network,Kasympt} we were able to give an asymptotically optimal solution to problem ${\cal R}^*_{max}$ as well as to problem ${\cal R}^*_{max}(k)$ on the probability  that $(G,p)$ is $k$-connected.
\\[1ex]
\indent
Since problems ${\cal R}_{max}$ and ${\cal R}_{min}$ have parameter $p \in [0,1]$, it is natural to introduce 
the following  relation on ${\cal G}_n^m$:
given $G, F \in {\cal G}_n^m$, let $G \succeq^r F$ if
$R(p,G) \ge R(p,F)$ for every $p \in [0,1]$.
In \cite{Kprob-network} we observed that there are graphs $G$ and $F$  of the same size that are not  $\succeq^r$-comparable as well as   non-isomorphic graphs $G$ and $F$  of the same size such that  $R(p,G) \equiv R(p,F)$. Therefore $\succeq^r$ is a partial quasi-order relation on ${\cal G}_n^m$.

Now the following  natural question on problem ${\cal R}_{max}$ is in order:
Does every ${\cal G}_n^m$ have a $\succeq^r$-maximum  
graph ? In other words,
does a most  reliable graph $M(p)$ in some ${\cal G}_n^m$ 
 depend on $p$ ? 
We gave a construction providing infinitely many pairs 
$(n,m)$ for which  ${\cal G}_n^m $ has no $\succeq^r$-maximum  \cite{Krandgr,Knotuniform} (see more details in {\bf \ref{non-uniform}} below), i.e. for each of these pairs $(n,m)$ there are non-isomorphic graphs $M(p_1)$ and $M(p_2)$ for some $0 <p_1 < p_2 < 1$. Some further interesting questions along this line are described below  in Section \ref{Rposet}.

Obviously, if $G$ is not connected, then $R(p,G) \equiv 0$.
Therefore, a similar non-trivial question about problem 
${\cal R}_{min}$ is:  
Does every ${\cal C}_n^m$ have a $\succeq^r$-minimum  
graph ? In other words,
does a least reliable graph $L(p)$  in some ${\cal C}_n^m$ depend on $p$ ? 
The answer to this question is not known. 
 More detailed about this question are given  
in Section \ref{Rposet}.
 \\[1ex]
\indent
The Matrix Tree Theorem (see, for example, 
\cite{BM,CDS} 
and 
{\bf \ref{MTT}} below) provides a beautiful algebraic formula for $t(G)$ and allows to find it in polynomial time.
On the other hand, the problem of finding $R(p, G)$ is $\#P$-hard.
For that reason the above problems 
on $R(p,G)$ are much harder to analyze than the problems on $t(G)$.
Therefore it was natural to expect that the Matrix Tree Theorem could be useful in developing adequate approaches to attack problem ${\cal T}_{max}$ on finding graphs of given size with the maximum number of spanning trees.
In early 60's motivated by this idea, we discovered that
the characteristic polynomial 
\\[1.5ex]
\indent
$L(\lambda ,G)
 = \sum \{(-1)^s c_s(G) \lambda ^{n - s}: 
s \in \{0, \ldots , n\}\}$
\\[1.5ex]
  of the Laplacian matrix from the Matrix Tree Theorem should play an important role in this regard. Accordingly, 
in\cite{Kpolinom, Ktree2,Ktree1,KuniGr}  
we undertook some study and established various combinatorial properties of   $L(\lambda ,G)$, its coefficients,  the Laplacian spectrum of $G$, and their relation with the above mentioned problems. 
\\[1ex]
\indent
Papers  \cite{Ktree1,Ktree2,Kpolinom} were published in Russian in the Soviet Union, \cite {Ktree1,Ktree2} were translated into English in 1966,  the results of our manuscript \cite{KuniGr} of 1963 were described by Cvetkovi\'c in  one of his papers in 1971 (see also \cite{CDS}), and later papers  
\cite{HK,Ktree-forest,KeigvMtr,KCh,Kcmpr,KoperProb,KoperTrees,
 Knotuniform,Kextr,Kmxtr,KgraphTrnsf} were published in the West, where according to \cite{CDS}
the interest to this topic became apparent in 70's. Unfortunately, till now some researchers working in this area are not aware of certain basic facts from those  papers and quite a few results from these papers have later been published 
again and again. 
\\[1.5ex]
\indent
 The following simple facts on $L(\lambda ,G)$ turned out to be pretty  useful. 
\bs {\em \cite{Kpolinom,Ktree1}}
\label{roots-trees}
Let $G \in \bar{{\cal G}}_n^m$ and
$(\lambda _0(G)\le \cdots \le \lambda _{n-1}(G))$ the list of all $n$ roots of $L(\lambda ,G)$
{\em (i.e. the spectrum of the Laplacian matrix of $G$)}. Then
\\[1ex]
$(a1)$ $0 = \lambda _0(G) \le  \cdots  \le \lambda _{n-1}(G) \le n$ and
\\[1ex]
$(a2)$ 
$ c_1(G) = \sum \{\lambda _i(G): i \in \{1, \ldots , n-1\}\} = 2 e(G) = 2m$ and 
\\[0.7ex]
\indent
$ c_{n-1}(G) = nt(G) = \prod \{\lambda _i(G): i \in \{1, \ldots , n-1\}\}$.
\es 

Let $K_n$ be the complete graph with $n$ vertices.
From {\bf \ref{roots-trees}} we have:
\bs {\em \cite{Kpolinom}}
\label{treesKn}
Let $G \in \bar{{\cal G}}_n^m$. Then
\\[1ex]
$(a1)$ $t(G) \le n^{-1}(2m/(n-1))^{n-1}$, and so
\\[1ex]
$(a2)$ if $m = (^n_2)$, then
$t(G) \le  n^{n-2} =  t(K_n)$.
\es

Thus, $(a2)$ in {\bf \ref{treesKn}} gives  a solution of problem ${\cal T}_{max}$ for 
$m =  (^n_2)$. 
Moreover, it turns out that 
$t(G) =  n^{n-2}$ if and only if $G = K_n$.
In \cite{LP} this solution was extended to
problem $R_{max}$ for $m =  (^n_2)$.
\\[1ex]
\indent
Here is another useful fact on $L(\lambda , G)$.
\bs {\em \cite{Ktree1}}
\label{inclu-exclu} Let $G \in {\cal G}_n^m$.
Then
$t(K_{n+s} - E(G)) = (n +s) ^{r-2} L(n+ s, G)$ for every 
$s \ge 0$.
\es

 It is interesting that the formula in 
{\bf \ref{inclu-exclu}} turns out to be the inclusion-exclusion formula for the number of spanning trees in $K_{n+s}$ avoiding the edges of its subgraph $G$ \cite{KCh} (see more details in {\bf \ref{incl-excl}} below).
Let $d(G, n +s) = t(K_{n +s}) - t(K_{n+s} - E(G)$.
Then $d(G, n +s)$ is the number of spanning trees in $K_{n+s}$ that are destroyed when the edges of $G$ are removed from $K_{n+s}$. For that reason we call $d(G, n +s)$ the 
{\em destroying ability of} an $n$-vertex graph $G$ in $K_{n+s}$.

Now by {\bf \ref{inclu-exclu}},  comparing simple $n$-vertex graphs by their number of spanning trees  is equivalent to comparing the destroying ability of their complements  in $K_n$. Furthermore, comparing the destroying abilities of $G$ in every complete graph $K_{n+s}$  is equivalent to comparing the polynomials 
$L(\lambda , G)$ in every integer point $\lambda = n + s$.
This suggests  the following partial quasi-order relation on 
${\cal G}^m$. 
Given $G, F \in {\cal G}^m$, let 
\begin{center}
$ G \succeq ^{\tau } F$ if 
$L(\nu + s ,G) \ge L(\nu + s ,F)$ 
\end{center}
for every integer $s \ge 0$, where $\nu = \max\{v(G),v(F)\}$.
Now it is clear that the study of properties of relation 
$\succeq ^{\tau }$ on ${\cal G}_n^m$ may help to solve problems 
${\cal T}_{max}$ and ${\cal T}_{min}$ for some classes of graphs.
\\[1.5ex]
\indent
In \cite{KoperProb} we found some $\succeq ^{\tau }$-increasing and   $\succeq ^r$-increasing operations on graphs in ${\cal G}_n^m$.
Using these operations we were able  to solve 
problems 
${\cal R}_{max}$ for $m \ge (^n_2) - \lfloor n/2   \rfloor $ and 
${\cal R}_{min}$ for  $m \ge (^n_2) - n +2 $ \cite{KoperProb} as well as problems ${\cal T}_{max}$ for 
$m \ge (^n_2) - \lfloor n/2   \rfloor $ 
and ${\cal T}_{min}$ for  $m \ge (^n_2) - n +2$ \cite{KCh}.
In all these cases there exist  a $\succeq ^r$-minimum graph ${\cal C}_n^m$ and a unique $\succeq ^r$-maximum graph in ${\cal G}_n^m$.
 Later we have found some more delicate $\succeq ^{\tau }$-increasing operations for some special classes of graphs 
 \cite{KgraphTrnsf} which allowed us to also solve 
problem ${\cal T}_{max}$ for $m \ge (^n_2) - n +2$ 
\cite{Kextr}. 
\\[1.5ex]
\indent
In this partly expository paper we  
discuss and describe some partial orders  on the set  
${\cal G}_n^m$  of graphs with $n$ vertices and $m$ edges and some operations on graphs within 
${\cal G}_n^m$ that are monotone with respect to these partial orders. The partial orders under consideration include those related  with some 
Laplacian characteristics of graphs as well as with some probabilistic characteristics of graphs having  randomly deleted edges.
\\[1.5ex]
\indent
In Section \ref{notation} we give necessary notions and notation, as well as  some simple observations.
In Section \ref{Preliminaries} we describe   some basic results  on the Laplacian polynomial of a graph. 
In Section \ref{operations}  we define various graph operation preserving the size of the graph 
and describe some simple and useful properties of these operations. 
In  Section \ref{partialOrders}  we introduce various partial orders on ${\cal G}_n^m$ related, in particular, with the Laplacian polynomial $L(\lambda ,G)$ 
and with the graph reliability $R(p, G)$ 
({\em Laplacian posets} and {\em reliability posets}) and establish some mutual properties of these relations. 
Section \ref{Lposets} contains some 
old and recent results on the 
$\succeq $-monotonicity of some graph  operations with respect to Laplacian posets. 
Section \ref{Rposet} includes some 
old and recent results on the 
$\succeq $-monotonicity of some graph operations
with respect to reliability posets and some open problems. 
In Section \ref{other-parameters} we consider 
linear orders $\succeq $ on ${\cal G}_n^m$ related with some other parameters of graphs and establish some results on $\succeq $-monotonicity of some graph  operations with respect to these  $\succeq $-orders.
The list of parameters considered in this section includes the numbers of Hamiltonian cycles or paths and the numbers of forests of special types, for example, the number of matchings of a given size.
Section \ref{weighted-graphs} contains   some  generalizations of the described results to weighted graphs.
\\[2ex]
\indent
Many results described in this paper were included in
our lectures on Algebra and Combinatorics in 
Rutgers University, 1992 - 1993 and
in University of Puerto Rico, 1995 - 2009.

\section{Notions, notation, and simple observations}
\label{notation}

\indent

All notions and facts on graphs that are  
used but not described here can be found in \cite{BM,D}.
\bs
\label{graph}
{\em Let $G = (V, E, \varphi )$ be a graph, where
$V = V(G)$ is the {\em set of vertices} of $G$,
$E = E(G)$ is the {\em set of edges} of $G$, and 
$\varphi =  \varphi _G$ is the function from $E$ to the set of unordered pairs of vertices of $G$ (the {\em incident function} of $G$). A graph is called {\em simple} if it has no loops and no parallel (or, the same, multiple) edges.

Let $v(G) = |V(G)|$ and $e(G) = |E(G)|$.
We  say that graphs $G$ and $F$ are {\em of the same size} if $v(G) = v(F)$ and $e(G) = e(F)$.
Let $Cmp(G)$ denote the set of components of $G$ and $cmp(G) = |Cmp(G)|$. 
Let $r(G) = v(G) - cmp(G)$ and 
$r^*(G) = e(G)- r(G) = e(G) - v(G) + cmp (G)$.
The parameter $r(G)$ is called the {\em rank} of $G$ and $r^*(G)$ is called the {\em corank} (or the {\em cyclomatic number}) of $G$.
Let $isl(G)$ denote  the number of isolated vertices of  $G$.}
\es
\bs
\label{graphset}
{\em 
Let $\bar{{\cal G}}$, $\bar{{\cal G}}_n$ and $\bar{{\cal G}}^m$ denote the sets of all graphs,  graphs with $n$ vertices, and graphs with $m$ edges, respectively, and let 
$\bar{{\cal G}}_n^m = {\bar{\cal G}}_n \cap {\bar{\cal G}}^m$. Replacing in the above notation $\bar{{\cal G}}$ by 
${\cal G}$ (${\cal C}$) gives the corresponding sets of simple graphs (respectively, connected simple graphs).
Let ${\cal T}_n$ denote the set of trees with $n$ vertices. 
Let ${\cal F}(G)$ and ${\cal T}(G)$ denote the sets of spanning forests and spanning trees of $G$, and 
accordingly, $f(G) = |{\cal F}(G)|$ and 
$t(G) = |{\cal T}(G)|$.
}
\es
\bs
\label{[X,Y]}
{\em 
For $X, Y \subseteq V(G)$ let $[X,Y]$ denote the set of 
edges of $G$ with one end-vertex in $X$ and the other end-vertex in $Y$.
For $x \in V(G)$, let $N(x,G) = \{v \in V(G): xv \in E(G)\}$, 
$D(x,G) = [x,V(G)]$, and $d(x,G) = |D(x,G)|$. We call 
$d(x,G)$ the {\em degree of a vertex $x$ in $G$}.
Let $\Delta (G) = \max\{d(x,G): x \in V(G)\}$ and 
$\delta (G) = \min\{d(x,G): x \in V(G)\}$.}
\es
\bs
\label{complement}
{\em 
Let $K = (V, E)$ be the graph such that $E = {V \choose 2}$. This graph is called the {\em simple complete graph with the vertex set $V$}. We put $K = K_n$ if $|V| = n$.
If $G$ is a subgraph of $K$, then 
$[G]^c = K - E(G)$  is called the {\em simple complement of} $G$.
In particular, $[K]^c$  is the {\em $($edge$)$ empty graph} with the vertex set $V$.}
\es
\bs
\label{plus,times}
{\em 
Given two disjoint graphs $G$ and $F$, let
$G + F = G \cup F$ and $G \times F$ be the graph obtained from $G \cup F$ by adding the set of edges 
$\{gf: g\in V(G), f \in V(F)\}$  \cite{Ktree1,Ktree2}.
Obviously, $[G + F]^c = [G]^c \times [F]^c$.
If $G$ consists of $k$ disjoint copies of a graph $F$, we write $G = kF$.
A graph $G$ is called {\em decomposable} if $G = A + B$ or 
$G = A \times B$ for some disjoint graphs $A$ and $B$.
We call a graph $G$  {\em totally decomposable} 
 \cite{Ktree1,Ktree2} (see also \cite{HK}) 
 if $G$ can be obtained from one vertex graphs by a series of operations $+$ and $\times $. The notion of a  totally decomposable graph turned out to be so natural that it was later reintroduced again and again under different names:  a {cograph} in \cite{Lrch}, a {\em hereditary Dacey graph} in \cite{Sum}, a {\em $D^*$-graph} in \cite{Jung}, etc.  
 A totally decomposable graph $G$ can be naturally described by the so called {\em $(+,\times)$-decomposition tree of}  $G$ introduced in \cite{Ktree1,Ktree2} (see also \cite{HK}) and later in  \cite{CLBS} under the name the {\em cotree} of $G$. Similar notions for weighted graphs were introduced in \cite{KeigvMtr}.
 
 In \cite{Ktree2} we gave a very simple procedure for finding formulas for the Laplacian polynomial and spectrum (and, in particular the number of spanning trees) of  a totally decomposable graph $G$ in terms of the parameters of the  $(+,\times)$-decomposition tree of $G$ 
 (see also \cite{HK})}.
 \es
 
 From this procedure we have, in particular: 
\bs {\em \cite{Ktree1} }
\label{LaplacOfDecmpsG} 
Every totally decomposable graph has an integral spectrum. 
\es

Similar results for multigraphs and weighted graphs
(and their Laplacian matrices) can be found in 
\cite{KeigvMtr}.
\bs
\label{thresholdgraph}
{\em 
We call a graph $G$ {\em vertex comparable} 
 \cite{KelLatvia,Kcourse,Kproceed} if 
 $N(x,G) - y \subseteq N(y)$
 or  $N(y,G) - x \subseteq N(x)$
 for every pair $(x,y)$ of vertices $x$ and $y$ in $G$.
 A graph $G$ is called {\em threshold} \cite{CH} if $G$ has no induced subgraph isomorphic to $P^3$, $O^4$, and 
 $2P^1$. Let ${\cal F}_n^m$ denote the set of threshold graphs with $n$ vertices and $m$ edges. It is easy to prove that  the following is true.}
 \es
\bs
\label{comparable-threshold} 
$G$ is vertex comparable if and only if $G$ is threshold. 
\es

It is also easy to prove the following:
\bs
\label{thresholdG-totalDecompos} 
A threshold is totally decomposable. 
Moreover, if $H$ is a threshold graph with $n \ge 2$ vertices, then there exists a threshold graph $H'$ with $n-1$ vertices such that $H = H' + g$ or $H = H' \times g$, where $g$ is a one vertex graph. 
 \es
 
 The above property provides a simple {\em recursive description} of a threshold graph with at least two vertices. 
 
 The Laplacian spectrum and Laplacian polynomial of a  threshold graph have some special properties. For example, by {\bf \ref{LaplacOfDecmpsG}}, every threshold graph has an integral spectrum. 
\bs
\label{extreme-threshold-graph}
{\em 
Now we will define some special threshold graphs
which we call {\em extreme} \cite{Kcourse,Kproceed}. 
Let $(k,r,s)$ be a triple of non-negative integers such that $r < s$. Let $F(k,r,s)$ denote the graph obtained from the complete graph $K_s$ with $s$ vertices as follows:
fix in $K_s$ a set $A$ of $r$ vertices and a vertex $a$ in $A$, add to $K_s$ a new vertex $c$ and the set 
$ \{cx: x \in A\}$ of new edges to obtain graph $C(r,s)$,
add to  $C(r,s)$ the set $B$ of $k$ new vertices and 
the set  $\{az: z \in B\}$ of new edge to obtain graph 
$F(k,r,s)$.
Let ${\cal H}(k,r,s)$ denote the set of  all graphs $H$ 
obtained from $C(r,s) \cup B$ by adding a tree on the vertex set $B \cup \{a\}$ (and so this tree has $k$ edges).
Clearly, $F(k,r,s)$ is a threshold graph and 
$F(k,r,s) \in {\cal H}(k,r,s)$.
Obviously, ${\cal H}(k,r,s) = \{F(k,r,s)\}$ if and only if  $k = 0$.
}
\es

\begin{figure}
\begin{center}
 \includegraphics[scale=0.80]
 {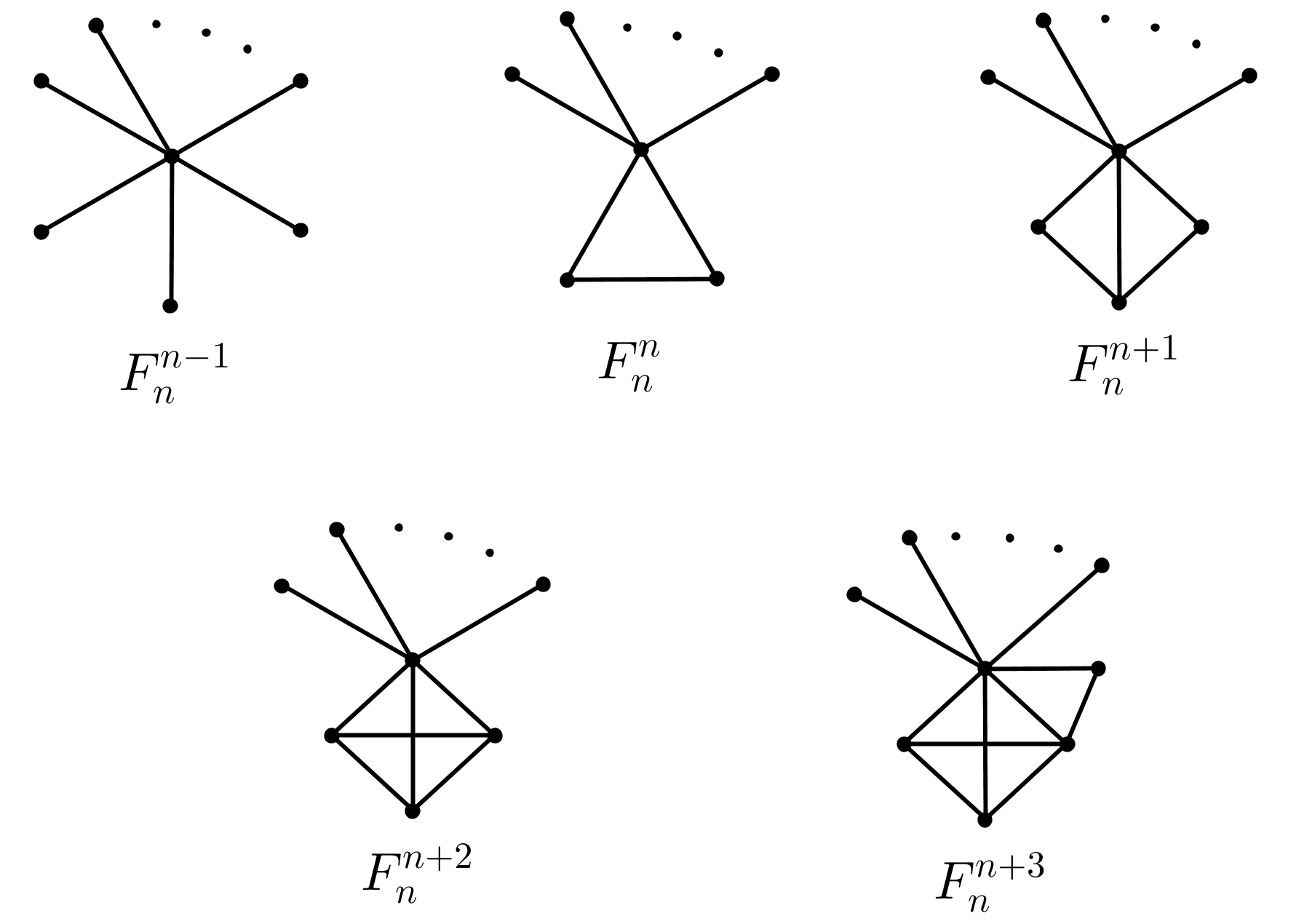}
 \end{center}
 \caption{Connected threshold graphs with $m \le  n+3$}
 \label {SmalThreshG}
\end{figure}

Let, as above, ${\cal C}_n^m$ be the set of simple connected graphs with $n$ vertices and $m$ edges.
It is easy to prove the following.
\bs
\label{T(k,s,r)} 
For every pair $(n, m)$ of integers such that
${\cal C}_n^m \ne \emptyset $ 
there exists  a unique triple 
$(k,r,s)$ of non-negative integers such that $r < s$ and 
$F(k,r,s) \in {\cal C}_n^m$.
\es
\bs
\label{F_n^m}
{\em 
If $F(k,r,s) \in {\cal C}_n^m$, we put  $F(k,r,s) = F_n^m$ and ${\cal H}(k,r,s) = {\cal H}_n^m$, and so  
${\cal H}_n^m \subset {\cal C}_n^m$.
We call $F_n^m$ the {\em extreme threshold graph} in 
${\cal C}_n^m$. 
 \\[0.7ex]
\indent
Obviously, $F_n^{n-1} = S_n$ for $n \ge 2$, 
$F_3^3 = \Delta $, 
$F_n^n$ with $n \ge 4$ is obtained from disjoint triangle 
$\Delta $ and the $(n - 3)$-edge star $S$ by identifying its center  with a vertex in $\Delta $ 
(and so $F_n^n$ with $n \ge 4$  is  $W_n$), $F_4^5 = K_4^-$, $F_n^{n+1}$ with $n\ge 5$ is obtained from disjoint $K_4^-$ and the
$(n - 4)$-edge star $S$ by identifying its center  with a vertex of degree three in  $K_4^-$, $F_4^6 = K_4$, 
$F_n^{n+2}$ with $n\ge 5$ is obtained from disjoint $K_4$ and the $(n - 4)$-edge star $S$ by identifying its center  
with a vertex  in  $K_4$ (see Figure \ref{SmalThreshG}).} 
\es

 It is easy to prove the following:
\bs
\label{ThreshGuniqueIfm<2n-3} 
Let $n$ and $m$ be natural numbers.
Suppose that $n-1\le m \le 2n - 3$. Then there exists only one connected threshold graph with $n$ vertices and $m$ edges, and so
${\cal F} _n^m = \{F_n^m\}$.
\es

\bs
\label{specialgraphs}
{\em 
We need notation for some special graphs:
$P$ is a path, $O$  is a cycle, 
$K_{1,n}$ is called a {\em star} $S$ (if $n\ge 2$, then a vertex of degree $n$ is the {\em center} of $S$ and every other vertex is a {\em leaf} of $S$; if $n = 1$, then every of two vertices of $S$ is a {\em leaf} and a {\em center}),
$Z$ is obtained from a star  $S$ with $e(S) \ge 2$ by adding a new vertex $x$ and a new edge between $x$ and a leaf 
of $S$, $K_4^-$ is obtained from $K_4$ by removing one edge, $\Delta $ is  the triangle, and
$W$ is obtained from a star $S$ with $e(S) \ge 3$ by adding an edge between two leaves of $S$.

Using the above operations ``+'' and `$`\times $'' on graphs (see {\bf \ref{plus,times}}) we have, in particular:
$mP^1$ is a matching with $m$ edges and 
$P^2 + (m-2)P^1$ is the disjoint union of the two-edge path and a matching with $m-2$ edges.}
\es
\bs
\label{D(r),Y}
{\em 
Let ${\cal D}(r)$ denote the sets of 
trees having diameter $r$.
We call a graph $F$ a {\em star-forest} if every component of $F$ is a star with at least one edge.
\\[0.7ex]
\indent
Obviously, $T \in {\cal D}(3)$ if and only if $T$ is obtained from a star-forest with two components by connecting their centers by an edge.
\\[0.7ex]
\indent
Also,  $T \in {\cal D}(4)$ if and only if
$T$ is obtained from a star-forest 
$F$ with  at least two components having two or more edges by specifying a leaf  for every star  and identifying all specified leaves with a new vertex. 
Let $S_1, \ldots, S_k$ be the components of $F$ with
$e(S_1) \le  \ldots \le e(S_k)$ and let 
$u(T) = (e(S_1), \ldots , e(S_k))$.
Obviously, every tree $T$ in ${\cal D}(4)$ 
is uniquely defined (up to isomorphism)  by $u(T)$.
\\[0.7ex]
\indent
Let $P$ be a path with $r \ge 2$ edges and $F$ a star-forest with at most $r - 1$ components.
Then there is an injection $\xi $ from $Cmp(F)$ to the set of non-leaf vertices of $P$.
Now let $Y$ be a tree  obtained from disjoint $P$ and $F$
by identifying the center of each component  $C$ of $F$
with vertex $\xi (C)$ in $P$.
Obviously, $Y \in {\cal D}(r)$.
A tree $Y$ obtained this way is called a {\em caterpillar}.
Let ${\cal K}(r)$ denote  the set of caterpillars having diameter $r$, and so ${\cal K}(r) \subseteq {\cal D}(r)$.
\\[0.7ex]
\indent
Let $K(r)$ be the graph obtained from disjoint path $P$ 
with $r \ge 2$ edges and a star $S$  by identifying 
a center vertex  of $P$ and a center of $S$.
Clearly, $K(r)$ is a caterpillar 
 and $K(r) \in {\cal D}(r)$.
Let $K_n(r)$ be graph $K(r)$ with $n$ vertices
(see Figure \ref{figK(r)}).
\begin{figure}
\vskip -5mm
\begin{center}
\includegraphics[scale=0.55]{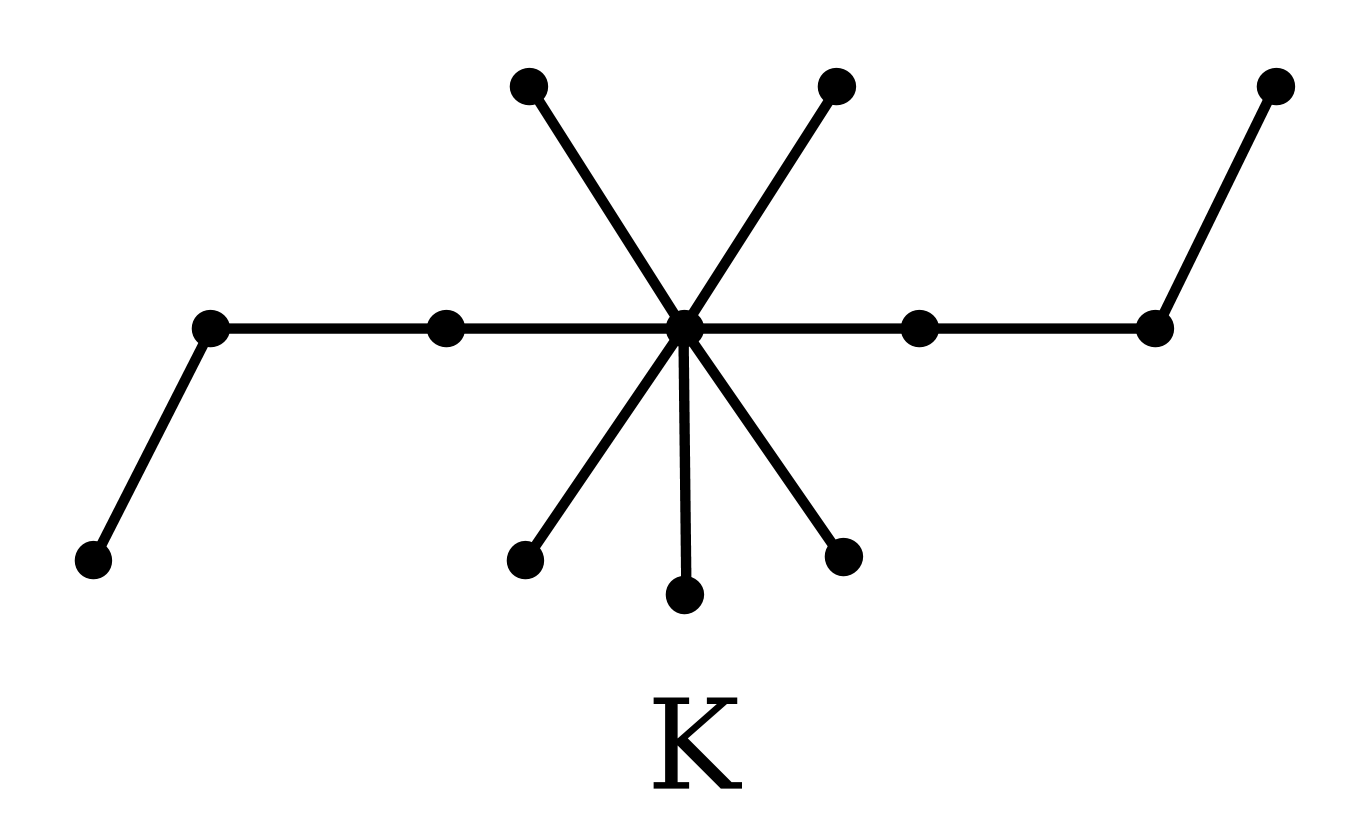}
\end{center}
\caption{$K ~= ~K_n(r)$, where 
~$v(K) ~= ~n = 12$ ~and ~
$diam(K) ~=~ r = 6$.}
\label{figK(r)}
\end{figure}
%
Let
${\cal D}_n(r)$ denote the set of $n$-vertex graphs in 
${\cal D}(r)$,  and ${\cal K}_n(r)$ the set of $n$-vertex caterpillars having diameter $r$.}
\es
\bs
\label{L(r),S(r)}
{\em 
Let  ${\cal L}(r)$ denote the sets of trees having $r$ leaves. 

Let ${\cal S}(r)$, $r \ge 3$,  denote the set of 
trees $T$ such that $T$ has exactly one vertex of degree $r$ and every other vertex in $T$ has degree at most two, and so ${\cal S}(r) \subseteq  {\cal L}(r)$.
If $r \ge 2$, then we call the vertex of degree $r$ in $T$ the {\em root}  of $T$. 
In other words, a tree $T \in  {\cal S}(r)$ 
if and only if it can be obtained from $r$ disjoint paths
$P_i$,  $1\le i \le r$, having at least two vertices by specifying one end-vertex of each path and identifying these specified  end-vertices of all paths.
Let $e(P_1) \le  \ldots \le e(P_r)$ and
put $w(T) = (e(P_1), \ldots , e(P_r)) $. 
Obviously, every tree $T$ in ${\cal S}(r)$ 
is uniquely defined (up to isomorphism)  by $w(T)$.

Let $M(r)$ denote the tree  $T$ in ${\cal S}(r)$ such that
every $e(P_i)$ in $T$, except possibly for $e(P_r)$, equals  one.

%
Let $L(r)$ denote the tree  $T$ in ${\cal S}(r)$ with the property: $|e(P_i) - e(P_j)| \le 1$ for every
 $1\le i, j \le r$. 
 
 %
Let ${\cal S}_n(r)$ denote the set of $n$-vertex graphs in 
${\cal S}(r)$. We denote the $n$-vertex trees $M(r)$ and 
$L(r)$ by $M_n(r)$ and $L_n(r)$, respectively. 
\\[0.7ex]
\indent
For two trees $T$ and $D$ in ${\cal S}_n(r)$, 
let $T > ^w  D$ if 
$w(T)$ lexicographically less than $w(D)$. 
Obviously, $>^w$ is a linear order on ${\cal S}_n(r)$.
In particular, if $M_n(r )>^w T_n >^w L_n(r)$ for every 
$T_n \in {\cal S}_n(r) \setminus \{M_n(r) , L_n(r)\}$.}
\es

\begin{figure}
\begin{center}
\includegraphics[scale=0.55]{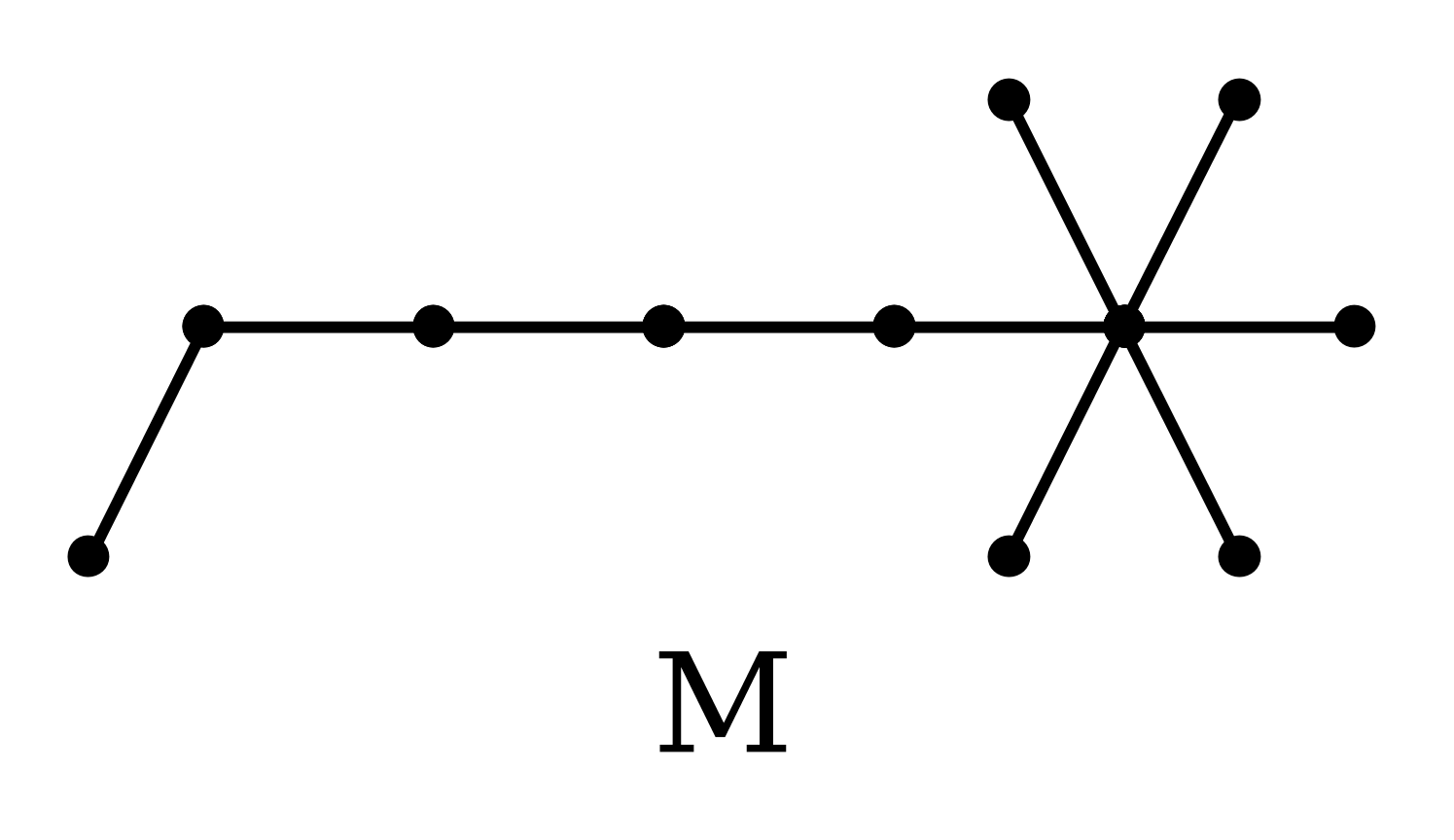}
\caption{$M~ = ~M_n(r)$, ~where
~ $v(M) ~= ~n = 11$ ~and~ 
$lv(M) ~= ~r = 6$.
}
\end{center}
\label{figM(r)}
\end{figure}


 \begin{figure}
\vskip -5mm
\begin{center}
\includegraphics[scale=0.55]{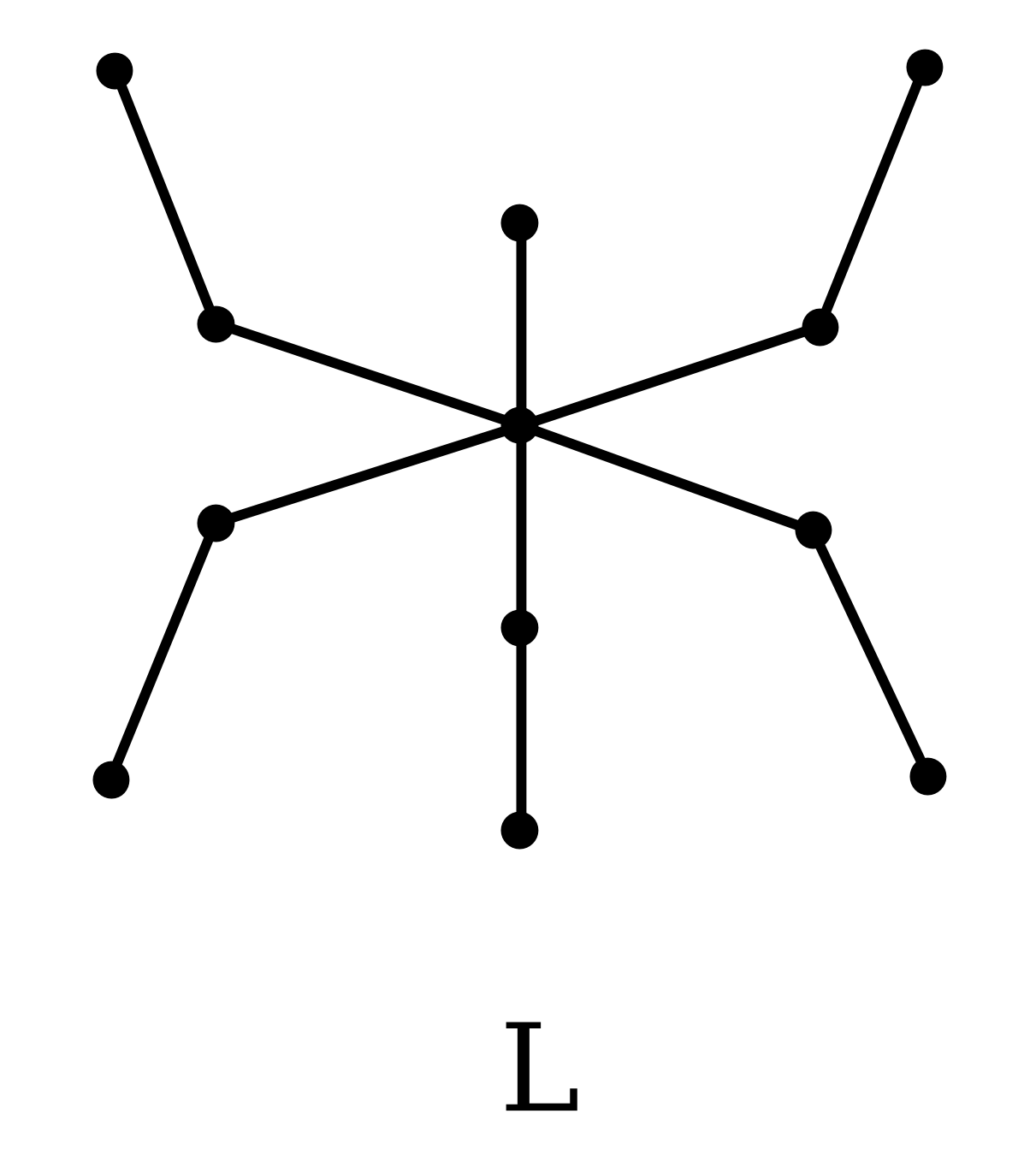}
\vskip 2mm
\end{center}\caption{$L = L_n(r)$, ~where  ~$v(L) = n = 12$ ~and ~$lv(L) ~= ~r = 6$.}
\label{figL(r)}
\end{figure}


Examples of graphs $M_n(r)$ and $L_n(r)$ are given on Figures \ref{figM(r)} and \ref{figL(r)}.
\bs
\label{symmetric-polynomials}
{\em 
Given  $x = (x_1, \ldots, x_n) \in \mathbb{R}^n$ and a permutation
$\sigma : \{1, \ldots , n\} \to \{1, \ldots , n\}$, let
$\sigma [x] = (x_{\sigma (1)}, \ldots, x_ {\sigma (n)})$.
 A function $f: \mathbb{R}^n \to \mathbb{R}$ is called {\em symmetric} if 
$f(x) = f(\sigma [x])$ for every $x \in \mathbb{R}$ and every permutation $\sigma : \{1, \ldots , n\} \to \{1, \ldots , n\}$. Let ${\cal SF}$ denote the set of symmetric functions.
\\[1ex]
\indent
Let $z$ and $n$ be positive integers, $z \le n$,  and 
$X = \{x_1, \ldots , x_n\} $, where each $x_i$ is a real number. 
Let
\\[1.5ex]
\indent
$\sigma _z(X) = \sum \{\prod\{x_i: i \in Z\}: Z \subseteq \{1, \ldots ,n\}, |Z| = z\}$. 
\\[1.5ex]
Function $\sigma _z(X)$ is called {\em the elementary symmetric polynomial of degree $z$ in 
the variables from $X$}. Let 
$\sigma (\emptyset) = 1$.}
\es
\bs
\label{convex,concave-functions}
{\em 
A symmetric functions $f: \mathbb{R}^n \to \mathbb{R}$ is said to be {\em concave} if  it has  
the following property 
\\[0.7ex]
(${\bf \cap }$):
for every $r, s \in \{1, \ldots , n\}$, $r \ne s$, and $\varepsilon \ge 0$, if $x_r \le x_s$ and  each $x_i \ge 0$, then
\\[0.7ex]
\indent
$f(\{x_i: i \in  \{1, \ldots , n\} \setminus \{x_r\} \cup \{x_r +\varepsilon\}\}) \ge 
f(\{x_i: i \in  \{1, \ldots , n\} \setminus \{x_s\} \cup \{x_s +\varepsilon\}\})$. 
\\[0.7ex]
It is easy to see that every elementary symmetric polynomial is concave. 
\\[1ex]
\indent
A symmetric functions $f: \mathbb{R}^n \to \mathbb{R}$ is said to be {\em convex} if  it has  
the following property  
\\[0.7ex]
(${\bf \cup }$):
for every $r, s \in \{1, \ldots , n\}$, $r \ne s$, and $\varepsilon \ge 0$, if $x_r \ge x_s$ and  each $x_i \ge 0$, then
\\[0.7ex]
\indent
$f(\{x_i: i \in  \{1, \ldots , n\} \setminus \{x_r\} \cup \{x_r +\varepsilon\}\}) \ge 
f(\{x_i: i \in  \{1, \ldots , n\} \setminus \{x_s\} \cup \{x_s +\varepsilon\}\})$. 
\\[1ex]
For  $x = (x_1, \ldots, x_n) \in \mathbb{R}^n$, let
$\delta _k(x) = 
\sum \{x_i^k:  i \in  \{1, \ldots , n\}\}$. 
Obviously, function 
\\[0.7ex]
$\delta _k : \mathbb{R}^n \to \mathbb{R}$ is convex. 
If $f  \equiv c$ for $c \in  \mathbb{R}$, then $f$ is both convex and concave.} 
\es

\section{Preliminaries on Laplacian parameters of graphs}
\label{Preliminaries}

\indent
  
Let $G$ be a  graph with possible parallel edges but with no loops and let
$V = V(G) = \{v_1, \ldots , v_n\}$.
Let $A(G)$ be the the symmetric $n\times n$-matrix $(a_{ij})$, where each $a_{ii} = 0$ and each $a_{ij}$, $i \ne j$,  is the number of parallel edges with the end-vertices $v_i$ and $v_j$. Let $D(G)$ be the (diagonal) $n\times n$-matrix
$(d_{ij})$, where  each $d_{ii} = d(v_i, G)$ and  
$d_{ij} = 0$ for $i \ne j$.
Let $L(G) = D(G) - A(G)$.
Matrix   $L(G)$ is called the {\em Laplacian matrix of} $G$.
Let $X \subset V$. 
We need the following notation: 
\begin{itemize}
\item $G_X$ is the graph obtained from $G$ by identifying all vertices from $X$ and removing all loops (that may appear as a result of such identification),

\item $L_X(G)$ denotes the matrix obtained from $L(G)$ by removing the rows and columns corresponding to every 
vertex $x \in X$, and so in particular, $L_x(G)$ is the matrix obtained from $L(G)$ by removing the row and column corresponding to vertex $x$ of $G$, and

\item $f(G, X)$ denotes the number of spanning forests $F$ of 
$G$ such that every component of $F$ has exactly one vertex in $X$, and so  the number $t(G)$ of spanning trees of $G$ is equal to $f(G, x)$ for every $x \in  V$.
\end{itemize}

We start with the following classical Matrix Three Theorem
(see, for example, \cite{BM,CDS}).
\bs
\label{MTT} Let $G$ be a graph  with possible parallel edges. 
Then
$t(G) = det (L_x(G))$ for every vertex $x$ in $G$.
\es

From {\bf \ref{MTT}} we have the following  generalization.
\bs {\em \cite{Kcourse,Kpolinom}}
\label{generMTT}
Let $G$ be a graph, $V = V(G)$,  and $X \subseteq V(G)$.
Then
\\[1ex]
\indent
$t(G_X) = f(G,  X) = det (L_X(G))$.
\es

Here is a more general version of the Matrix Tree Theorem. Let $L_{(i,j)}(G)$ denote the matrix obtained from $L(G)$ by removing $i$-th row and $j$-th column.
\bs {\em \cite{Kcourse}}
\label{MTT2} Let $G$ be a graph  with possible parallel edges. 
Then
\\[1ex]
\indent
$t(G) = (-1)^{i+j}det (L_{(i,j)}(G))$ for every  $i,j \in V(G)$,
\\[1ex]
 i.e. $t(G)$ equals every cofactor of $L(G)$.
\es

Here is yet another version of the Matrix Tree Theorem. 
Let, as before,  $V = V(G) = \{v_1, \ldots , v_n\}$ and let 
$\bar{L}(G)$ be the $(n+1) \times (n+1)$-matrix obtained from $L(G) = (l_{i,j})$ by adding the elements 
$l_{i, n+1} = 1$ if $i \in \{1, \ldots , n\}$, 
$l_{n+1, j} = 1$ if $j \in \{1, \ldots , n\}$,  and
$l_{n+1, n+1} = 0$.
\bs {\em \cite{HKS}}
\label{MTT3}
Let $G$ be a graph  with possible parallel edges. Then
$t(G) = n^{-2}det (\bar{L}(G))$.
\es

Let $L(\lambda ,G) = det (\lambda I - L(G)$ 
and $Spctr(G)$ be the multi-set of the eigenvalues of $L(G)$.
It is easy to see the following.
\bs {\em \cite{Kcourse,Kpolinom}}
\label{posemdef}
Let $G \in \bar{{\cal G}}_n$. Then
\\[0.7ex]
$(a1)$ 
$L(G)$ is a positive semi-definite matrix and
 $det (L(G)) = 0$, and so
\\[1ex]
$(a2)$ all eigenvalues of $L(G)$ are non-negative real numbers$:$
\\[1ex]
\indent
$Spctr(G) = (0 = \lambda _0(G) \le \lambda _1(G) \ldots  \le \lambda _{n-1}(G) )$, and
\\[1ex]
$(a3)$
$L(\lambda , G) = \lambda P(\lambda , G)$, where $P(\lambda , G)$ is a polynomial of degree $n-1$
with the root sequence
$(\lambda _1(G) \le \ldots  \le \lambda _{n-1}(G))$ {\em (we denote it by $Sp(G)$)}$:$ 
\\[1ex]
\indent
$P(\lambda , G) = \sum \{(-1)^s c_s(G) \lambda ^{n- 1 - s}: 
s \in \{0, \ldots , n-1\}\}$, 
\\[1ex]
where
$c_s(G) = \sigma _s(Sp(G))$ for  $0 \le s\le n-1$.
\es

Let, as above,  $\Delta (G_n)$ and $\delta (G_n)$ denote  the maximum and the minimum vertex degree of $G_n$, respectively,  and $\lambda (G_n) = \lambda _{n-1}(G_n)$.
\bs {\em \cite{Kcourse,KCh,Kpolinom}}
\label{lambda-bounds}
$G \in \bar{{\cal G}}_n$. Then
\\[1ex]
$(a1)$ 
$\lambda (G) \le 
\max \{d(x,G) + d(y,G): x,y \in V(G), x \ne y\}$,
\\[1ex]
$(a2)$ 
$\lambda _i(G) \ge \lambda _i(G - e)$ for every
$e \in E(G)$ and $i \in \{1,  \ldots , n-1\}\}$, and so
\\[1ex]
$(a3)$ 
$\lambda (G) \ge \Delta (G) + 1$ and
$\lambda _1 (G) \le \delta (G)$. 
\es

We call $P(\lambda , G)$ the {\em Laplacian polynomial of $G$}.
 \\[1ex]
\indent
For a graph  $F$, let 
$\gamma (F) = \prod \{ v(C): C \in Cmp(F)\}$ if $F$ is a forest and $\gamma (F) = 0$, otherwise.
Recall that ${\cal F}(G)$ is the set of spanning forests of $G$.
\\[1ex]
\indent
Using {\bf \ref{generMTT}}, we obtained the following 
combinatorial interpretation of the coefficients of 
$P(\lambda , G)$.
\bs {\em \cite{Kcourse,KK,Kpolinom,KCh}} 
\label{coefficients}
Let $G$ be a  graph with $n$ vertices, $V = V(G)$, $s$ an integer,  and
$0 \le s \le n - k$, where $k$ is the number of components of $G$.
Then
\\[1ex]
\indent
$c_s(G) = \sum\{t(G_ {V - S}) : S \subseteq V, |S| = s\} =
\\[1ex]
\sum\{f(G, V - S) : S \subseteq V, |S| = s\} = 
\sum\{\gamma (F): F \in {\cal F}(G), e(F) = s\}$.
\es

Obviously,  $c_s(G) = 0$ for $s \ge \min \{e(G), v(G) - k\}$.
\\[1.5ex]
\indent
Let
$\nabla (G)$ denote the number of triangles of $G$ 
and $\delta _i[G] = \sum \{d(v,G)^i: v \in v(G)\}$.
From  {\bf \ref{coefficients}} we have, in particular:
\bs {\em \cite{Kcourse,KCh}}
\label{somecoefficients}
Let $G$ be a  graph with  $n$, and $m$ edges.  Then 
\\[1ex]
$(a0)$
$c_0(G) = 1$, 
\\[1ex]
$(a1)$
$c_1(G) = \delta _1[G] = 2m$, 
\\[1ex]
$(a2)$
$c_2(G) = 2m^2 - m - \frac{1}{2}\delta _2[G]$, and
\\[1ex]
$(a3)$
$c_3(G) = \frac{4}{3}m^3 - 2m^2 - (m-1)\delta _2[G]  + \frac{1}{3}\delta _3[G] - 2\nabla (G)$, 
and
\\[1ex]
$(a4)$
$c_{n-1}(G) = n t(G) = (-1)^{n-1}P(0,G)$.
\es
 
The coefficients of the polynomial $P(\lambda , G)$ satisfy the following recursion.
\bs {\em \cite{Kcourse,KCh}}
\label{recursion}
Let $G$ be a  graph with $m$ edges.
Then 
\\[1ex]
\indent
$(m - s) c_s(G) = \sum \{c_s(G - e): e \in E(G)\}$
for $s \in \{0, \ldots , m\}$.
\es

Given $G \in {\cal G}_n^m$,
let $\Phi (\lambda , G) = \lambda ^{m - n +1} P(\lambda , G) =  \lambda ^{m - n} L(\lambda , G)$. This modification of the Laplacian polynomial of a graph has the following useful property.
\bs {\em \cite{Kcourse,KCh}}
\label{propertyPhi}
Let $G$ be a graph and $G'$ obtained from $G$ by adding some isolated vertices.
Then $\Phi (\lambda , G) = \Phi (\lambda , G')$.
\es

Using {\bf \ref{recursion}}, we obtained the following recursion for $\Phi (\lambda , G)$.
\bs {\em \cite{Kcourse,KCh}}
\label{recursionPhi} 
$\Phi (\lambda , G) = \Phi (a , G) +
\sum \{\int_a^\lambda \Phi (t, G - u)dt : u \in E(G)\}$.
\es 

In \cite{Ktree1} we proved the following important and frequently used Reciprocity Theorem (that was later rediscovered and published several times).
\\[1ex]
\indent
Let, as above,  $[G]^c = K_n - E(G)$, where $G \in {\cal G}_n$.
\bs 
{\em \cite{Kcourse,Ktree1}}
\label{reciprocity}
Let  $G$ be a simple graph with $n$ vertices.
Then
\\[1ex]
$(a1)$ $\lambda _i(G) + \lambda _{n-i}([G]^c) = n$
for every $i \in \{1, \ldots , n-1\}$ or, equivalently,
\\[1ex]
$(a2)$
$P(\lambda , [G]^c) = (-1)^{n-1} P(n - \lambda , G)$. 
\es

Reciprocity Theorem {\bf \ref{reciprocity}} is a particular case of the following Reciprocity Theorem for so called directed balanced graphs \cite{KelLatvia}.
\\[1ex]
\indent
A {\em simple directed graph $($or digraph$)$} $D$ is a pair $(V, E)$, where $V$ is a non-empty set 
 and
$E \subseteq [V]^2$, where $[V]^2 = \{V\times V \setminus \{(x,x): x \in V\}$ (and so  $D$ has no parallel edges 
\\[0.3ex]
and no loops).
A digraph $\vec{K}_\circ = (V,  [V]^2)$ is called a 
{\em simple complete digraph} with the vertex set $V$.
A digraph 
\\[0.3ex]
$D^c = (V,  [V]^2 \setminus E) = \vec{K}_\circ \setminus E$ is called 
the {\em complement of } $D = (V, E)$.
\\[0.3ex]
\indent
Let $V = V(D) = \{v_1, \ldots , v_n\}$.
Let $In(D)$ be the (diagonal) $(V\times V)$-matrix $(d_{ij})$ 
\\[0.3ex]
such that $d_{ii} = d_{in}(v_i)$ and $d_{ij} = 0$ for $i \ne j$.
Let
$A(D)$ be the $(V\times V)$-matrix $(a_{ij})$ such 
\\[0.3ex]
that 
$a_{ij} = 1$ if $(v_i,v_j) \in E(D)$ and 
$a_{ij} = 0$ if $(v_i,v_j) \not \in E(D)$.
Let 
\\[0.3ex]
$L_{in}(D) = In(D) - A(D)$ and 
$L_{out}(D) = Out(D) - A(D)$.
Clearly, $[L_{in}(D)]^{\top} = L_{out}(D^{-1})$.
\\[0.3ex]
Let 
$L_{in}(\lambda, D) = \det (\lambda I - L_{in}(D))$
 and 
$L_{out}(\lambda, D) = \det (\lambda I - L_{out}(D))$.
We put 
\\[0.3ex]
$L(D) = L_{in}(D)$ and $L(\lambda, D) = L_{in}(\lambda, D)$.
Obviously, $L(D)$ has an eigenvalue $0$.
Let 
\\[0.3ex]
$L(\lambda, D) = \lambda P(\lambda , D)$, and so $P(\lambda , D)$ is a polynomial.
Let $Sp(D)$ denote the set of all 
\\[0.3ex]
$n-1$ roots of $P(\lambda , D)$.
A digraph $D$ is called {\em balanced} if
$d_{in}(v, D) = d_{out}(v, D) \ne 0$ for every $v \in V(D)$.
\bs
{\em \cite{KelLatvia}}
\label{di-reciprocity}
Let  $D$ be a  simple balanced digraph with $n$ vertices.
Then
\\[1ex]
$(a1)$ there exists a bijection
$\sigma :
Sp(D) \to Sp(D^c)$ such that $x + \sigma ( x) = n$
for every $x \in Sp(D)$ or, equivalently,
\\[1ex]
$(a2)$
$P(\lambda , D^c) = (-1)^{n-1} P(n - \lambda , D)$.

Moreover, the matrices $L(D)$ and $L(D^c)$ are simultaneously diagonalizable.
\es

Furthermore, the following Reciprocity Theorem is  true for all simple digraphs.
\bs
{\em \cite{KelLatvia,Ktree1}}
\label{di-reciprocity}
Let  $D$ be a simple  digraph with $n$ vertices.
Then
\\[0.8ex]
$(a1)$
$P(\lambda , D^c) = (-1)^{n-1} P(n - \lambda , D)$
or, equivalently,
\\[0.8ex]
$(a2)$ there exists 
a bijection $\varepsilon : Sp(D) \to Sp(D^c)$ such that $e + \varepsilon ( e) = n$
for every 
\\[0.3ex]
$e \in Sp(D)$.
\es

The  Reciprocity Theorem {\bf \ref{reciprocity}} can also be generalized as follows.
\bs 
{\em \cite{Kcourse}}
\label{general-reciprocity}
Let 
$D = (V, E)$,  $D_1 = (V, E_1)$, and $D_2 = (V, E_2)$
be a simple digraphs such 
\\[0.3ex] 
that 
$E_1 \cup E_2 = E$ and $E_1 \cap E_2 = \emptyset $  
{\em (and so  $D = D_1 \cup D_2$ 
is {\em decomposed in two parts} $D_1$ 
\\[0.3ex]
and $D_2$)}.
Suppose that 
$L(D_1)~ L(D_2) = L(D_2)~ L(D_1)$.
Then there exist  bijections 
\\[0.3ex]
$\alpha _j: Sp(D) \to Sp(D_j)$, $j \in \{1,2\}$, such that
$\alpha _1(e) + \alpha _2(e) = e$ for every $e \in Sp(D)$.
\es

Theorem {\bf \ref{general-reciprocity}} can be further generalized for the case when $G$ is decomposed into 
$p $ parts $D_1, \ldots , D_p$, $p \ge 2$.
\indent
Here is an extension  the  Reciprocity Theorem {\bf \ref{reciprocity}}  to the class of bipartite graphs.
\bs 
{\em \cite{JK}}
\label{bipartite-reciprocity}
Let $X$ and $Y$ be finite disjoint sets, $|X| = |Y| = s$.
Let $B$ be the complete $(X,Y)$-bipartite graph
{\em (and so $v(B) = 2s = n$)}. 
Suppose that $B_1$ and $B_2$ be  $(X,Y)$-bipartite graphs such that 
\\[0.7ex]
$(c1)$ $B_1$ and $B_2$ are edge disjoint and 
$B_1\cup B_2 = B$ and
\\[0.7ex]
$(c2)$ $B_1$ is an $r$-regular graph, and so $B_2$ is  an $(s-r)$-regular graph, and so
\\[0.7ex] 
$V(B_1) = V(B_2) = V(B) = X \cup Y$.
\\[0.7ex]
\indent
Let  
$Sp(B_j) = 
(\lambda _1(B_j) \le \ldots \le \lambda _{n-1}(B_j))$ be the Laplacian spectra of  $B_j$, $j = 1, 2$.
\\[0.7ex]
Then 
\\[0.7ex]
\indent
$(a1)$
$\lambda _{2s-1}(B_1) + \lambda _{2s-1}(B_2) = 2s$
and
\\[0.7ex]
\indent
$\lambda _i(B_1) + \lambda _{2s -1 - i}(B_2) = s$
for every $i \in \{1, \ldots , 2s -2\}$  and, equivalently,
\\[0.7ex]
\indent
$(a2)$
$P(\lambda , B_2) =
(-1)^{2s-1}(\lambda - 2s + \lambda _{2s-1})
((\lambda - s + \lambda _{2s-1})^{-1}
P(s - \lambda , B_1)$.
\es

For example, let $B = K_{3,3}$, $B_1$ be a 6-vertex cycle  
in $B$, $B_2 = B \setminus E(B_1)$, and so  $B_2$ is 
a 3-edge matching,
$B = B_1 \cup B_2$, $B_1$ is 2-regular, $B_2$ is 1-regular, and $s = 3$. 
Then 
\\[0.7ex]
\indent
$Sp(B_1) =  (1 ,1 , 3 , 3 , 4)$ and  
$Sp(B_2) = \{0,  0 , 2 , 2 , 2\}$. 
\\[0.7ex]
Therefore
\\[0.7ex]
\indent
$\lambda _5(B_1) + \lambda _5(B_2) = 4 + 2 = 6 = 2s$ and
\\[0.7ex]
\indent
$(\lambda _1(B_1), \lambda _2(B_1),\lambda _3(B_1),  (\lambda _4(B_1)) +
(\lambda _4(B_2), \lambda _3(B_2), \lambda _2(B_2), (\lambda _1(B_2)) =
\\[0.9ex]
\indent
(1 ,1 , 3 , 3) + (2 ,2, 0, 0) = (3, 3, 3, 3) = (s, s, s, s)$. 
\\[1.5ex]
\indent
From {\bf \ref{reciprocity}} we have:
\bs 
{\em \cite{Kcourse,Ktree1}}
\label{AtimesB}
Let $A$, $B$, and $G$ be simple graphs,  $v(G) = n$,
$v(A) = a$, and $v(B) = b$.
Then
\\[1ex]
$(a1)$
$0 \le \lambda _1(G) \le \ldots  \le \lambda _{n-1}(G) \le n$,
\\[1ex]
$(a2)$ $cmp(G) = \mu (0)$ and $cmp(G^c) = \mu (n) +1$,
where $\mu (z)$ is the multiplicity of the eigenvalue 
$z$ of $L(G)$,
\\[1ex]
$(a3)$ 
$P(\lambda , A + B) = 
\lambda P(\lambda , A) P(\lambda , B)$, i.e.
$\Phi (\lambda , A + B) = \Phi (\lambda , A) \Phi (\lambda , B)$, and
\\[1ex]
$(a4)$
$P(\lambda , A \times B) = 
(\lambda - a - b) P(\lambda - b, A) P(\lambda - a, B)$.
\es

Since $nt(G) = c_{n-1}(G) = (-1)^{n-1}P(0, G)$, we have from 
{\bf \ref {reciprocity}} $(a2)$:
\bs {\em \cite{Kcourse,Ktree1}}
\label{P(n+r)=}
Let $G, F \in {\cal G}_n$.
Then
\\[1ex]
$(a1)$
$t(K_{n+r} - E(G)) = (n + r) ^{r-1} P(n+ r, G)$, and so
\\[1ex]
$(a2)$
$t(K_n - E(G)) = n ^{-1} P(n, G) =
 \sum \{(-1)^s c_s(G) n ^{n- 2 - s}: 
s \in \{0, \ldots , n-1\}\}$.
\es

\bs {\em \cite{Kcourse,KCh}}
\label{incl-excl}
Equality $(a2)$ in {\bf \ref{P(n+r)=}} is the inclusion-exclusion formula for the number $t(K_n - E(G))$ of spanning trees of $K_n$ avoiding the edges of its subgraph $G$.
\es

{\bf Proof.}~~
Let $F$ be a forest in $K_n$ and $t(K_n, F)$ denote the number of spanning trees in $K_n$ containing $F$.
Let 
\\[1ex]
\indent
$\Sigma _s(K_n)  = \sum \{t(K_n, F): F \in {\cal F} (K_n), e(F) = s\} $. 
\\[1ex]
As we have shown in   \cite{Ktree-forest}, 
$t(K_n, F^s) = \gamma (F) n ^{n- 2 - s}$. Therefore   
$\Sigma _s(K_n) = c_s(G) n ^{n- 2 - s}$.
By  {\bf \ref{P(n+r)=}} $(a2)$,
\\[1ex]
\indent
$t(K_n - E(G)) = \Sigma _0 - \Sigma _1 + \cdots + (-1)^s \Sigma _s + \cdots + (-1)^m \Sigma _m$,
where $m = e(G)$.
\ep
\\[2ex]
\indent
Let $G$ be a simple graph.
Let $G^l$ denote the {\em line graph of} $G$, i.e. 
$V(G^l) = E(G)$ and $(a,b) \in E(G^l)$ if and only if
edges $a$ and $b$ in $G$ are adjacent. 
Let $\dot{G}$ denote the graph obtained from $G$ by subdividing every edge $e$ of $G$ by exactly one vertex $s(e)$.
Let $\check{G}$ denote the graph obtained from $\dot{G}$ by adding edge $s(a)s(b)$ if and only if edges $a$ and $b$ in $G$ are adjacent. 
\bs {\em \cite{Kpolinom}}
\label{`line'graphs} 
Let $G$ be an $r$-regular graph with $n$ vertices and 
$mr$ edges $($and so $m = \frac{1}{2}nr $$)$. Then
\\[1ex]
$(a1)$
$P(\lambda , G^l) = (\lambda - 2r)^{m - n} P(\lambda , G)$, 
\\[1ex]
$(a2)$
$P(\lambda , \dot{G}) =
(-1)^n(\lambda - 2)^{m - n}
P(\lambda (r + 2 - \lambda)), G)$,
and
\\[1ex]
$(a3)$ 
$P(\lambda , \check{G}) =
(\lambda -  r - 1)^n(\lambda - 2r - 2)^{m - n}
P(\frac{\lambda ^2 -  (r + 2)\lambda}{\lambda- r - 1}, G)$.
\es

From 
{\bf \ref{somecoefficients}} $(a4)$ and {\bf \ref{`line'graphs}} 
we have:
\bs {\em \cite{Kpolinom}}
\label{`line'graphs-trees} 
Let $G$ be an $r$-regular graph with $n$ vertices and 
$mr$ edges (and so $m = \frac{1}{2}nr $). Then
\[t (\check{G}) =
\frac{n}{m+n} 2^{m - n}(r + 1)^{m-1}(r+2)   t(G),\]
\[ t(\dot{G}) = \frac{n}{m+n}2^{m - n} (r+2)  t( G),\]
and
\[ t(G^l) = \frac{n}{m} 2^{m - n} r^ {m - n}  t(G).\]
\es

We will see below  that threshold graphs (see definition {\bf \ref{thresholdgraph}}) play a special role in problems  ${\cal T}_{min}$ and ${\cal R}_{min}$ as well as some other optimization problems. It is known that a threshold graph is uniquely defined by its degree sequence. 

For a connected threshold graph $G$ with $n$ vertices, there exists a 
partition of $V(G) = S \cup K$ into two disjoint sets  
$S$ and $K$ (with $|S| = s$, $|K| = k$, and so 
$s + k = n$) such that $K \ne \emptyset $, 
$G[S]$ has no edges, $G[K]$ is a maximum complete subgraph of $G$, and
there is an ordering $K = (x_1,  \ldots , x_k)$ of the vertices in $K$ and 
an ordering
$S = (x_{k+1}, \ldots , x_n)$ of vertices  in $S$ such that 
$i > j \Rightarrow N(x_i,A) \subseteq N(x_j,A)$. Let $d_i = d(x_i,G)$. 
Then $d_{1} \ge \ldots \ge d_n$.
\\[1ex]
\indent
The recursive description of a threshold graph $G$ allows 
to give the following explicit formulas for the Laplacian spectrum $Sp(G)$ 
and the  number of spanning trees $t(G)$ in terms of the degree sequence of $G$.
\bs {\em 
\cite{HK,Kcourse}}
\label{t(G),thresholdG}
Let $G$ be a connected threshold graph.
Then 
\\[1.5ex]
\indent
$(a1)$
$Sp(G) = (d_1 + 1\ge \ldots \ge d_{k-1} +1\ge d_{k+1}\ge \ldots \ge d_n)$, where $d_1 + 1 = n$,
\\[1.5ex]
\indent
$(a2)$ $t(G) = n ^{-1} \prod \{d_i + 1: i = 1, \ldots k-1\}
\times \prod \{d_j : i = k+1, \ldots ,n\}$, and
\\[1.7ex]
\indent
$(a3)$
$t(G) = d_n (\frac{n}{n -1})^{d_n - 1} t(G - x_n) = 
n^{-1} k^{k-1} \prod 
\{d_{k+i} (\frac{k + i}{k+i-1})^{d_{k + i}}: i = 1, \ldots, s\}$.
\es
\vskip 0.6ex

~~~It is also easy to prove the following.
\bs {\em \cite{HK,Kcourse}}
\label{threshold-G-defined-by-Lspectrum}
Every threshold graph is uniquely defined by its Laplacian spectrum as well as by its degree function.
\es

\bs {\em \cite{KuniGr}, see also \cite{Kpolinom}}
\label{sK_n-defined-by-Lspectrum}
Let $G, F \in {\cal G}_n^m$.
Suppose that  
\\[0.7ex]
$(h1)$
each component of $G$ is a complete graph
and 
\\[0.7ex]
$(h2)$
$isl(G) \le isl(F)$.
\\[0.7ex]
\indent
Then $G$ and $F$ are isomorphic if and only if 
$P(\lambda , G) = P(\lambda , F)$.
\es

The following example shows that condition $(h2)$ in 
{\bf \ref{sK_n-defined-by-Lspectrum}} is essential.
Let $G = K_{6} + K_{10}$ and $F$ a graph obtained from the line graph of $K_{6}$ by adding an isolated vertex.
Then 
$P(\lambda , G) = P(\lambda , F) = \lambda (\lambda - 6)(\lambda - 10)$ and, obviously, $G$ and $F$ are not isomorphic.  

Notice that  the description of our theorem  
{\bf \ref{sK_n-defined-by-Lspectrum}} in \cite{CDS}, page 163, is incorrect, namely, condition $(h2)$ is missing.
\\[1.5ex]
\indent
It turns out that among totally decomposable graphs there are infinitely many non-isomorphic $L$-cospectral graphs.
\bs {\em \cite{KuniGr} (see also \cite{Kpolinom})}
\label{non-isom,cospectr,decomposG}
For every $n \ge 16$ there exist totally decomposable graphs $G$ and $F$ with $n$ vertices such that
$G$ and $F$ are not isomorphic and 
$P(\lambda ,G) = P(\lambda , F)$.
\es

{\bf Proof.} 
Let $A = g^{s-1} + (2g)g^{s-1}$ and $B = 2g^s$, where $s \ge 2$. 
Let $G = A \times B^c$ and $F = A^c \times B$.
Then 
$G$ and $F$ are totally decomposable,
$v(G) = v(F) = 8s$,
$G$ and $F$ are not isomorphic, and
$P(\lambda ,G) = P(\lambda , F)$.
\ep
\bs {\em \cite{Kcourse,Ktree1}}
\label{LaplacIntegerG}
Suppose that a graph $G$ is obtained from 
disjoint simple graphs $F_1,  \ldots , F_k$
 by a series of operations $+$ and $\times$.
If every $F_i$ has an integer Laplacian spectrum, then $G$ has also an  integer Laplacian spectrum. In particular, every totally decomposable graph has an integer spectrum.
\es

In \cite{KeigvMtr} we proved a more general result of this type for weighted graphs.

\bs {\em \cite{Kcourse}}
\label{non-isom,cospectr,decomposG}
Let $G$ be a simple graph having $n$ vertices
and the Laplacian spectrum $\{1, \ldots, n-1\}$.
Then $n \ge 16$.
\es

Our proof of theorem 
{\bf \ref{non-isom,cospectr,decomposG}} does not use any results obtained by computer except for the table in \cite{RW} of some special graphs on 7 vertices.
\\[1ex]
\indent
Let $T$ be a tree and $d(x,y,T) = d(x,y)$ be the number of edges in the path of $T$ with the ends $x$ and $y$.
Let 
\\[1ex]
\indent
$W(T) = \sum \{d(x,y): \{x,y\} \subseteq V(T), x \ne y\}$.
\\[1ex]
The parameter $W(T)$ is called the {\em Weiner index of a tree} $T$. The Weiner index turns out to be a useful notion in organic chemistry \cite{W}.
\\[1ex]
\indent
Let $T(x,y)$ be the graph obtained from $T$ by identifying vertices $x$ and $y$, $x \ne y$. Then obviously,
$d(x,y,T) = t(T(x,y))$. Therefore from 
{\bf \ref{posemdef}} and
 {\bf \ref{coefficients}}
we have:
\bs {\em \cite{Kproceed,Kcourse,W}} 
\label{Weiner}
Let $T$ be a tree with $n$ vertices.
Then 
\\[1ex]
$W(T) = c_{n-2}(T) = \sum\{t(T_ {V - S}) : S \subseteq V, |S| = n - 2\} =
\\[1ex]
\sum\{f(T, V - S) : S \subseteq V, |S| = n - 2\} = 
\sum\{\gamma (F): F \in {\cal F}(T), e(F) = n - 2\} =
\\[1ex] 
\sigma _{n-2} (\lambda _1(G), \ldots , \lambda _{n-1}(T))$.
\es

\section{Some operations on graphs}

\label{operations}

\indent

Let $\succeq $ be a partial order relation on a subset  
${\cal A}$ of ${\cal G}_n^m$ and $Q: {\cal A} \to {\cal A}$
a function on  ${\cal A}$. 
We say that $Q$ is an {\em operation on} ${\cal A}$ and 
that {\em $Q$ is $\succeq $-increasing on ${\cal A}$ 
$($$\succeq $-decreasing on ${\cal A}$$)$} 
if $Q(G) \succeq G$ 
$($respectively, $ G \succeq Q(G)$$)$ for every 
$G \in {\cal A}$.
We say that {\em $Q$ is a $\succeq $-monotone 
operation on ${\cal A}$} if either $Q$ is $\succeq $-increasing or $\succeq $-decreasing on ${\cal A}$.

A function $f: {\cal A} \to \mathbb{R}$ induces the following quasy-linear order 
$\succeq ^f$: for $G, G' \in {\cal A}$, let $G\succeq ^f F$ if 
$f(G) \ge f(G')$. In this case instead of {\em $\succeq ^f$-increasing}, {\em $\succeq ^f$-monotone}, etc. we say simply
{\em $f$-increasing}, {\em $f$-monotone}, etc.

A set ${\cal Q}$ of operations on ${\cal A}$ induces a partial order 
relation as follows:
given $G, F \in {\cal A}$, we define $G \succeq ^{\cal Q} F$ if
$F$ can be obtained from $G$ by a series of operations from ${\cal Q}$. 
Obviously, a graph $A$ in ${\cal A}$ is $\succeq ^{\cal Q}$-minimal if and only if no operation from ${\cal Q}$ can be applied to $A$.
\\[1ex]
\indent
We will use the following simple observation.
\bs
\label{observation}
If  $Q$ is $\succeq $-increasing $($$\succeq $-decreasing$)$ operation on ${\cal A}$, then for every $G \in {\cal A}$ there exists $A \in {\cal A}$ such that 
$A \succeq G$ $($respectively, $G \succeq A$$)$ and operation $Q$ cannot be applied to $A$ $($i.e. $A$ is 
 $\succeq ^Q$-minimal$)$.
\es

\begin{figure}
\begin{center}
\includegraphics[scale=0.3]{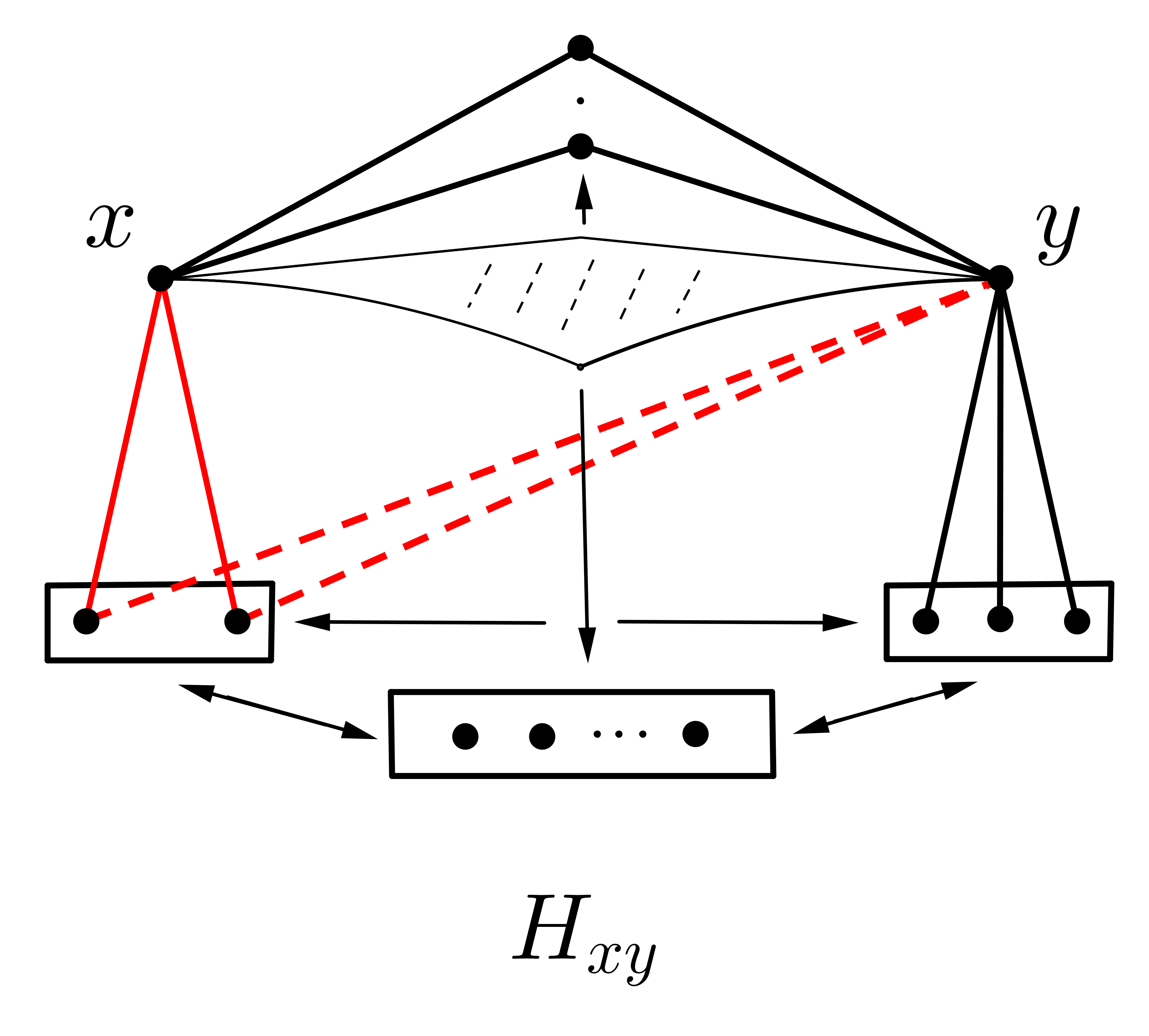}
\end{center}
\caption{$H_{xy}$-operation }
\label{hammock}
\end{figure}

In \cite{KgraphTrnsf,Kmxtr,KoperProb,KoperTrees,Kcmpr} we introduced various operations on graphs that preserve the number of vertices and edges of a graph and that are  
monotone with respect to some graph parameters.
Here are some of these operations \cite{KoperTrees,KoperProb}.
\\[1ex]
\indent
Let $H$ be a graph, $x,y \in V(H)$, and $x \ne y$.
We call $xHy$ a {\em two-pole with poles $x$ and $y$}.
Let $uAv$ be another two-pole.
Let $G$ be obtained from disjoint $xHy$ and $uAv$ by
identifying $x$ with $u$  and  $y$ with $v$.
Let $H_{xy}(G) = (G - [x,X]) \cup [y,X]$, where 
$X = N_x(A) \setminus (N_y(A) \cup y)$ and 
$Y = N_y(A) \setminus (N_x(A) \cup x)$.
We call the two-pole $xHy$ an $(x,y)$-{\em hammock in} $G$ and  call this operation the $H_{xy}$-{\em operation} or just a {\em hammock operation} (see 
Figure \ref{hammock}).

A two-pole $xHy$ is called {\em symmetric} if 
$H$ has an automorphism $\alpha : V(G) \to V(G)$ such that $\alpha (x) = y$ and  $\alpha (y) = x$.
We call the $H_{xy}$-operation {\em symmetric} if
$xHy$ is a symmetric two-pole.
Obviously, if $u$ or $v$ is an isolated vertex in $A$, then
$H_{xy}(G)$ is isomorphic to $G$. Therefore, when applying this operation,  we will always assume that both 
$u$ and $v$ are not isolated vertices in $A$.

We call $H_{xy}$-operation an {\em $(x,y)$-path operation} or {\em $xPy$-operation} if $xHy$ is an $(x,y$)-path. 
Obviously, an $(x,y)$-path operation is symmetric.

If $d(x,G) = d(x,H) +1$, then we call the $H_{xy}$-operation 
a {\em close-door $H_{xy}$-operation}. 
The reverse of 
a close-door $H_{xy}$-operation is called an  
{\em open-door $H_{xy}$-operation}.

A hammock operation is a particular case of more general operations in \cite{KoperTrees,KoperProb} 
(see also {\bf \ref{weightedDiamondTree}} and 
{\bf \ref{weightedDiamondRandom}} below).

\begin{figure}
\begin{center}
\includegraphics[scale=0.22]{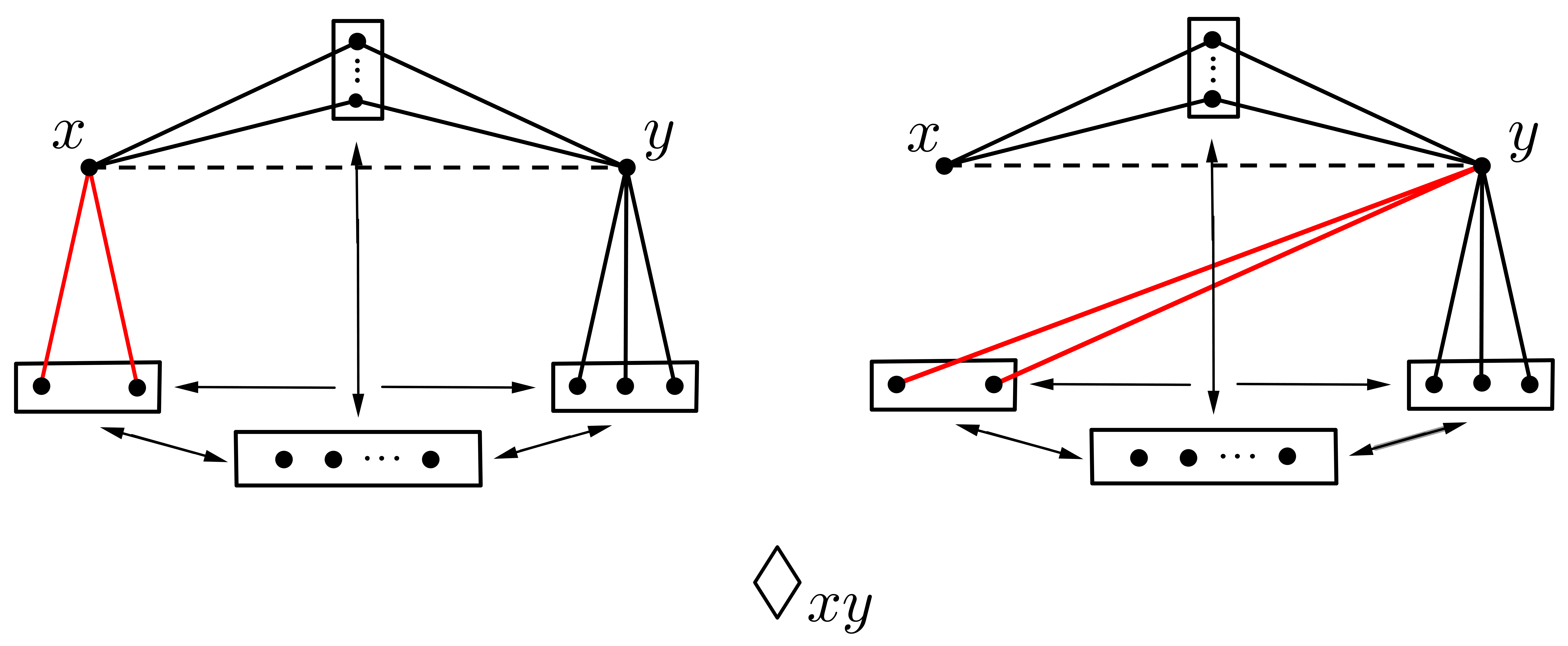}
\end{center}
\caption{$\diamondsuit _{xy}$-operation }
\label{diamond}
\end{figure}

One of possible  specifications of the $H_{xy}$-operation is
when $V(H) = \{x,y\}$ (see Figure \ref{diamond}).

We call this simpler operation the {\em $\diamondsuit _{xy}$-operation} 
\cite{KoperTrees,KoperProb}.
Obviously, the {\em $\diamondsuit _{xy}$-operation} is symmetric. 

\begin{figure}
\begin{center}
 \includegraphics[scale=0.3]{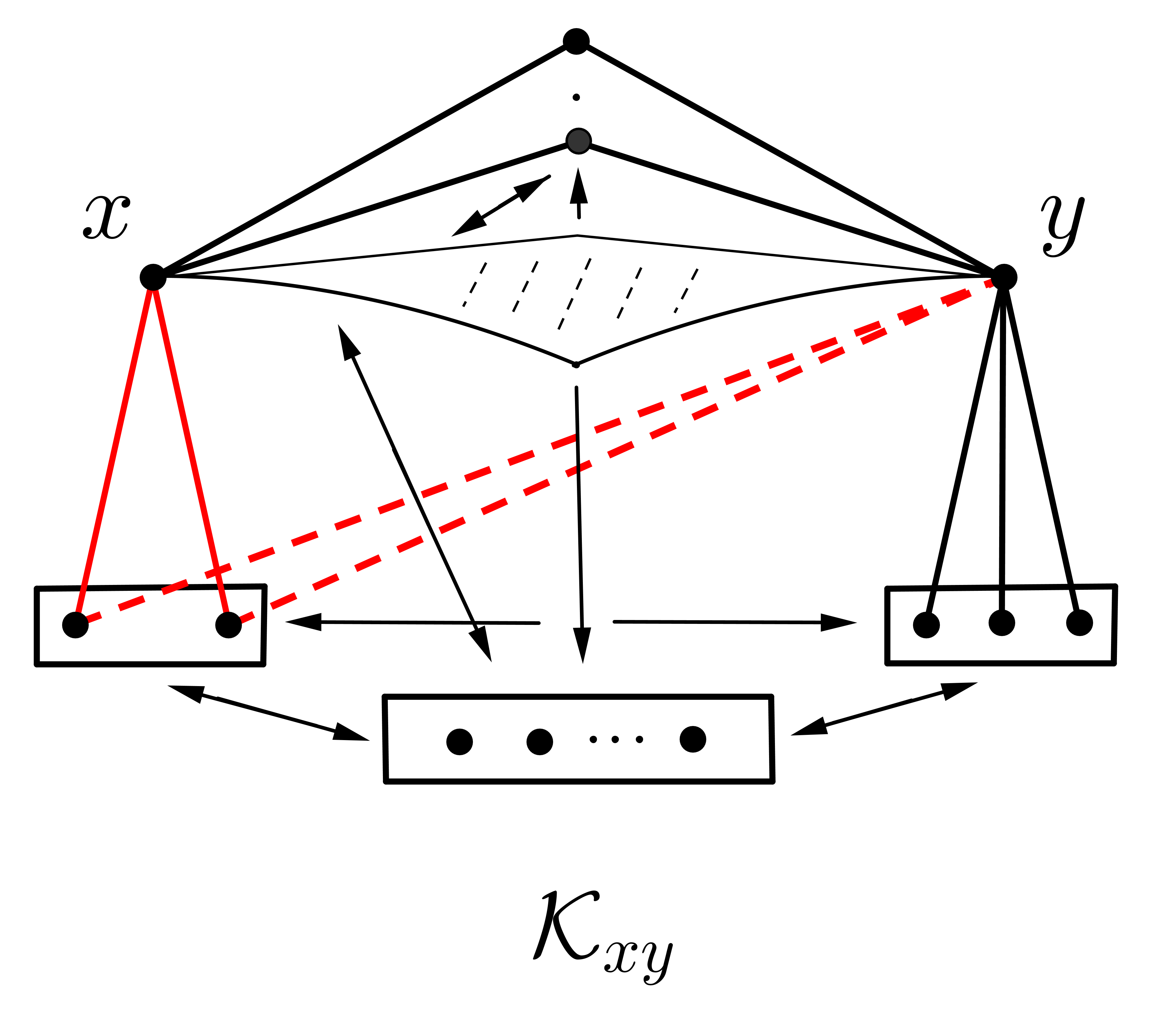}
 \end{center}
 \caption{${\cal K}_{xy}$-operation}
 \label{K-operation}
\end{figure}

We will also use the following  particular case of 
the close-door $\diamondsuit _{xy}$-operation 
\cite{KoperTrees,KoperProb}.
Let $G$ be a graph with three vertices $x$, $y$, $z$  such that $xz \in E(G)$ and $yz \not \in E(G)$.
Let $G' = G - xz + yz$. 
We put $G' = D_{xyz}(G)$  if $xz$ is the only edge in $G$ incident to $x$ and say that $G'$ {\em is obtained from 
$G$ by the $D_{xyz}$-operation}.
We also call a $D_{xyz}$-operation a {\em close-door operation} (it `closes' the `door' $xz$). Accordingly, the reverse of 
a close-door operation is called an  {\em open-door operation}.
\\[1ex]
\indent
The following is a natural generalization of the $H_{xy}$-operation \cite{Kcourse,Kproceed}.
Let $G$ be a graph,  $x,y \in V(G)$, and $x \ne y$.
Let ${\cal K}$ be an induced  subgraph of $G$ containing $x$ and  $y$, and so $x{\cal K}y$ is a two-pole.
Let 
$X = N_x(G) - (V({\cal K})\cup N_y(G)) $ and 
$Y = N_y(G) - (V({\cal K})\cup N_x(G))$.
Let $[x,X] = \{xv: v \in X\}$ and $[y,X] = \{yv: v \in X\}$, 
and so $[x,X] \subseteq E(G)$ and $[y,X] \cap  E(G) = \emptyset $.
Let 
$G' ={\cal K}_{xy}(G) =  (G - [x,X]) \cup [y,X]$, and 
so $[y,X] \subseteq E(G')$ and 
$[x,X] \cap  E(G) = \emptyset $.
We say that 
{\em ${\cal K}_{xy}(G)$  is obtained from $G$ 
by the ${\cal K}_{xy}$-operation} (see Figure 
\ref{K-operation}).

Obviously, if $X$ or $Y$ is empty, then 
${\cal K}_{xy}(G) = G$. Therefore, when applying this operation,  we will always assume that both $X$ and $Y$ are not empty sets.
\\[1ex]
\indent
It turns out that under certain conditions
on $(G, {\cal K},x,y)$
some graph parameters are 
`monotone' with respect to this operation.
\\[1ex]
\indent
A  ${\cal K}_{xy}$-operation on $G$ is called {\em $\alpha $-symmetric} or, simply, {\em symmetric} if
$G -([x,X]) \cup [y,Y])$ has an automorphism 
$\alpha : V(G) \to V(G)$ such that
$\alpha (x) = y$ and 
$\alpha (y) = x$,  $\alpha (z)  \in V({\cal K})$ for every 
$z \in V({\cal K})$ and $\alpha (v) = v$ for every 
$v \in X \cup Y$, and so
$\alpha (u) \in V(G - {\cal K})$ for every
$u \in V(G - {\cal K})$.
\\[1ex]
\indent
Clearly, $H_{xy}(G)$ is obtained from $G$ by the  
${\cal K}_{xy}$-operation, where  $x{\cal K}y = xHy$.
\\[1ex]
\indent
The above described operations proved to be very useful for ``improving'' some graph characteristic and finding graphs with some extremal properties.
It turns out that many results  on the $H_{xy}$-operation are also true for the ${\cal K}_{xy}$-operation.
\\[1ex]
\indent
Here are some useful properties of 
the ${\cal K}_{xy}$-operation.
\bs 
\label{properties-of-Kxy} 
Let $(G,{\cal K}, x, y)$ be as described above and let $G_{xy} = {\cal K}_{xy}(G)$.
Then 
\\[0.7ex]
$(a1)$
$V(G_{xy}) = V(G_{yx})$, $E(G) = E(G_{xy})$,
and there exists an isomorphism $\omega $ from $G_{xy}$ to $G_{yx}$ such that
$\omega (x) = y$ and $\omega (y) = x$,
\\[0.7ex]
$(a2)$ for a simple graph $G$, $[{\cal K}_{xy}(G)]^c =  ({\cal K}^c_{yx})([G]^c)$ and
if the ${\cal K}_{xy}$-operation on $G$ is $\alpha $-symmitric, then the $[{\cal K}]^c_{xy}$-operation on $[G]^c$ is also 
$\alpha $-symmitric;
in particular, 
$[{\diamondsuit }_{xy}(G)]^c = {\diamondsuit }_{yx}([G]^c)$,
 and
\\[0.7ex]
$(a3)$ 
If ${\cal K}^*$ is the subgraph of $G$ induced by 
$V({\cal K}) \cup (N_x(G) \cap N_y(G)$, then
${\cal K}_{xy}(G) = {\cal K}^*_{xy}(G)$ and
the ${\cal K}_{xy}$-operation on $G$ is symmetric if and only if the ${\cal K}^*_{xy}$-operation on $G$ is symmetric.
 \es

Obviously, the $\diamondsuit $-operation cannot be applied to a graph $G$ 
(i.e. $G$ is $\succeq ^{\diamondsuit }$-minimal) if and 
only if for every two distinct vertices $x$, $y$ of $G$
either $N(x,G) - \{x,y\} \subseteq N(y,G)$ 
or $N(y,G) - \{x,y\} \subseteq N(x,G)$, i.e. if and only if $G$ is 
a   vertex comparable graph (or, the same, a threshold graph). Therefore we have:
\bs {\em \cite{Kcourse,Kproceed}}
\label{GtoTreshGraph}
If $G$ is a  non-threshold graph, then there exists a threshold graph $F$  obtained from $G$ by a series of $\diamondsuit $-operations.
\es

It is easy to prove the following strengthening of 
{\bf \ref{GtoTreshGraph}}.
\bs {\em \cite{Kcourse,Kproceed}}
\label{conGtoConTreshGraph}
If $G$ is a connected non-threshold graph, then there exists a connected threshold graph $F$  obtained from $G$ by 
a series of $\diamondsuit $-operations.
\es

Thus, from  {\bf \ref{observation}} and 
{\bf \ref{conGtoConTreshGraph}}
we have:
\bs {\em \cite{Kcourse,Kproceed}}
\label{GtoTreshGraph}
Suppose that the $\diamondsuit $-operation is
$\succeq $-decreasing. 
Then for every $G$ in ${\cal G}_n^m$ 
$($in ${\cal C}_n^m$$)$ and every clique $K$ in $G$ there exists  $F$ in ${\cal G}_n^m$ 
$($respectively, in ${\cal C}_n^m$$)$ such that  $F$ is a  
threshold graph in ${\cal G}_n^m$ 
$($respectively, in ${\cal C}_n^m$$)$, $F$ contains $K$,
and $G \succeq F$.
\es

\bs {\em \cite{Kcourse,Kproceed,KoperTrees}}
\label{T>G}
Suppose that the close-door operation is
$\succeq $-decreasing.
Let $G \in {\cal C}_n^m$ and let $\ddot{G}$ be the graph obtained from $G$ by adding $m -n +1$ isolated vertices. Then for every spanning tree $T$ of $G$ there exists a tree $D$ with $m$ edges such that $T$ is a subgraph of $D$ and $D \succeq  \ddot{G}$.
\es

Theorem {\bf \ref{T>G}} follows from the fact that 
there exists a series of open-door operations that transforms $\ddot{G}$ to a tree $D$ with $e(G)$ edges and with $T \subseteq D$.

\bs {\em \cite{Kcourse,Kproceed,KoperTrees}}
\label{K-F}
Suppose that the $\diamondsuit $-operation is
$\succeq $-decreasing. Let $F$ be a simple graph with no isolated vertices, with $r$ edges, and with at most $n$ vertices. 
Then 
\\[1ex]
$(a1)$
$K_n - E(rP^1) \succeq K_n - E(P^2 + (r-2)rP^1)   \succeq K_n - E(F)$ 
for every $F$  not isomorphic 
\\[1ex]
to $rP^1$ or $P^2 + (r-2)P^1$ and for $r \le n/2$ and 
\\[1ex]
$(a2)$ $K_n - E(F)  \succeq K_n - E(S^r)$ for every $F$  not isomorphic to $S^r$ and for $r \le n-1$.
\es

The similar claims are true for  a $\succeq $-increasing
$\diamondsuit $-operation.

\section{Definitions of some   posets of graphs}
\label{partialOrders}

\indent

We will first describe  some Laplacian posets, namely some  partial order relations on 
${\cal G}^m$ and ${\cal G}_n^m$ related with the Laplacian polynomials of graphs 
\cite{KgraphTrnsf,Kmxtr,Kcmpr,KCh}.
Let $\nu (G,F) = \max \{v(G), v(F)\}$.
We remind that
\\[1.5ex]
\indent
$P(\lambda , G) = \sum \{(-1)^s c_s(G) \lambda ^{n- 1 - s}: 
s \in \{0, \ldots , n-1\}\}$ and 
\\[1.5ex]
\indent
$\Phi (\lambda , G) = \lambda ^{m - n +1} P(\lambda , G)$.
where $n = v(G)$.
\\[1.5ex]
\indent
Let $s$ and $r$ be non-negative integers and $x$ be a real number. We write:
\\[1.7ex]
$(\tau )$ $G \succeq ^{\tau } _s F$ if 
$t(K_{s+r} - E(G)) \ge t(K_{s+r} - E(F))$ for every integer 
$r \ge 0$ and  $s \ge \nu(G, F)\}$,
\\[0.7ex]
\indent
$G \succ ^{\tau }_s F$ if $G \succeq ^{\tau }_s F$ and 
$t(K_{s+r} - E(G)) > t(K_{s+r} - E(F))$ for some $r \ge 0$, and
\\[0.7ex]
\indent
$G \succ\succ  ^{\tau }_s F$ 
if $G \succeq ^{\tau }_s F$ and 
$t(K_{s+r} - E(G)) > t(K_{s+r} - E(F))$ for every $r \ge 0$.
\\[1.7ex]
$(p)$ $G \succeq ^p_x F$ if $ v(G) = v(F)$ and
$P (\lambda, G) \ge P (\lambda, F)$ for every
$\lambda \ge x$,
\\[0.7ex]
\indent
$G \succ ^p_x F$ if $G \succeq ^p_x F$ 
and $P (\lambda, G) > P (\lambda, F)$ for some 
$\lambda  \ge x$, and 
\\[0.7ex]
\indent
$G \succ\succ  ^p_x F$ if $G \succeq ^p_x F$ and 
$P (\lambda, G) > P(\lambda, F)$ for every
$\lambda  \ge x$,
\\[1.7ex]
$(\phi)$ 
$G \succeq ^{\phi}_x F$ if  
$\Phi (\lambda ,G) \ge \Phi (\lambda ,F)$ for every
$\lambda \ge x$,
\\[0.7ex]
\indent
$G \succ ^{\phi }_x F$ if $G \succeq ^{\phi }_x F$ and
$\Phi(\lambda, G) > \Phi (\lambda, F)$ for some 
$\lambda  \ge x$,
and
\\[0.7ex]
\indent
$G \succ\succ  ^{\phi}_x F$ if $G \succeq ^{\phi}_x F$ and 
$\Phi(\lambda, G) > \Phi(\lambda, F)$ for every 
$\lambda  \ge x$,
\\[1.7ex] 
$(c)$ $G \succeq ^c F$ if $v(G) = v(F) = n$, $c_s(G) \ge c_s(F)$ for every $s \in \{2, \ldots , n-1)\}$ and 
\\[0.7ex]
\indent
$G \succ ^c F$ if $G \succeq ^c F$ and 
$c_s(G) > c_s(F)$ for some $s \in \{2, \ldots , n-1\}$,
\\[1.7ex]
$( t )$
$G \succeq ^{t } F$ if $t (G) \ge 
t (F)$ and
$G \succ  ^{t } F$ if 
$t(G) > t(F)$,
\\[1.7ex]
$( \lambda )$
$G \succeq ^{\lambda } F$ if $\lambda (G) \ge 
\lambda (F)$ and
$G \succ  ^{\lambda } F$ if 
$\lambda (G) > \lambda (F)$,
\\[1.7ex]
$( \infty )$
$G \succeq ^{\infty} F$ if 
there exists a number $N$ such that
$P(\lambda ,G) \ge P(\lambda ,F)$ for $\lambda \ge N$. 
\\[1.7ex]
\indent
If $v(G) = v(F) = n$, we write 
$\succeq ^p$ instead of $\succeq ^p_n$ and 
$\succeq ^{\tau }$ instead of $\succeq ^{\tau }_n$.
If $G \succeq ^{\phi}_x F$ and $x =  \max \{\lambda (G), \lambda (F)\}$, we write 
$ \succeq ^{\phi }$ instead of  $\succeq ^{\phi }_x$, and so in this case
$\lambda (G) \le \lambda (F)$.
Notice that if $v(G) = v(F)$, then relations 
$\succeq ^p_x$ and $\succeq ^{\phi }_x$ are the same and 
~$ \succeq ^{\phi } ~\Rightarrow ~\succeq ^p$.
\\[1.7ex]
 \indent
Now we will describe some other posets  on graphs.  
Let $R^k(p,G)$ denote the probability that the random graph $(G,p)$ has at most $k$ components. 
Let 
$a^k_s(G)$ is the number of spanning subgraphs of $G$ with $s$ edges and at most $k$ components (and so 
$f_k(G) = a^k_{n-k}(G)$ is the number of spanning forests of $G$ with exactly $k$ components).
Then obviously, 
\\[1.5ex]
\indent
$R^k(p,G) = 
\sum \{a^k_s(G) p ^s q^{m-s}: s \in \{n-k, \ldots , m\}\}$. 
\\[1.5ex]
Let $A(\lambda , G)$ denote the characteristic polynomial of the adjacency matrix $A(G)$ of a graph $G$.
Let $h_0(G)$ and $h_1(G)$ denote the number of Hamiltonian cycles and 
Hamiltonian paths in $G$, respectively.
\\[1ex]
\indent
Let $G, F \in {\cal G}_n^m$. We write:
\\[1.7ex]
$(a)$ $G \succeq ^{a(k)} F$ if $a^k_s(G) \ge a^k_s(F)$ for every 
$s \in \{n-k, \ldots , m\}$,
\\[1.7ex]
$(r)$ $G \succeq ^{r(k)} F$ if $R^k(p,G) \ge R^k(p,F)$ for every $p \in [0,1]$, 
\\[1ex]
\indent
$G  \succ\succ  ^{r(k)} F$ if $R^k(p,G) > R^k(p,F)$ for every $p \in (0,1)$, 
\\[1.7ex]
$(\alpha)$ 
$G \succeq ^{\alpha } F$ if  $\alpha(G) \le 
\alpha (F)$ and
$A (\lambda ,G) \ge A(\lambda ,F)$ for 
$\lambda \ge \alpha (G)$,
\\[1ex]
\indent
$G \succ ^{\alpha } F$ if $G \succeq ^{\alpha } F$ and
$A(\lambda, G) > A (\lambda, F)$ for some 
$\lambda  \ge \alpha (F)$, and
\\[1ex]
\indent
$G \succ\succ  ^{\alpha } F$ if $G \succeq ^{\alpha } F$ and 
$A(\lambda, G) > A(\lambda, F)$ for all 
$\lambda  \ge \alpha (F)$,
and
\\[1ex]
$(h)$ $G \succeq ^{h_i} F$ if $h_i(G) \ge h_i(F)$ for 
$i \in \{1,2\}$.
\\[1.7ex]
\indent
Put $a_s(G) =  a^1_s(G)$, $R(p,G) = R^1(p,G)$,  and 
let  $\succeq ^r $ be $\succeq ^{r(1)}$.
Notice that  
\\[1ex]
$a_{n-1}(G) = t(G) = n^{-1}c_{n-1}(G)$, where $(-1)^{n-1}c_{n-1}(G)$ is the last coefficient of $P(\lambda ,G)$. 
\\[1.7ex]
\indent
For $G, F \in {\cal G}_n^m$ and $z \in \{c, p, \lambda , t, a, r, \infty \}$, we write $G \asymp ^z F$ if
$G \asymp ^z F$ if $G \succeq ^z F$ and $F \succeq ^z G$.
We also write $G \succeq ^{\alpha,\beta} F$ instead of
$G \succeq ^\alpha F$ and $G \succeq ^\beta F$.
\\[1.7ex]
\indent
In order to define relation $\succ\succ  ^c$ on ${\cal G}_n^m$ we need the following simple observation.
For $G, F \in {\cal G}_n^m$,  
$G$ and $F$ are called {\em  $Cmp$-cospectral}  if both 
$G$ and $F$ are forests and  there exists a bijection  
$\sigma : Cmp(G) \to Cmp (F)$ such that $v(C) = v(\sigma (C))$ for every $C \in Cmp(G)$.
It is easy to see that if $G$ and $F$ are $Cmp$-cospectral, then  
\\[1ex]
$c_m(G) = c_m(F) = \gamma (F)$.
\\[1ex]
\indent
 Now in view of {\bf \ref{somecoefficients}} 
  we define $\succ\succ  ^c$ on ${\cal G}_n^m$ as follows:
$G \succ\succ  ^c F$  if  $c_s(G) > c_s(F)$ for every 
$s \in \{2, \ldots , m-1\}$ in case  $G$ and $F$ are $Cmp$-cospectral, and $c_s(G) > c_s(F)$ for every $s \in \{2, \ldots , n - 1\}$ in case  $G$ and $F$ are not $Cmp$-cospectral.
\bs
\label{partial-orders}
Let $G, F \in {\cal G}_n^m$.
Then
\\[1.5ex]
$(a1)$
$G \succeq ^{\phi}F \Rightarrow G 
\preceq ^ {\lambda} F$ and $G \succ\succ ^{\phi}F \Rightarrow G \prec ^ {\lambda} F$,
\\[1.5ex]
$(a2)$ 
$G \succeq ^{\phi}F \Rightarrow 
G \succeq ^p F \Rightarrow 
G \succeq ^{\tau } F \Rightarrow 
G \succeq ^{\infty} F$,
\\[1.5ex]
$(a3)$ 
$G \asymp ^{\phi} F \Leftrightarrow G \asymp ^c F \Leftrightarrow G \asymp ^p F 
\Leftrightarrow G \asymp ^{\tau } F \Leftrightarrow G \asymp ^{\infty} F \Leftrightarrow P(\lambda ,G) \equiv P(\lambda , F)$, 
\\[1.5ex]
$(a4)$ 
$\succeq ^{\phi}$, 
$ \succeq ^c $, 
$ \succeq ^p $,   
$ \succeq ^{\tau } $, 
$ \succeq ^{a(k)}$,
$ \succeq ^r $
 are partial quasi-order relations and   
 $\succeq ^t $,
 $\succeq ^\lambda $,
$ \succeq ^{\infty} $, $ \succeq ^{h_i} $ are linear 
quasi-order relations on   ${\cal G}_n^m$,
and
\\[1.5ex]
$(a5)$ 
$G \succeq ^{a(k)} F \Rightarrow G \succeq ^{r(k)} F$,
$G \succeq ^{r(k)} F \Rightarrow a^k_{n-k}(G) \ge a^k_{n-k}(F)$, and
\\[1ex]
$G \asymp ^{r(k)} F \Leftrightarrow R^k(p,G) \equiv R^k(p,F)$.
\es

{\bf Proof.}
By {\bf \ref{reciprocity}} $(a2)$, $G \succeq ^{\phi} F \Rightarrow G \succeq^ p F$.
By  {\bf \ref{P(n+r)=}}, $G \succeq ^p F \Rightarrow 
G \succeq ^{\tau } F$.
All other claims above are obvious.
\ep
\\[1.5ex]
\indent
It is interesting to compare relations 
$\succeq ^c$ and $\succeq ^p$.
Let 
\\[1.5ex]
$P^+(\lambda , G) = \sum \{(-1)^{2r} c_{2r}(G) \lambda ^{n- 1 - 2r}:  r \in \{0, \ldots , \lceil\frac{1}{2}n \rceil\}$ and
\\[1ex]
$P^-(\lambda , G) = \sum \{(-1)^{2r+1} c_{2r+1}(G) \lambda ^{n- 2r}:  r \in \{0, \ldots , \lceil \frac{1}{2}n \rceil - 1\}$.
\\[1.5ex]
\indent
We call $P^+(\lambda , G)$ and $P^-(\lambda , G)$ {\em the positive and the negative part of} $P(\lambda , G)$, respectively.
If $G, F \in {\cal G}_n^m$ and $G \succeq ^c F$, then
both positive and negative parts of $P(\lambda , G)$ are greater or equal to the positive and negative parts of $P(\lambda , F)$, respectively, for $\lambda \ge 0$.
Therefore it is very possible and not a surprise at all that
there are pairs $G, F \in {\cal G}_n^m$ such that 
$G \succeq ^p F$ but $G \not \succeq ^c F$ or
$G \succeq ^c F$ but $G \not \succeq ^p F$ or, moreover,
$G \succeq ^p F$ but $F \succeq ^c G$.
It is more surprising 
when both $G \succeq ^c F$ and  $G \succeq ^p F$.
Here are two examples illustrating the above observation.
\\[1ex]
\indent
{\bf Example 1.}
Let $F$ be a star with the center vertex $c$ and with at least three edges and let $G$ be obtained from $F - u$ by  by adding a new edge between two vertices adjacent to $c$. Then 
$G, F \in {\cal G}_n^m$ for some $n = m +1 \ge 4$,
and $G$ exactly one triangle, an isolated vertex,  exactly two components.
It is easy to show  that $G \succ ^p F$.
It is also easy to show that $c_2(G) - c_2(F) = 2m - 6 > 0$,
$c_{n-1}(G) = 0$, and $c_{n-1}(F) = n$.
Therefore  $G \not \succ ^c F$.
\\[1ex]
\indent
{\bf Example 2.}
Let $G$ be the graph having three components, namely, a triangle and two one edge components.
Let $F$ be the graph two components, namely, $P_3$ and $P_4$.
\\[0.7ex]
Then 
$P(\lambda ,G) = 
\lambda ^2 (\lambda - 2)^2 (\lambda - 3)^2 =
\lambda ^2 (\lambda ^4 - 10 \lambda + 37 \lambda ^2 - 60 \lambda + 36)$
and
\\[1ex]
$P(\lambda ,F) = 
\lambda (\lambda -1) (\lambda -2) (\lambda - 3) 
((\lambda - 2)^2 - 2) =
\lambda (\lambda ^5 - 10 \lambda ^4 + 37 \lambda ^3 - 62 \lambda ^2 + 46 \lambda - 12)$.
\\[0.7ex]
Therefore $G \succ ^p F$ and $F \succ ^c G$.
 \\[1ex]
\indent
As we will see below,  the  symmetric operations described in Section \ref{operations} turn out to be  monotone with respect to almost all  relations on graphs mentioned 
in this Section
 and 
in Section \ref{other-parameters}. 
 \\[1ex]
\indent
Let $\succeq \in 
\{ \succeq ^{a(k)},
\succeq ^c, 
\succeq ^p,   
\succeq ^{\tau }, 
\succeq ^{r(k)},   
\succeq ^t ,
\succeq ^\alpha ,
 \succeq ^{h_i} ,
 \succeq ^{\infty} 
 \}$.
In particular, we will have:
\bs
\label{succeq}
 Let $G'$ be a graph obtained from a graph $G$ by a 
 $\diamondsuit $-operation.
 Then $G \succeq G'$. 
\es

 \section{On Laplacian posets of graphs}
 \label{Lposets}
 
 \indent
 
In this section we will describe some results on the monotonicity of the operations  defined in Section 
\ref{operations} with respect to some Laplacian posets in 
Section \ref{partialOrders} and on the problems 
${\cal T}_{max}$ and  ${\cal T}_{min}$ of finding graphs with the maximum and minimum  number of spanning trees among the graphs of the same size, respectively. 
Since $t(G)$ is proportional  to $R(p,G)$ for $p$ close to zero,
some results in Section \ref{Rposet} on $R(p,G)$ provide the corresponding results on $t(G)$ 
(see, for example, {\bf \ref{maxR(G),cycl(G)<4}}).
\\[2ex]
\indent
Using {\bf \ref{posemdef}} $(a3)$ we can prove the following inequalities.
\bs {\em \cite{Kcourse,Kpolinom}}
\label{c_s(G)<c_s(K)}
Let $G \in \bar{{\cal G}}_n^m$. Then
\\[1ex]
$(a1)$
$c_s(G) \le (2m/n)^s {{n-1} \choose s}$, and in particular,  
\\[1ex]
$(a2)$
$c_s(G) \le c_s(K_n)$ for every $s \in \{0, \ldots , n-1\}$,  and so
\\[1ex]
$(a3)$ $K_n \succeq ^c G$, and moreover, 
$K_n \succ\succ ^c G$ for $G$ not isomorphic to $K_n$.
\es

Given a symmetric function  
$\delta : \mathbb{R}^n \to \mathbb{R}$
and a graph $G$,  let  
$\delta [G] =\delta (\{d(v,G): v \in V(G)\}$.
In particular, let $\delta _s[G] = \sum \{d(v,G)^s: v \in V(G)\}$.
\bs  {\em \cite{Kcourse,Kproceed}}
\label{delta(G),H-operation}
Let $G, F \in {\cal G}^m_n$ and 
$F$ be obtained from $G$ by a symmetric  hammock-operation.
Let $\delta : \mathbb{R}^n \to \mathbb{R}$ be a convex symmetric function.
Then $\delta [G] \le \delta [F]$, and so, in particular,
$\delta _s [G] \le \delta _s[F]$ for every positive integer 
$s$.
\es

\subsection{Some results on relations  
$\succeq ^{t }$,  $\succeq ^{\tau }$, and  $\succeq ^p$} 



%
\bs {\em \cite{Kcourse,KoperTrees}}
\label{t(G),H-operation}
Let $G, G' \in {\cal G}^m_n$ and 
$F$ be obtained from $G$ by a symmetric 
$H_{xy}$-operation.
Then 
\\[1ex]
$(a1)$
$t(G) \ge t(G')$ and
\\[1ex]
$(a2)$
$t(G) = t(G')$ if and only if $G$ is isomorphic to $G'$ or 
$x$ is a cut vertex in $G$.
\es 

From  {\bf \ref{GtoTreshGraph}}, 
{\bf \ref{delta(G),H-operation}}, and
{\bf \ref{t(G),H-operation}}
 we have:
\bs  {\em \cite{Kcourse,Kproceed}}
\label{thresh,H-operation}
For every  $G \in {\cal G}^m_n$ and a clique $K$ in $G$ there exist  threshold graphs $F'$ and $F''$ in ${\cal G}^m_n$ containing $K$ and such that 
$t(G) \ge t(F')$ and $\delta [G] \le \delta [F'']$, where 
$\delta : \mathbb{R}^n \to \mathbb{R}$ is a convex symmetric function.
\es

Since a $\diamondsuit $-operation is a particular case of an $H$-operation, the claim of {\bf \ref{t(G),H-operation}} is also true for a $\diamondsuit $-operation.
Therefore from {\bf \ref{properties-of-Kxy}} $(a2)$ and 
{\bf \ref{t(G),H-operation}} we have:
\bs {\em \cite{Kcourse,KoperTrees}}
\label{>t,D-operation}
Let $G, F \in {\cal G}_n^m$ and 
$F$ be obtained from $G$ by a  $\diamondsuit $-operation. Then
\\[1ex]
$t([G]^c) \ge t ([F]^c)$.
\es

From {\bf \ref{recursionPhi}} and {\bf \ref{>t,D-operation}}
we have: 
\bs {\em \cite{Kcourse,Kproceed}}
\label{Gsucceq^bDiamond(G)} 
Let $G$ be a graph and $G'$ be obtained from $G$ by 
a    $\diamondsuit $-operation.
\\[1ex]
Then $G \succeq ^{\tau } G'$. 
\es

The following is a generalization of 
{\bf \ref{Gsucceq^bDiamond(G)}} for 
an $H$-operation and $\succeq ^p$.
\bs {\em \cite{Kcourse}}
\label{Gsucceq^pH(G)} 
Let $G$ be a graph and $G'$ be obtained from $G$ by a symmetric   $H$-operation.
Then $G \succeq ^p G'$. 
\es

From {\bf \ref{K-F}} and {\bf \ref{>t,D-operation}}
 we have, in particular, the following characterization of graphs  with $n$ vertices  and 
$e(K_n) - r $ edges having  the maximum number of spanning trees  provided $r \le n/2$.
We recall  that
$rP^1$ is a graph-matching with $r$ edges and
$P^2  + kP^1$ is a disjoint union of the 2-edge path $P^2$ and the $k$-edge matching. 
\bs
\label{maxtree-graph} 
{\em \cite{Kcourse,KCh,KoperTrees}}
Let $F$ be a simple graph  with no isolated vertices and  with $r$  edges.
Suppose that $F$ is not isomorphic to 
$rP^1$ or $P^2  + (r-2)P^1$. 
Then 
\\[1.3ex]
\indent
$t (K_n - E(rP^1)) > 
t(K_n - E(P^2  + (r-2)P^1)>  t(K_n - E(F))$ 
\\[1.3ex]
for every $n$ and $r$  such that  $r \le n/2$.
\es

In \cite{Kextr} we were able to obtain more general result by giving a complete characterization of  graphs with $n$ vertices  and ${n\choose 2} - r$
edges having the maximum number of spanning trees,  where  $r \le n$. The proof of this result uses essentially  some of the above operations and  some more delicate  $\succeq ^p$-monotone operations for  special classes of graphs \cite{KgraphTrnsf}.
\\[1ex]
\indent
Here is the description of all graphs with $n$ vertices, ${n\choose 2} - n$ edges and  the maximum number of spanning trees.
\bs {\em \cite{Kextr}}
\label{maxtree-graph}
Let $Q \in {\cal G}_n^n$.
Suppose that 
$t(K_n - E(Q)) \ge t(K_n - E(F))$ for every 
$F \in {\cal G}_n^n$. Then 
\\[1ex]
$(a0)$ if $n = 0\bmod 3$, say $n = 3k$ and $k \ge 1$, then $Q = kO^3$,
\\[1ex]
$(a1)$ if $n = 1\bmod 3$, say $n = 3k +1$ and 
$k \ge 1$, then $Q = O^4 + (k-1)O^3$, and
\\[1ex]
$(a2)$ if $n = 0\bmod 3$, say $n = 3k+2$ and 
$k \ge 1$, then $Q = O^5 + (k-1)O^3$.
\es

Let, as above, $S^r$ be the graph-star with $r$ edges. 
From {\bf \ref{K-F}} and {\bf \ref{>t,D-operation}}
we have, in particular:
\bs
\label{mintree-graph} {\em \cite{Kcourse,KCh}}
Let $F$ be a simple graph with no isolated vertices, with $r$  edges, and with at most $n$ vertices. 
Suppose that $F$ is not isomorphic to $S^r$.
Then 
$t(K_n - E(F)) > t(K_n - E(S^r))$ 
for every $n$ and $r$ such that  $r \le n- 1$.
\es

From {\bf \ref{T>G}} and 
{\bf \ref{>t,D-operation}} we have, in particular:
\bs {\em \cite{Kcourse,KoperTrees}} 
\label{T>^bG}
Let $G$ be a connected graph. Then for every spanning tree $T$ of $G$ there exists a tree $D$ with $e(G)$ edges such that $D$ contains $T$ and $D\succeq ^{\tau }_n G$, where $n = e(G) +1$.
\es
 
As above, $P$ is a path, $O$  is a cycle, $S$ is a star, and
$K_4^-$ is obtained from $K_4$ by deleting one edge.
Also let $Z$ be obtained from 
a star  $S$ by adding a new vertex $x$ and a new edge between $x$ and a leaf of $S$ and let
$W$ be a windmill, i.e. 
$W$ is obtained from a star $S$ with at least two edges 
by adding an edge between two leaves of $S$.
\bs {\em \cite{Kextr}}
\label{M>>^p}
Let $G^m$ be a graph with $m$ edges.
Suppose that in each claim  below $G^m$ is  
not isomorphic to any of the special  graphs listed in this claim. We write $\succ\succ ^p$ instead of  $\succ\succ ^p_n$, where $n = 2m$.
\\[1ex]
$(a1)$ 
If $m \in \{2,3,4\}$, then
\\[1ex]
$2P^1 \succ\succ ^p P^2$,
\\[1ex]
$3P^1 \succ\succ ^p P^1 + P^2 \succ\succ ^p P^3 
\succ\succ ^p  O^3  \succ\succ ^p S^3$, and
\\[1ex]
$4P^1 \succ\succ ^p 2 P^1 + P^2 \succ\succ ^p P^1 + P^3
\succ\succ ^p  P^1 + O^3 \succ\succ ^p P^4 \succ\succ ^p P^1 + S^3 \succ\succ ^p  O^4  \succ\succ ^p  Z^4
\\[1ex]
\succ\succ ^p  W^4
\succ\succ ^p S^4$.
\\[1ex]
$(a2)$ For $m = 5$, $Z^5 \succ ^{\infty} K_4^-$ and 
$(Z^5, K_4^-)$ is the only 
$ \succ\succ ^p$-incomparable pair of graphs.
\\[1ex]
$(a3)$ 
If $m \ge 6$, then
\\[1ex]
$mP^1 \succ\succ ^p 
(m - 2)P^1 + P^2  \succ\succ ^p
(m - 4) P^1 + 2 P^2  \succ\succ ^p  
(m - 3) P^1 + P^3 \succ\succ ^p  
\\[1ex]
(m - 3)P^1 + O^3
\succ\succ ^p
 (m - 6) P^1 + 3P^2   \succ\succ ^p 
 (m - 5) P^1 + P^2 + P^3 \succ\succ ^p
(m - 4)P^1 + P^4 
\\[1ex]
\succ\succ ^p
G^m \succ\succ ^p
W^m \succ\succ ^p S^m$. 
\\[1ex]
$(a4)$ 
If $m \ge 7$, then $G^m  \succ\succ ^p
Z^m
\succ\succ ^p
W^m \succ\succ ^p S^m$.
\es

The following is a generalization of 
{\bf \ref{Gsucceq^bDiamond(G)}} 
for $({\cal G}_n^m, \succ ^p)$.
\bs {\em \cite{Kproceed,Kcourse}}
\label{Gsucceq^pH(G)} 
Let $G$ be a simple graph and $G'$ be obtained from $G$ by a  $\diamondsuit$-operation.
Then $G \succeq ^p G'$. 
\es

Obviously, $G - e \succ ^\tau G$.   
It turns out that this inequality remains true 
for $\succ ^p$.
\bs {\em \cite{Kproceed,Kcourse}}
\label{Gsucceq^pD(G)} 
Let $G$ be a simple graph and $e$ an edge of $G$.
Then 
\\[1ex]
\indent
$(\lambda -1) P(\lambda , G - e) > \lambda P(\lambda , G) $ for $ \lambda \ge n$, 
and so $G - e \succ ^p G$. 
\es

\begin{figure}
\vskip -2mm
\begin{center}
\includegraphics[scale=0.365]{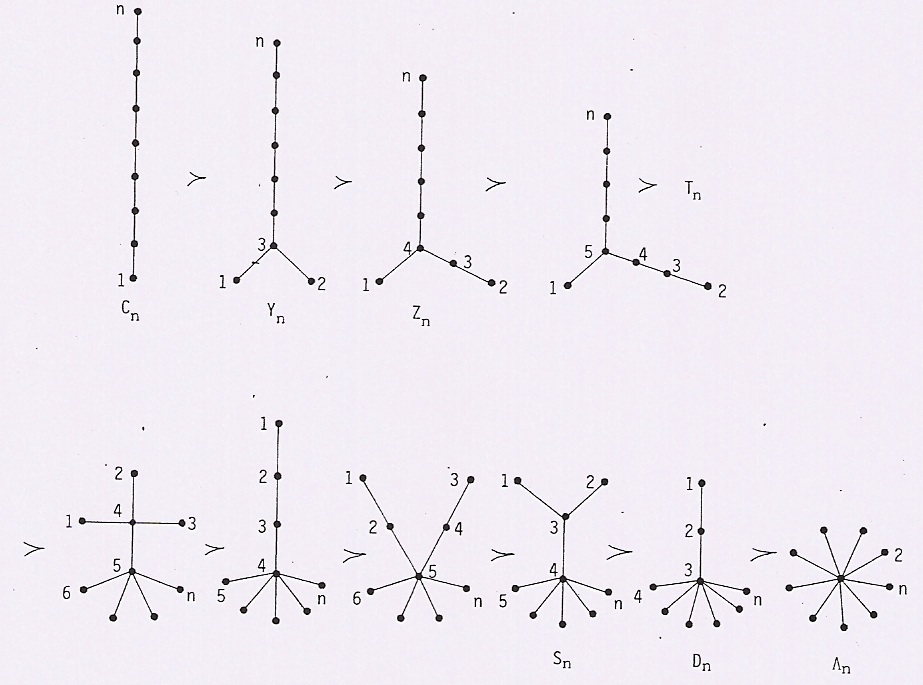}
\caption{Four   $\succ\succ ^{p}$-``best'' and six  $\succ\succ ^{p}$-``worst'' trees in ${\cal T}_n$, $n \ge 9$. }
\label{sucsucT}
\end{center}
\end{figure}

Here are some results  on 
$\succ\succ ^p$-comparison of $n$-vertex trees that we were able to prove using the tree operations from \cite{Kcmpr}. 
Notice that if $G$ and $F$ are  trees with the same number of vertices, then  $P(\lambda , G) = {\Phi } (\lambda , G) $.
\\[1ex]
\indent
Let  $T(a,b,c)$ denote the tree $T$ from 
${\cal S}(r)$ with $w(T) = (a, b, c)$, where $a \le b \le c$.
\bs {\em \cite{Kcmpr}}
\label{bestTr}
 Let  $T_n$ be a tree with $n$ vertices.   The trees are numbered according to 
  the Harary list of ${\cal T}_{n}$  with $n \le 10$ {\em (see \cite{H})}.
 \\[1ex]
 $(a1)$ If $1 \le n \le 9$,
 then $\succ\succ ^ p$ is a linear order on ${\cal T} _n$. 
 In particular, we have:
 \\[1ex]
 $(a1.1)$  
 the $\succ\succ ^ p$-order of the trees in ${\cal T} _7$ is 
\\[1ex]
$1,~2,~3,~4,~9,~10,~5,~6,~11,~7,~8$,
\\[1ex]
 $(a1.2)$
 the $\succ\succ ^ p$-order of the trees in ${\cal T} _8$ is 
\\[1ex]
$1,~2,~3,~5,~4,~13,~14,~17,~16,~15,~23,~6,~7,~8,~18,~20,~19,~22,~9,~10,~21,~11,~12 $,
and
\\[1ex]
 $(a1.3)$ 
 the $\succ\succ ^ p$-order of the trees in ${\cal T} _8$ is
 $\succ\succ ^ p$-order is 
\\[1ex]
$1,~2,~3,~4,~5,~6,~19,~20,~23,~22,~21,~24,
~26,~25,~27,~44,~42,~45,~7,~8,~9,~10,~11,~28,~30
\\[0.5ex]
29,~40,~39,~31,~38,~32,~46,~47,~33, ~41, ~12, ~13, ~14, ~37, ~35, ~36, ~42, ~15, ~16, ~34, ~17, ~18$.
\\[1ex]
 $(a2)$ There are  $\succ ^ p$-non-comparable trees in ${\cal T} _{10}$.
\\[1ex]
 $(a3)$ 
 If $n \ge 10$, then  in ${\cal T}_n$ there exist
 four  successively $\succ\succ ^p$-best and six successively 
\\[0.5ex] 
 $\succ\succ ^{\phi}$-worst trees $($see  Figure \ref{sucsucT}$) $, and in particular,
 \\[1ex]
$P_n \succ\succ ^ p
 T(1,1, n-3) \succ\succ ^p
T(1,2, n-4) \succ\succ ^p
T(1,3, n-5) \succ\succ ^p
T_n\succ\succ ^p
Z_n 
\succ\succ ^p S_n$,  
\\[0.5ex]
where $T_n$  is not isomorphic to any of the special trees listed above. 
\es

Here are some more details on  
$\succ\succ ^p$-comparison of $n$-vertex trees with $n \le 10$.
\bs 
\label{bestTr}
 Let  $T_n$ be a tree with $n$ vertices.  The trees are numbered according
 to the Harary list of ${\cal T}_{n}$  with $n \le 10$ {\em (see \cite{H})}.
 \\[1.5ex]
 $(a1)$ If $1 \le n \le 7$, then $\succ\succ ^ p ~=~ \succ\succ ^ {\phi } ~=~ \succ\succ ^ c$.
 \\[1.5ex]
 $(a2)$ In ${\cal T} _8$~ 
 $\succ\succ ^ p ~= ~\succ\succ ^ p_x$  for $x = 5,6639$ ~but~
$\succ\succ ^ p ~\ne~ \succ\succ ^ {\phi }$, namely, there are two 
$\succ\succ ^ {\phi }$-non-comparable  pairs 
$(4, 13)$
and 
$(8, 18)$ of trees in ${\cal T} _8$, 
namely, $\lambda (8) = 5.236$, 
\\[0.5ex]
$\lambda (18)= 5.125$, 
and
$0 \le L(\lambda , 8) < L(\lambda , 18)$ for $\lambda (8) \le \lambda < x$,
$L(x , 8) = L(x , [18])$, and 
\\[0.5ex]
$L(\lambda , 8) > L(\lambda , 18)$ for $\lambda > x$. There is one 
$\succ\succ ^ c$-non-comparable  pair 
$(16, 6)$ 
in ${\cal T} _8$, 
namely, $c_0(16) = c_0(6) = 1$, $c_1(16) = c_1(6) = 14$,  $c_5(16) = c_5(6) = 204$, $c_7(16) = c_7(6) = 8$, $c_i(16) < c_i(6)$ for 
$2 \le i \le 4$, and $c_6(16) > c_6(6)$.
\\[1.5ex]
 $(a3)$
In ${\cal T} _9$~ 
$\succ\succ ^ p ~= ~\succ\succ ^ p_x$~ for $x = 6.842$~ but~ 
$\succ\succ ^ p ~\ne~ \succ\succ ^ {\phi }$. There are 26
\\[0.5ex] 
$\succ\succ ^ {\phi }$-non-comparable pairs of trees in ${\cal T} _9$, for example,  
pair $(14,37)$ is one of them,
namely, 
$x = 6.84129$ is the ``crossing point'' for pair $(14,37)$,
$\lambda (14) = 6.147$, 
\\[0.4ex]
$\lambda (37) = 6.062$, 
$0 \le L(\lambda , 14) < L(\lambda , 37)$ for $\lambda (14) \le \lambda < x$,
$L(x , 14) = L(x , 37)$, and 
\\[0.5ex]
$L(\lambda , 14) > L(\lambda , 37)$ for $\lambda > x$.
\\[1.5ex]
 $(a4)$  
 $(71,82)$, 
$(70,82)$,
and 
 $(69,82)$ 
  are  $\succ ^ p$-non-comparable pairs of trees in ${\cal T} _{10}$.
In particular, 
\\[0.5ex]
$x = 11.4772 > 10$ is the ``crossing point'' for pair 
 $(71,82)$,
 namely, $\lambda (71) = 7.119$, 
 \\[0.5ex]
 $\lambda (82) = 6.702$,
 $0 \le L(\lambda, 71) < L(\lambda, 82)$ for $\lambda (71)  \le \lambda < x$, $L(x, 71 ) = L(x, 82 )$, 
\\[0.5ex]
 and  $L(\lambda, 71) > L(\lambda, 82)$ for $\lambda > x  $. 
\es

\subsection{Some results on relation $\succeq ^\phi $}
 \label{phi}

\indent

Using the recursion {\bf \ref{recursionPhi}}
we  proved by induction the following inequalities.
\bs {\em \cite{Kcourse,KCh}}
\label{Msuccsucc^phi...} 
Let $G^m$ be a graph with $m$ edges not isomorphic to any of the special  graphs listed below.
Then  $mP^1 \succ\succ ^{\phi} 
P^2 + (m - 2)P^1 \succ\succ ^{\phi}  G^m \succ\succ ^{\phi} S^m$.
\es

\bs {\em \cite{Kcourse,Kcmpr}}
\label{P^m>G^m}
Let $G^m$ be  a connected graph with $m\ge 4$ edges. 
If $G^m$ is not a path, 
then $P^m \succ\succ ^{\phi} G^m$.
\es

In \cite{Kcmpr} 
we  were able to find an inductive proof of {\bf \ref{P^m>G^m}}  using  {\bf \ref{recursionPhi}}. 
On the other hand, we found some $\succeq ^{\phi}$-increasing operations on trees that 
allow us to give another  proof 
of {\bf \ref{P^m>G^m}}  and some other results on the poset $({\cal T}_n, \succeq ^{\phi})$. It turns out that the revers of an $(x,y)$-path operation is one of $\succeq ^{\phi}$-increasing operations on trees and that every tree $T_n$ which is not a path can be transformed to the path $P_n$ by a series of such operations.
Later we used these methods to prove similar results for
the poset $({\cal G}_n^m, \succeq ^c)$. In particular, we found that some of the above mentioned  $\succ\succ ^{\phi}$-increasing operations on trees are also   $\succ\succ ^c$-increasing operations 
(see  {\bf \ref{Psucc^cG}}  below).
\\[1ex]
\indent
Using {\bf \ref{recursionPhi}} we obtain  from {\bf \ref{P^m>G^m}}:
\bs {\em \cite{Kcmpr}}
\label{O^m>G^m}
Let $G^m$ be  an edge 2-connected graph with $m\ge 5$ edges.  If $G^m$ is not a cycle, then  $O^m \succ\succ ^{\phi} G^m$.
\es

\bs {\em \cite{Kcmpr}}
\label{Pi+Pj} Let $i$ and $j$ be positive integers such that $2 \le j+1 \le i - 1$.
Then 
\\[1ex]
$\Phi (\lambda , P_{i - 1}  + P_{j + 1}) > 
\Phi (\lambda , P_{i }  + P_{j })$ for $\lambda \ge 4$, i.e.
$ P_{i - 1}  + P_{j + 1} \succ\succ ^\phi_4 
P_i   + P_j $.
\es

In \cite{Kcmpr2} (which is a continuation of \cite{Kcmpr}) we obtained similar $\succeq ^{\phi}_x$-comparison  results for some other  classes of graphs. 
In particular, we extended {\bf \ref{Pi+Pj}} to the class of graphs whose each component is either a path or a cycle. 
Here are some of these results.
\\[1ex]
\indent
We remind that if $G \succ ^{\phi}_x F$ for $x =  \max \{\lambda (G), \lambda (F)\}$, we write 
$ \succ ^{\phi }$ instead of  $\succ ^{\phi }_x$,
and so in this case $\lambda (G) \le \lambda (F)$.
Also notice that $3 \le \lambda (O_{2s +1}) <  \lambda (O_{2k}) = 4$.
\bs {\em \cite{Kcmpr2}}
\label{O,4n}
Let $n$ and $k$  be  integers such that $n \ge 2$ and $k \ge 0$. Then
\\[1ex]
$(a0)$ 
$O_{2n - 1}+ O_{2n + 1} \succ\succ ^{\phi} ~~
O_{4n}~\succ\succ ^{\phi} ~ 2O_{2n}$,
\\[1ex]
$(a1)$ $
O_{2n - 2k - 3}+ O_{2n + 2k + 3}  \succ\succ ^{\phi} _x~
O_{2n - 2k - 1}+ O_{2n + 2k+ 1}
\succ\succ ^{\phi} _x~ O_{2n - 1}+ O_{2n + 1}$ for 
\\[1ex] 
$n \ge k +3$ and $x = 4.05$, 
and 
\\[1ex]
$(a2)$ $ 2O_{2n} \succ\succ ^{\phi} ~
O_{2n - 2k}+ O_{2n + 2k}  \succ\succ ^{\phi} ~
O_{2n - 2k - 2}+ O_{2n + 2k+ 2}$  
for $n \ge k + 3$.
\es
\bs {\em \cite{Kcmpr2}}
\label{O,4n+2}
Let $n$ and $k$  be  integers such that $n \ge 1$ and 
$k \ge 0$. Then
\\[1ex]
$(a0)$ 
$2O_{2n + 1} \succ\succ ^{\phi} ~~
O_{4n + 2}~\succ\succ ^{\phi} ~ O_{2n} + O_{2n +2} $,
\\[1ex]
$(a1)$ $
O_{2n - 2k - 1}+ O_{2n + 2k+3}  \succ\succ ^{\phi} _x~
O_{2n - 2k +1 }+ O_{2n + 2k+ 1}
\succ\succ ^{\phi} _x~ 2O_{2n +1}$ for  
$n \ge k +1$ and 
\\[1ex]
$x = 4.05$,
and 
\\[1ex]
$(a2)$ $O_{2n} + O_{2n +2} \succ\succ ^{\phi} ~
O_{2n - 2k}+ O_{2n + 2k}  \succ\succ ^{\phi} ~
O_{2n - 2k - 2}+ O_{2n + 2k+ 2}$  
for $n \ge k + 3$.
\es

\bs {\em \cite{Kcmpr2}}
\label{O,4n +1}
Let $n$ and $k$  be  integers such that $n \ge 2$  and 
$k \ge 0$. Then
\\[1ex]
$(a0)$ 
$O_{2n - 1} + O_{2n+2} \succ\succ ^{\phi} _x~~
O_{4n +1}~\succ\succ ^{\phi} ~ O_{2n} + O_{2n +1}$ for $x = 4.133$,
\\[1ex]
$(a1)$ $
O_{2n - 2k -1}+ O_{2n + 2k + 2}  \succ\succ ^{\phi} ~
O_{2n - 2k +1}+ O_{2n + 2k}
\succ\succ ^{\phi} ~ O_{2n - 1}+ O_{2n + 2}$ for 
$n \ge k + 2$,
\\[1ex]
and 
\\[1ex]
$(a2)$ $O_{2n - 1} + O_{2n+2} \succ\succ ^{\phi} ~O_{2n} + O_{2n +1} \succ\succ ^{\phi} ~
O_{2n - 2k }+ O_{2n + 2k+1}  \succ\succ ^{\phi} ~
O_{2n - 2k -2}+ O_{2n + 2k +3}$  for
$n \ge k + 2$.
\es

\bs {\em \cite{Kcmpr2}}
\label{O,4n-1}
Let $n$ and $k$  be  integers such that $n \ge 2$  and 
$k \ge 0$. Then
\\[1ex]
$(a0)$ 
$O_{2n - 1} + O_{2n} \succ\succ ^{\phi} _x~~
O_{4n -1}~\succ\succ ^{\phi} ~ O_{2n-2} + O_{2n +1}$ for $x = 4.325$,
\\[1ex]
$(a1)$ $
O_{2n - 2k -1}+ O_{2n + 2k}  \succ\succ ^{\phi} ~
O_{2n - 2k +1}+ O_{2n + 2k - 2}
\succ\succ ^{\phi} ~ O_{2n - 1}+ O_{2n}$ for 
$n \ge k + 2$,
\\[1ex]
and 
\\[1ex]
$(a2)$ $O_{2n - 1} + O_{2n} \succ\succ ^{\phi} 
O_{2n-2} + O_{2n +1} \succ\succ ^{\phi} ~
O_{2n - 2k }+ O_{2n + 2k -1}  \succ\succ ^{\phi} ~
O_{2n - 2k -2}+ O_{2n + 2k +1}$  for
$n \ge k + 2$.
\es

We also obtain similar results for graphs $O_n + P_k$.
Here are some of them.
\bs {\em \cite{Kcmpr2}}
\label{P(2k)O(2n+1)?P(2k+1)O(2n)}
Let $n$ and $k$  be  integers such that $n \ge 2$ and 
$k \ge 2$. Then
\\[1ex]
$(a1)$ if $n \ge 2k+1$, then
$P_{2k} + O_{2n+1} \succ\succ ^{\phi} _x~~ P_{2k+1} + O_{2n +1}$ for 
$x \ge 4.74$,
\\[1ex]
$(a2)$ 
 if $k \le n < 2k+1$, then
$P_{2k+1} + O_{2n +1} \succ\succ ^{\phi} _x~~ P_{2k} + O_{2n+1}$ for $x > 4$
and 
\\[1ex]
$(a3)$
 if $n <  k$, then
$P_{2k} + O_{2n+1} \succ\succ ^{\phi} _x~~ P_{2k+1} + O_{2n +1}$ for $ x \ge 4 $.
\es

\bs {\em \cite{Kcmpr2}}
\label{P(2k+1)O(2n)?P(2k-1)O(2n+2)}
Let $n$ and $k$  be  integers such that 
$n \ge 2$ and  $k \ge 2$. 
Then
\\[1ex]
$(a1)$ if $n <  k$ or $n \ge 2k$, then
$P_{2k+1} + O_{2n} \succ\succ ^{\phi} _x~~ 
P_{2k-1} + O_{2n + 2}$ for $x > 4$ and
\\[1ex]
$(a2)$ 
 if $k \le n < 2k$, then
 $P_{2k-1} + O_{2n +2} \succ\succ ^{\phi} _x~~ 
P_{2k+1} + O_{2n}$ for $x >4$.
\es

\bs {\em \cite{Kcmpr2}}
\label{P(2k+2)O(2n)>P(2k)O(2n+1)}
Let $n$ and $k$  be  integers such that 
$n \ge 2$ and  $k \ge 1$. 
Then
  \\[1ex]
  \indent
$P_{2k+2} + O_{2n} \succ\succ ^{\phi} _x~~ 
P_{2k} + O_{2n + 1}$ for $x > 3.5$.
\es

\bs {\em \cite{Kcmpr2}}
\label{P(2n-3)O(4)>P(1)O(2n)}
Let $n$ and $k$  be  integers such that 
$n \ge 2$. 
Then
  \\[1ex]
  \indent
$P_{2n - 3} + O_{4} \succ\succ ^{\phi} _x~~ 
P_{1} + O_{2n}$ for $x > 4$.
\es

From  {\bf \ref{O,4n}} - {\bf \ref{O,4n-1}} we have the following corollaries.
 \bs {\em \cite{Kcmpr2}}
\label{O,2n,>p}
All inequalities in {\bf \ref{O,4n}} - {\bf \ref{O,4n-1}} remain true if the relations $\succ\succ ^{\phi} $ and 
$\succ\succ ^{\phi} _x$ are replaced by $\succ\succ ^p$.
\es

\bs {\em \cite{Kcmpr2}}
 \label{>Max} 
 Let $M_n$ be a $\succeq ^{\phi}_x$-maximal 2-regular graph in 
 ${\cal G}_n^n$ for $x = 4.05$.
 Then
 \\[1ex]
$(a0)$ if $n = 0\bmod 3$, say $n = 3k$ and $k \ge 1$, then $M_n = kO_3$,
\\[1ex]
$(a1)$ if $n = 1\bmod 3$, say $n = 3k +1$ and 
$k \ge 1$, then $M_n = O_4 + (k - 1)O_3$, and
\\[1ex]
$(a2)$ if $n = 2\bmod 3$, say $n = 3k+2$ and 
$k \ge 1$, then $M_n = O_5 + (k-1)O_3$.
 \es

 \bs {\em \cite{Kcmpr2}}
 \label{>Min} 
 Let $W_n$ be a $\succeq ^{\phi}$-minimal  2-regular graph in  ${\cal G}_n$.
 Then
 \\[1ex]
$(a0)$ if $n = 0\bmod 4$, say $n = 4k$ and $k \ge 1$, then $W_n = kO_4$,
\\[1ex]
$(a1)$ if $n = 1\bmod 4$, say $n = 4k +1$ and 
$k \ge 1$, then $W_n = O_5 + (k-1)O_4$,
\\[1ex]
$(a2)$ if $n = 2\bmod 4$, say $n = 4k +2$ and 
$k \ge 1$, then $W_n = O_6 + (k-1)O_4$, and
\\[1ex]
$(a3)$ if $n = 3\bmod 4$, say $n = 4k+3$ and 
$k \ge 1$, then $W_n = O_7 + (k-1)O_4$.
 \es

\bs {\em \cite{Kcmpr2}}
 \label{>MaxBip} 
 Let $BM_{2n}$ be a $\succeq ^{\phi}$-maximal bipartite 2-regular graph  in  ${\cal G}_{2n}$.
 \\[1ex]
 \indent
 Then 
 $BM_{2n} = O_{2n}$.
 \es

\bs {\em \cite{Kcmpr2}}
 \label{>MinBip} 
 Let $BW_{2n}$ be a $\succeq ^{\phi}$-minimal bipartite graph 
 in  ${\cal G}_{2n}$.
 Then
 \\[1ex]
$(a0)$ if $n = 0\bmod 2$, say $2n = 4k$ and $k \ge 1$, then $BW_{2n}= kO_4$ and
\\[1ex]
$(a1)$ if $n = 1\bmod 2$, say $2n = 4k + 2$ and 
$k \ge 1$, then $BW_{2n} = O_6 + (k-1)O_4$.
 \es

Now using {\bf \ref{somecoefficients}}$(a4)$ and  Reciprocity Theorem {\bf \ref{reciprocity}} we have from  
{\bf \ref{O,2n,>p}} - {\bf \ref{>Min}}:
 \bs {\em \cite{Kcmpr2}}
 \label{MaxMinTreeG}
 Let $F$ be a 2-regular  graph with $2n$ vertices
{\em  (and so every component of $F$ is an even cycle)}.
Suppose that $F$ is  not isomorphic to $M_n$ in  
{\bf \ref{>Max}} or to $W_n$ in {\bf \ref{>Min}}.
Then 
\\[1ex]
\indent
$t(K_{n+s} \setminus E(M_n)) >
t(K_{n+s} \setminus E(F)) > t(K_{n+s} \setminus E(W_n))$ 
\\[1ex]
for every  non-negative integer $s$. 
\es

Also  using {\bf \ref{somecoefficients}}$(a4)$ and  Reciprocity Theorem {\bf \ref{bipartite-reciprocity}} for regular bipartite graphs we obtain from  {\bf \ref{O,2n,>p}}, {\bf \ref{>MaxBip}}, and   {\bf \ref{>MinBip}} the description of regular bipartite graphs with $2n$ vertices and 
$n^2 - 2n$ edges having the maximum and the minimum number of spanning trees:
 \bs {\em \cite{Kcourse}}
 \label{MaxMinTreebipG}
 Let $B$ be a 2-regular bipartite graph with $2n$ vertices
{\em  (and so every component of $B$ is an even cycle)}.
Suppose that $B$ is not a cycle and not isomorphic to 
$BW_{2n}$ in {\bf \ref{>MinBip}}.
Then 
\\[1ex]
\indent
$t(K_{n,n} \setminus E(BM_{2n})) >
t(K_{n,n} \setminus E(B)) > t(K_{n,n} \setminus E(BW_{2n}))$.
\es

We were also able to prove the following inequalities in addition to  {\bf \ref{bestTr}} $(a3)$.
\bs {\em \cite{Kcourse}}
\label{P^m>^fO^m>}
If $m \ge 4$, then $P^m \succ\succ ^ {\phi}
 O^m \succ\succ ^ {\phi}
 T(1,1, m-2)$.
\es

Using some $\succ\succ ^\phi $-monotone operations on $n$-vertex trees from \cite{Kcmpr}, we have obtained the following two theorems on the posets 
$({\cal D}_n(r), \succeq ^\phi )$ and 
$({\cal L}_n(r), \succeq ^\phi )$ (see the definitions in 
{\bf \ref{D(r),Y}}   and   {\bf \ref{L(r),S(r)}}). 
It turns out that similar results hold for $\succeq ^c $ as well (see {\bf \ref{Yr>cDr}} and {\bf \ref{Yr>^cLr}}).
\\[1ex]
\indent
Let $\delta '(G)$ be the minimum non-leaf vertex degree in $G$.
\bs {\em \cite{Kproceed,Kcourse}}
\label{Dr<Yr} 
Let $r \ge 3$ and $n \ge r + 2$.
\\[1ex]
$(a1)$ for every $D \in {\cal K}_n(r) \setminus {\cal K}_n(r)$ there exists 
$Y \in {\cal K}_n(r)$ such that $D \succ\succ  ^{\phi} Y$,  
\\[1ex]
$(a2)$  $D \succ\succ  ^{\phi} K_n(r)$ for every 
$D \in {\cal K}_n(r) \setminus \{K_n(r)\}$, and therefore, 
{\em (from $(a1)$ and $(a2)$)},
\\[1ex]
$(a3)$ 
$D  \succ\succ  ^{\phi} K_n(r)$ for every 
$D \in {\cal D}_n(r) \setminus \{K_n(r)\}$,
\\[1ex]
$(a4)$
$({\cal D}_n(3), \succeq ^\phi )$ is a linear order poset, namely, for $T, T' \in {\cal D}_n(3)$ we have:
\\[1ex]
 $T \succ\succ ^{\phi} T' \Leftrightarrow 
 \delta '(T) > \delta ' (T')$, and
\\[1ex]
$(a5)$
$({\cal D}_n(4), \succeq ^\phi )$ is a linear  poset, namely, for $T, T' \in {\cal D}_n(4)$ we have:
\\[1ex]
 $T \succ\succ ^\phi T' \Leftrightarrow  T > ^u T' $.
\es

\bs {\em \cite{Kproceed,Kcourse}}
\label{Tr<Yr} 
Let $r \ge 3$ and $n \ge r + 2$.   
Then 
\\[1ex]
$(a0)$ $L_n(r) \succ\succ  ^\phi  L_n(r + 1)$ for every 
$r \in \{2, \dots , n - 2\}$, 
\\[1ex]
$(a1)$
$({\cal S}_n(r), \succeq ^\phi )$ is a linear  poset, namely, for $T, T' \in {\cal S}_n(r)$ we have:
\\[1ex]
 $T \succ\succ ^\phi T' \Leftrightarrow  T > ^w T' $,
\\[1ex]
$(a2)$
$M_n(r) \succ\succ  ^{\phi} L $ for every 
$L \in {\cal S}_n(r) \setminus \{M_n(r)\}$,
\\[1ex]
$(a3)$ for every $L \in {\cal L}_n(r) \setminus {\cal S}_n(r)$ there exists $Z \in {\cal S}_n(r)$ such that $L \succ\succ  ^{\phi} Z$, and
\\[1ex]
$(a4)$ 
$L  \succ\succ  ^{\phi} L_n(r)$ for every 
$L \in {\cal S}_n(r) \setminus \{ L_n(r)\}$, 
therefore
{\em (from $(a3)$ and $(a4)$)}
\\[1ex]
$(a5)$
$L  \succ\succ  ^{\phi} L_n(r)$ and, 
in particular, $\lambda (L) >  \lambda (L_n(r))$
for every 
$L \in {\cal L}_n(r) \setminus \{L_n(r)\}$.
\es

Let $T$ be an $n$-vertex tree of maximum degree $r$.
Then $T$ can be transformed by some $\succ^{\phi}$-increasing operations from  from \cite{Kcmpr},  to an $n$-vertex star-tree $S$ with $r$ leaves.  Therefore by {\bf  \ref{Tr<Yr}} 
$(a2)$, we have:
 \bs {\em \cite{Kcourse}}
 \label{TwithMaxDg=r,succ^phi}
 Let $T$ be an $n$-vertex tree of maximum degree $r$ and $T$ is not isomorphic to $M_n(r)$.
 Then $M_n(r)  \succ\succ  ^{\phi} T$.
\es

As we mentioned in {\bf \ref{partial-orders}} $(a1)$, 
$G \succeq ^{\phi} F \Rightarrow G \preceq ^{\lambda} F$ and 
$G \succ\succ ^{\phi} F \Rightarrow 
G  \prec ^{\lambda}F$. 
Therefore, the  replacement of  
$\succeq ^{\phi}$ by $\preceq ^{\lambda}$ and
$ \succ\succ ^{\phi}$ by $ \prec ^{\lambda}$
in any theorem results in another correct theorem.
For example, from
{\bf \ref{Msuccsucc^phi...}} and {\bf \ref{P^m>G^m}},
we have, in particular:
\bs {\em \cite{Kcourse,Kcmpr}}
\label{l(P)<l(G)<l(S)}
Let $G^m$ be a  graph with $m$ edges 
and with no isolated vertices.
 \\[1ex]
$(a1)$ if $G^m$ is not isomorphic to $S^m$, then 
$\lambda (G^m) < \lambda (S^m)$,
\\[1ex]
$(a2)$ if $m \ge 4$ and $G^m$ has a vertex of degree at least three, then
$\lambda (P^m) < \lambda (G^m)$,  
\\
and so
\\[1ex]
$(a3)$ if $G^m$ is a tree not isomorphic to $P^m$ or $S^m$, then 
$\lambda (P^m) < \lambda (T^m) < \lambda (S^m)$.
\es

Here is an alternative proof of {\bf \ref{l(P)<l(G)<l(S)}}.
\\[1ex]
\indent
{\bf Proof.}
We prove $(a1)$. 
Obviously, $\lambda (S^m) = m+1$. Let $Q$ be a component of $G^m$ with the maximum number of vertices. Let $v(Q) = n'$. Then $n' \le m +1$.
By {\bf \ref{AtimesB}} $(a1)$, $\lambda (G^m) \le n'$.
If $n' < m +1$, then we are done. Therefore let $n' = m +1$. Then $Q$ is a tree with $m$ edges. Since $G^m$ has no isolated vertices, clearly $Q = G^m$. Hence $Q$ is a tree not isomorphic to $S^m$. Then the complement of $Q$ is connected. Thus, by  {\bf \ref{reciprocity}}, 
$\lambda (Q) < n' = m+1$.
\\[0.5ex]
\indent
Now we prove $(a2)$.
Since $G^m$ has a vertex of degree at least three, 
$\lambda (G^m) \ge 4$, 
by {\bf \ref{lambda-bounds}} $(a3)$.
By {\bf \ref{lambda-bounds}} $(a1)$, 
$\lambda (P^m) \le 4$, 
and so  $\lambda (P^m) \le \lambda (G^m)$.
It is known \cite{CDS} that
\\[1ex]
$\lambda (P^m) = 4~ sin ^2 (\frac{m}{2m+2} \pi)$.
Thus,  $\lambda (P^m) < 4 \le \lambda (G^m)$.
\ep
\\[1.5ex]
\indent
Claim $(a2)$ in {\bf \ref{l(P)<l(G)<l(S)}}  also follows from   {\bf \ref{bipartite,diamondsuit}} below.
\\[1.5ex]
\indent
There are infinitely many examples showing  that theorem {\bf \ref{Gsucceq^pD(G)}} on the 
$\succeq ^p$-monotonicity of the $\diamondsuit $-operation
is not  true if $\succeq ^p$ is replaced by 
$\succeq ^\phi$. However, it may be true under some additional condition.
\bs {\em \cite{Kcourse}}
\label{bipartite,diamondsuit}
Let $G \in {\cal G}$ and $G'$ be obtained from $G$ by the 
$\diamondsuit _{x,y}$-operation.
Suppose that  $G$ is a connected graph and $G'$ is a bipartite graph. Then
\\[1ex]
$(a1)$ 
$\lambda (G) \le \lambda (G')$, and moreover,
\\[1ex]
$(a2)$
$G \succeq ^\phi   G'$.
\es

{\bf Proof} 
(uses 
{\bf \ref{lambda-bounds}}, 
{\bf \ref{recursionPhi}}, 
{\bf \ref{alpha(As)<alpha(Az)}},
and
{\bf \ref{bip=odd}}).
Let $L^+(G) = D(G) + A(G)$
and $\lambda ^+(G)$ the maximum eigenvalue of 
$L^+(G)$.
\\[1ex]
\indent
Since $G$ is a connected graph, we have  from Theorem 8.4.5 in \cite{HJ}:
\\[1ex]
{\sc Claim 1.}
{\em 
$\lambda (G) \le \lambda ^+(G)$.}
\\[1.5ex]
${\bf (p1)}$
We prove claim $(a1)$ in our theorem.
By {\sc Claim 1}, $\lambda (G) \le \lambda ^+(G)$.
By {\bf \ref{alpha(As)<alpha(Az)}} below, 
$\lambda ^+(G) \le \lambda ^+(G')$.
Since $G'$ is a bipartite graph, by {\bf \ref{bip=odd}}, 
$\lambda ^+(G') = \lambda (G')$. Thus,
\\[1ex]
\indent
$\lambda (G) \le \lambda ^+(G) \le \lambda ^+(G') = \lambda (G')$.
\\[1.5ex]
${\bf (p2)}$ Now we prove  claim $(a2)$ in our theorem.
\\[0.7ex]
\indent
Let $X = N(x,G) \setminus (N(y,G) \cup y)$ and 
$Y = N(y,G) \setminus (N(x,G) \cup x)$.
Obviously, $e(G) = e(G')$.
For  $e \in E(G)$ let $\alpha (e) = e$  if $e \not \in [x,X]$ and $\varepsilon (e) = e'$ if $e = xs \in  [x,X]$ and $e' = ys$ for some $s \in X$. Since $G'$ is obtained from $G$ by   
the $\diamondsuit _{x,y}$-operation, $\varepsilon $ is a bijection from $E(G)$ to $E(G')$.
By definition of $\succeq ^\phi $,
$G \succeq ^\phi   G'$ if and only if 
$\lambda (G) \le \lambda (G')$ and $\Phi (\lambda , G) \ge \Phi (\lambda , G')$ for $\lambda \ge \lambda (G')$.

We prove claim $(a2)$ by induction on $e(G)$.
Our claim is obviously true for $e(G) = 0$.
Now we suppose that our claim is true for $e(G) = m - 1$ and we have to prove that it is also true for $e(G) = m$, where $m \ge 1$.

Let $E = E(G)$, $\lambda ' = \lambda (G')$, and
$\Delta (\lambda  , G) = 
\Phi (\lambda  , G) - \Phi (\lambda  , G')$.
Obviously, 
\\
$ \diamondsuit_{x,y} (G - e) = G' - \varepsilon (e)$. 
Therefore 
$\Delta (t , G - e) = 
\Phi (t , G - e) - \Phi (t , G' -  \varepsilon (e))$.
\\[1ex]
By {\bf \ref{recursionPhi}}, 
$\Phi (\lambda , G)  = \Phi (s , G)  + 
\sum \{\int  _s^\lambda \Phi (t , G - e) dt : e \in E\}$.
\\[1ex]
Therefore 
$\Delta (\lambda , G)  = \Delta (s , G)  + 
\sum \{\int  _s^\lambda \Delta (t , G - e) dt : e \in E\}$.
\\[1ex]
Since $G'$ is bipartite, $G' - \varepsilon (e)$ is also bipartite. Therefore by the induction hypothesis,
\\[1ex]
$\Delta (t , G - e) \ge 0$ for 
$ t \ge \lambda (G' - \varepsilon (e))$.
By {\bf \ref{lambda-bounds}}, 
$\lambda (G' - \varepsilon (e)) \le \lambda (G') = 
\lambda '$. 
\\[1ex]
Therefore from the above inequality we have:
$\Delta (t , G - e) \ge 0$ for 
$ t \ge \lambda (G') = \lambda '$.
\\[1ex]
By our claim $(a1)$ (that we  already proved in 
${\bf (p1)}$), 
$ \lambda (G) \le  \lambda (G') = \lambda '$.
\\[1ex]
Therefore 
$\Delta (\lambda ', G) = \Phi ((\lambda ', G) \ge 0$.
Thus, if $\lambda \ge \lambda '$, then
\\[1ex]
\indent
$\Delta (\lambda , G)  = \Delta (\lambda ' , G)  + 
\sum \{\int  _{\lambda '}^\lambda \Delta (t , G - e) dt : e \in E\} \ge 0$.
\ep
 \\[1.5ex]
 \indent
We can also prove the following generalization of 
{\bf \ref{bipartite,diamondsuit}}.
\bs {\em \cite{Kcourse}}
\label{bipartite,K_xy}
Let $G \in {\cal G}$ and $G'$ be obtained from $G$ by 
a symmetric $K_{xy}$-operation.
If $G'$ is a bipartite graph, then $G \succeq ^\phi   G'$.
\es

It turns out that 
a symmetric $K_{xy}$-operation 
(and in particular, the $\diamondsuit $-operation) 
is also $\succeq ^{\alpha }$-monotone.
\bs {\em \cite{Kcourse}}
\label{A(l,G),diamondsuit}
Let $G$ be  a connected graph and $G'$ be the graph obtained from $G$ by   a symmetric $K_{xy}$-operation. 
Then
$G \succ\succ ^{\alpha } G'$.
\es

Here is another useful $\succeq ^{\alpha }$-inequality.
\bs {\em \cite{Kcourse}}
\label{A(l,G),subgraph}
Let $G$ be  a connected graph. 
Then
$G - e \succ\succ ^{\alpha } G$ for every $e \in E(G)$.
\es

\subsection{Some results on relation $\succeq ^c $}
\label{c}
\bs {\em \cite{Kproceed,Kcourse}}
\label{discon,coef<=}
Let $G$ be a simple  graph, $x, y \in V(G)$, $x \ne y$, and $C_x$, $C_y$  different components of $G$ such that
$x \in V(C_x)$, $y \in V(C_y)$, $v(C_x) \ge 2$, and
 $v(C_y) \ge 2$. 
 Suppose that  $G'$ is obtained from $G$ by the
$\diamondsuit _{xy}$-operation.
Then $G \succ\succ ^c G'$.
\es

\bs {\em \cite{Kproceed,Kcourse}}
\label{coef<=}
Let $G$ be a simple connected graph and 
$G'$ be obtained from $G$ by a symmetric
$H_{xy}$-operation, where $H$ is connected.
Then 
\\[0.7ex]
$(a1)$ $c_s(G) > c_s(G')$ for every $s \in \{2, \ldots, n-2\}$
and
\\[0.7ex]
$(a2)$ 
$c_{n-1}(G) = c_{n-1}(G')$ if and only if $G$ is isomorphic to $G'$ or $x$ is a cut vertex in $G$.
\es

 Alternative proofs of {\bf \ref{coef<=}} as well as  more general results  and some other $\succeq ^c$-monotone operations on graphs are given below 
 (see {\bf \ref{coef<K-operation}} and 
{\bf \ref{coef<}} - {\bf \ref{A,F,aHb}}).
\\[1.5ex]
\indent
From {\bf \ref{coef<=}} we have, in particular:
\bs {\em \cite{Kcourse}}
\label{A,B,xHy-operation}
Let $A$, $B$, and $H$ be disjoint connected graphs,
$a \in V(A)$, $b \in V(B)$, and $x, y \in V(H)$, $x \ne y$.
Let $G$ be obtained from $A$, $B$, and $H$ by identifying $a$ with $x$ and $b$ with $y$ and let
$G'$ be obtained from $A$, $B$, and $H$ by identifying 
$a$ and $b$ with $x$.
Suppose that the two-pole $xHy$ is symmetric.
Then $c_s(G) > c_s(G')$ for every $s \in \{2, \ldots, n-2\}$
and
$c_{n-1}(G) = c_{n-1}(G')$.
\es

Since a close-door operation is a particular case of a 
symmetric hammock-operation, we have from 
{\bf \ref{coef<=}}:
\bs {\em \cite{Kcourse}}
\label{door}
Let $G$ be a connected graph and $G'$ be obtained from $G$ by a  close-door operation.
Let $s \in \{2, \ldots v(G) - 1\}$.
Then 
$c_s(G) > c_s(G')$ and
$c_s(G) = c_s(G')$ if and only if  $G$ and $G'$ are isomorphic.
\es

From {\bf \ref{T>G}} and 
{\bf \ref{door}}
we have, in particular:
\bs {\em \cite{Kproceed,Kcourse}} 
\label{T>^cG}
Let $G$ be a connected graph and let $\ddot{G}$ be the graph obtained from $G$ by adding $e(G) - v(G) +1$ isolated vertices. Then for every spanning tree $T$ of $G$ there exists a tree $D$ with $e(G)$ edges such that $D$ contains $T$ and $D\succeq ^c \ddot{G}$. Moreover,
if $e(G) \ge v(G)$, then $D\succ\succ ^c \ddot{G}$.
\es

From {\bf \ref{GtoTreshGraph}} and {\bf \ref{coef<=}} we have:
\bs {\em \cite{Kproceed,Kcourse}}
\label{threshcoef<=}
For every graph  $G$ in ${\cal C}_n^m$ there exists a threshold graph $F$ in ${\cal C}_n^m$ such  that
$G \succeq ^c F$.
\es

In the next theorem we will use the notions of an extreme threshold graph $F_n^m$ and the corresponding set of graphs ${\cal H}_n^m$ defined in 
{\bf \ref{extreme-threshold-graph}} and {\bf \ref{F_n^m}}.
\bs {\em \cite{Kproceed,Kcourse}}
\label{threshcoef<=some (m,n)}
Let $G \in {\cal C}_n^m$ and $G \ne F_n^m$.
\\[1.7ex]
$(a0)$
If $m = n - 1 \ge 3$, then
$c_s(G) > c_s(F_n^{n-1})$ for every 
$s \in \{2, \ldots , n-2\}$ and
\\[0.5ex]
$c_{n-1}(G) = c_{n-1}(F_n^{n-1}) = n$.
\\[1.7ex]
$(a1)$
If $m = n \ge 3$, then
$c_s(G) > c_s(F_n^n)$ for every 
$s \in \{2, \ldots , n-2\}$ and
\\[0.3ex]
$c_{n-1}(G) = c_{n-1}(H) = 3n$ for every 
$H \in {\cal H}_n^n$, and so $c_{n-1}(G) = c_{n-1}(F_n^n)$.
\\[1.7ex]
$(a2)$ If $n \ge 4$ and $m = n+1$, then
$c_s(G) > c_s(F_n^{n+1})$ for every 
$s \in \{2, \ldots , n-2\}$ and
\\[0.3ex]
$c_{n-1}(G) = c_{n-1}(H) = 8n$ for every 
$H \in {\cal H}_n^{n+1}$, and so $c_{n-1}(G) = c_{n-1}(F_n^{n+1})$.
\\[1.5ex]
$(a3)$ If $n \ge 5$ and $m = n+2$, then 
$c_s(G) > c_s(F_n^{n+2})$ for every 
$s \in \{2, \ldots , n-2\}$ and
\\[0.3ex]
$c_{n-1}(G) = c_{n-1}(H) = 16n$ for every 
$H \in {\cal H}_n^{n+2}$, and so $c_{n-1}(G) = c_{n-1}(F_n^{n+2})$.
\\[1.7ex]
$(a4)$ If $n \ge 6$ and $n+2 \le m \le 2n - 4$, 
then 
$c_s(G) > c_s(F_n^m)$ for every 
$s \in \{2, \ldots , n-2\}$
\\[0.3ex]
and
$c_{n-1}(G) = c_{n-1}(H)$ for every 
$H \in {\cal H}_n^m$, and so $c_{n-1}(G) = c_{n-1}(F_n^m)$.
\\[1.7ex]
$(a5)$ If  $m = 2n - 3$, 
then for every $n \ge 6$ there exists  $G \in {\cal C}_n^m$ such that $G \not \succeq ^c F_n^m$.
\es

Claims $(a0)$ - $(a4)$ in {\bf \ref{threshcoef<=some (m,n)}} follow basically from {\bf \ref{threshcoef<=}}.
\\[1ex]
\indent
Recall that an $P_{xy}$-operation is an $H_{xy}$-operation, where $xHy$ is an $xy$-path. This operations was one of the $\succeq ^{\phi }$-increasing operation on trees in \cite{Kcmpr}.
\\[1ex]
\indent
It is easy to see that any tree $T_n$ can be transformed to a path $P_n$ by a series of reverse $P_{xy}$-operations.
Therefore we have from  
{\bf \ref{coef<=}} and {\bf \ref{T>^cG}}:
\indent
\bs {\em \cite{Kproceed,Kcourse}}
\label{Psucc^cG}
Suppose that $G \in {\cal G}_n^m$, $G$ is connected, and graph $\ddot{G}$ is obtained from $G$ by adding $m - n +1$ isolated vertices.
If $G$ is not a path, then 
$P^m \succ\succ ^c \ddot{G}$.
\es

As we have mentioned before, theorem 
{\bf \ref{Psucc^cG}} on the poset 
$({\cal G}_n^m,  \succeq ^c)$ is similar to theorem 
{\bf \ref{P^m>G^m}} on the poset  
$({\cal G}_n^m,  \succeq ^{\phi })$.
\\[1ex]
\indent
Now we want to demonstrate  another  proof of 
{\bf \ref{Psucc^cG}}   that uses {\bf \ref{recursion}} and that is  similar to the proof of {\bf \ref{P^m>G^m}} in 
\cite{Kcmpr} using {\bf \ref{recursionPhi}}.
\\[1ex]
\indent
We need the following claim that can be easily proved  by induction using {\bf \ref{recursion}}.
\bs {\em \cite{Kcourse}}
\label{A1cpB1>A2cpB2} 
Let $i \in \{1,2\}$ and $(A_i, B_i)$ be a pair of disjoint paths.
Suppose that
\\[0.7ex]
$(h1)$ $v(A_1\cup B_1) = v(A_2\cup B_2) = n$,
\\[0.7ex]
$(h2)$ $v(A_i) \le v(B_i)$ for every $i \in \{1,2\}$, and
\\[0.7ex]
$(h3)$ $v(A_1) > v(A_2)$.
\\[0.7ex]
\indent
Then 
$c_s((A_1 \cup B_1) > c_s(A_2\cup B_2)$
for every $s \in \{2, \ldots , n - 1\}$. 
\es

Notice that {\bf \ref{A1cpB1>A2cpB2}} is also a simple particular case of {\bf \ref{A,B,aHb}} below.
\\[1ex]
\indent
We also need the following result interesting in itself.
Let $mx(G)$ denote the number edges of a  component  with the maximum number of edges.
\bs {\em \cite{Kcourse,Kcmpr}}
\label{Tr-e}
Let $P$ be a path, $T$  a tree, and $e(P) = e(T)$. Then there exists a bijection 
$\varepsilon : E(T) \to E((P)$ such that
$mx(T - u) \ge mx(P -  \varepsilon (u))$.
\es

Now we are ready to prove  
{\bf \ref{Psucc^cG}}. Obviously,  because of {\bf \ref{T>^cG}} theorem 
{\bf \ref{Psucc^cG}} follows from the  theorem for trees below.
\bs {\em \cite{Kproceed,Kcourse}}
\label{succ^c-upbounds}
Suppose that 
$G$ is a tree with $n$ vertices,
$G$ is not a path,
and $n \ge 4$.
Then $P_n \succ\succ ^c G$.
\es

{\bf Proof} 
~(uses
{\bf \ref{recursion}}, 
{\bf \ref{A1cpB1>A2cpB2}},  
and 
{\bf \ref{Tr-e}}).
We prove  our claim  by induction on  
$m = n-1$.
Recall that for two trees $Q$ and $R$ of $m$ edges,
$Q \succ\succ  ^c R$ if $c_s(Q) > c_s(R)$ for 
$2 \le s \le m-1$.

Suppose that $m = 3$. Then 
${\cal T}^{3} = \{P^3, S^3\}$ and $c_2(P^3) = 10 > 9 = c_3(S^3)$. Therefore our claim is true for $m = 3$.
We assume that our claim is true for $e(P) < m$ and prove that it is also true for $e(P) = m \ge 5$.
Let $s \in \{2, \ldots , m-1\}$. 
By {\bf \ref{Tr-e}},  there exists a bijection 
$\varepsilon : E(T) \to E((P)$ such that
$mx(T - u) \ge mx(P -  \varepsilon (u))$.
By {\bf \ref{recursion}}, 
\bqn
\label{recT}
(m - s) c_s(T) = \sum \{c_s(T - u): u  \in E(T)
\eqn
and
\bqn
\label{recP}
(m - s) c_s(P) = \sum \{c_s(P - \varepsilon (u)): u  \in E(T).
\eqn

Let $A_u$ and $B_u$ be the two components of 
$T - u$ and $e(A_u) \le e(B_u)$. 
Similarly, let $A'_u$ and $B'_u$ be the two components of 
$P - \varepsilon (u)$  and $e(A'_u) \le e(B'_u)$.
Then $e(B_u) \ge e(B'_u)$.
If $A$ and $B$ are disjoint graphs, then
\bqn
\label{c(AcupB)}
c_s(A \cup B) = \sum \{c_i(A) c_j(B): i + j  = s\}.
\eqn
Therefore
\\[1.5ex]
\indent
$c_s(T - u) = c_s(A_u \cup B_u) =\sum \{c_i(A_u) c_j(B_u): i + j  = s\}$. 
\\[1.5ex]
Let $A''_u$ and $B''_u$ be two disjoint paths such that
$e(A''_u) = e(A_u)$ and $e(B''_u) = e(B_u)$.
Since $e(A_u) < m$ and $e(B_u) < m$ and since $A_u$ and $B_u$ are trees, we have by the induction hypothesis:
\\[1.5ex]
\indent
$c_i(A''_u) \ge c_i(A_u)$ and $c_jB''_u) \ge c_j(B_u)$.
\\[1.5ex]
 Therefore by (\ref{c(AcupB)}),
\bqn
\label{c(A''cupB''>AcupB)}
c_s(A''_u \cup B''_u) \ge c_s(A_u \cup B_u) = c_s(T - u).
\eqn
Since $e(A''_u) = e(A_u)$  
and  $e(B_u) \ge e(B'_u)$, we have  by 
{\bf \ref{c(A''cupB''>AcupB)}} and
{\bf \ref{A1cpB1>A2cpB2}}:
\\[1.5ex]
\indent
$c(s(P - \varepsilon (u)) = c_s(A'_u \cup B'_u) \ge c_s(A''_u \cup B''_u) = c_s(T - u)$.
\\[1.5ex]
Since $T$ is not a path, there is $u \in E(T)$ such that
either $A_u$ or $B_u$ is not a path. 
Let $D_u$ 
be one of  $A_u$, $B_u$ which is not a path.
Then by the induction hypothesis,
$c_i(D''_u) > c_i(D_u)$ 
for some $i \le s$.
Now by (\ref{recT}) and (\ref{recP}) we have:
$c_s(P) > c_s(T)$.
\ep
\\[1.7ex]
\indent
The next theorem on $\succ\succ ^c$ is similar to theorem 
{\bf \ref{O^m>G^m}} on $\succ\succ ^\phi $ in \cite{Kcmpr}.
\bs {\em \cite{Kproceed,Kcourse}}
\label{Osucc^cG}
Suppose that $G \in {\cal G}_n^m$, $G$ is edge 2-connected,  and graph $\ddot{G}$ is obtained from $G$ by adding $m - n +1$ isolated vertices.
If $G$ is not a cycle, then 
$O^m \succ\succ ^c G$.
\es

{\bf Proof} ~(uses 
{\bf \ref{recursion}}
and
{\bf \ref{Psucc^cG}}). 
We prove that $O_n \succ\succ  ^c G$.
Since $e(O_n) = e(G)$, we can assume that 
$E(O_n) = E(G)$. Let $s \in \{2, \ldots , n-1\}$.
By {\bf \ref{recursion}},
\\[1ex]
$(n - s) c_s(O_n) = \sum \{c_s(O_n - u): u \in E(C_n)\}$ and
$(n - s) c_s(G) = \sum \{c_s(G - u): u \in E(G)\}$.
\\[1ex]
Obviously, $O_n - u$ is a path with $n$ vertices for every 
$u \in E(O_n)$.
Since $G$ is  edge 2-connected, 
$G - u$ is  connected for every $u \in E(G)$.
Therefore by {\bf \ref{Psucc^cG}}, 
$c_s(O_n - u) \ge  c_s(G - u)$ for every $u \in E(G)$.
Since $G$ is not a cycle, $G - t$ is not a path for some $t \in E(G)$. Hence by {\bf \ref{Psucc^cG}}, 
$c_s(O_n - t) >  c_s(G - t)$. 
Thus from  the above recursions we have:
$c_s(O_n) >  c_s(G)$. 
\ep 
\bs {\em \cite{Kproceed,Kcourse}}
\label{unicyclicG}
Let $O, G \in {\cal C}_n^n$, where  $O$ is a cycle. 
If $G$ is not a cycle, then $O \succ\succ ^c G$.
 \es

\subsection
{More  on Laplacian  posets of graphs}
\label{newLposets}

\indent

It turns out that a matching $M$ with $m$ edges is not only the $\succ\succ ^\phi $-maximum but also the $\succ\succ ^c $-maximum in $G^m$.
 \bs
\label{matching}
Let $M $ be a graph-matching,
$G$ a  simple graph not isomorphic to $M$, and 
$e(G) = e(M)$.
Then 
$M \succ\succ  ^{c, \phi} G$.
\es

{\bf Proof}~ (uses 
{\bf \ref{coefficients}}, 
{\bf \ref{recursion}}, 
{\bf \ref{propertyPhi}}, and
{\bf \ref{Msuccsucc^phi...}}).
Let $e(G) = m$.
By {\bf \ref{Msuccsucc^phi...}}, 
$M \succ\succ  ^{\phi} G$. 
Thus, it suffices  to prove 
that $c_s(M) > c_s(G)$ for every 
$s \in \{2, \ldots , m\}$.
We prove our claim by induction on $m$. 
By {\bf \ref{propertyPhi}}, $\Phi (\lambda , H)$ does not depend on the number of isolated vertices of a graph $H$.
\\[1ex]
\indent
Suppose that  $m = 2$. Then $G$ is a 3-vertex path plus an isolated vertex. Now
$\Phi (\lambda , M) = (\lambda - 2)^2 = 
\lambda ^2 - 4\lambda + 4$ and
$\Phi (\lambda , G) = (\lambda - 3)(\lambda - 1) = 
\lambda ^2 - 4\lambda + 3$. 
It follows that for $m = 2$ our claim is true.
Now we assume that our claim is true for $e(G) = m - 1$ and we prove that it is also true for $e(G) = m \ge 3$.
Since $e(M) = e(G)$, there is a bijection 
from $E(M)$ to $E(G)$. We can assume that this bijection is the identity, i.e. that $E(M) = E(G)$.
By {\bf \ref{recursion}},
\\[1ex]
$(m - s) c_s(M) =  \sum \{c_s(M - u): u \in E(G)\}$ and
$(m - s) c_s(G) = \sum \{c_s(G -u): u \in E(G)\}$
\\[0.7ex]
for $2 \le s \le m - 1$.
\\[1ex]
\indent
Suppose that $2 \le s \le m - 1$.
Obviously,   $e(M - u) = e(G - u) = m - 1$ and 
$M - u$ consists of a matching $M_u$ plus two isolated vertices,  and so $e(M_u) = e(M - u)= m - 1$. 
Since $\Phi (\lambda , M - u)$ does not depend on the number of isolated vertices of a graph $M - u$,
each $c_s(M_u) = c_s(M - u)$.
By the induction hypothesis,
$c_s(M - u) \ge c_s(G - u)$ for every $u \in E$.
Since $G$ is not a matching, $G$ has an edge $t$ such that the graph obtained from $G - t$ by removing the isolated vertices  is not a matching. Then by the induction hypothesis, $c_s(M - t) > c_s(G -t)$. Thus,  our claim for
 $m$ follows from the above recursions.
\\[1ex]
\indent
Finally, suppose that $s = m$. 
Then  by {\bf \ref{coefficients}}, $c_m(M) = \gamma (M) = 2^m$.
If $G$ is not a forest, then $c_m(G) = 0$.
So we assume that $G$ is a forest, and so by {\bf \ref{coefficients}}, $c_m(G) = \gamma (G)$. Since $G$ is not a matching, it has a component $C$ with $r = e(C) \ge 2$.
 Consider in $M$ the the subgraph
$N$ induced by the edge subset $E(C)$. Then 
$N$ is a  matching with $r$ edges and $\gamma (N) = 2^r$.
If $r = m$, then  $\gamma (G) = m+1$, and so
$\gamma (M) = 2^m > m+1 = \gamma (G)$.
So we assume that $r < m$. 
Then $\gamma (M) = \gamma (N)~ \gamma (M - N)$ and
$\gamma (G) = \gamma (C) ~\gamma (G - C)$.
Since $r \ge 2$, $\gamma (N) = 2^r > r+1 = \gamma (C)$ and by the induction hypothesis,
$\gamma (M - N) \ge \gamma (G - C)$. Therefore
$c_m(M) = \gamma (M) > \gamma (G)$.
\ep
\\[1.5ex]
\indent
Using {\bf \ref{recursion}}, {\bf \ref{Msuccsucc^phi...}}, and
{\bf \ref{matching}}
we can also  prove the following stronger result.
\bs
\label{star}
Let $G$ be a simple graph with $m$ edges  
not isomorphic to $mP^1$, $P^2 + (m-2)P^1$, and  
$S^m$.
Then 
$mP^1 \succ\succ ^{c, \phi} P^2 + (m-2)P^1 \succ\succ ^{c, \phi} G \succ\succ ^{c, \phi} S^m$.
\es

Now we need to recall the definition of a ${\cal K}_{xy}$-operation on a graph $G$.
Let $x, y \in V(G)$ and ${\cal K}$ be an induced  subgraph of $G$ containing $x$  and $y$. 
Let
$X = N_x(G) \setminus  (V({\cal K})\cup N_y(G))$,
$Y = N_y(G) \setminus  (V({\cal K})\cup N_x(G))$,
$[x,X] = \{xv: v \in X\}$, and $[y,X] = \{yv: v \in X\}$, and so $[x,X] \subseteq E(G)$ and $[y,X] \cap  E(G) = \emptyset $.
Let
$G' = {\cal K}_{xy}(G) =  (G - [x,X]) \cup [y,X]$, and so 
$[y,X] \subseteq E(G')$ and $[x,X] \cap E(G') = \emptyset $.
We say that  ${\cal K}_{xy}(G)$ {\em is obtained from $G$ by the ${\cal K}_{xy}$-operation}.
We call ${\cal K}_{xy}$-operation on $G$ {\em symmetric} if
$G -([x,X]) \cup [y,Y])$ has an automorphism 
$\alpha : V(G) \to V(G)$ such that
$\alpha (x) = y$, $\alpha (y) = x$,
$\alpha [{\cal K}] = {\cal K}$,
$\alpha (v) = v$ for every $v \in X \cup Y$, and so 
$\alpha [G - ( {\cal K} \cup X \cup Y)] = G - ( {\cal K} \cup X \cup Y)$.
\\[1.5ex]
\indent
Let, as above,  $\delta _s[G] = \sum \{d(v,G)^s: v \in V(G)\}$.
It is easy to prove the following generalization of 
{\bf \ref{delta(G),H-operation}}.

\bs 
{\em \cite{Kcourse,Kproceed}}
\label{delta(G),K-operation}
Let $G, G' \in {\cal G}^m_n$ and 
$G'$ be obtained from $G$ by a symmetric  ${\cal K}_{xy}$-operation.
Let $f : \mathbb{R}^n \to \mathbb{R}$ be a convex symmetric function and $f [G] = f (\{d(v,G): v \in V(G)\})$.
Then $f[G] \le f[G']$, and so, in particular,
$\delta _s [G] \le \delta _s[G']$ for every positive integer 
$s$.
\es

\bs
\label{coef<K-operation}
Let $G$ be a graph and 
$G'$  be obtained from $G$ by a symmetric ${\cal K}_{xy}$-operation.
Then $c_s(G) \ge c_s(G')$ for every 
$s \in \{0, \ldots , v(G) - 1\}$.
\es

{\bf Proof} ~ (uses {\bf \ref{somecoefficients}}, {\bf \ref{recursion}}, and 
{\bf \ref{properties-of-Kxy}} $(a3)$).
If $\alpha : V(G) \to V(G)$ is an automorphism of a graph $G$ and $F$ is  a subgraph of $G$, then let
$\alpha [F]$ denote the image of $F$ under the automorphism $\alpha $. In particular, if $e = pq \in E(G)$, then  $\alpha [e] = \alpha (p)\alpha (q)$ and if 
$A \subseteq E(G)$, then let $\alpha [A] = \{\alpha [a]: a \in A\}$.
By {\bf \ref{properties-of-Kxy}} $(a3)$, we can assume that
$N_x(G) \cap N_y(G) \subseteq V({\cal K})$.
Now, since $G'$ is obtained  from $G$ by 
the ${\cal K}_{xy}$-operation, we have: 
\\[1ex]
\indent
$G' =  {\cal K}_{xy}(G) = (G - [x,X]) \cup [y,X]$, 
where 
$X = N_x(G) \setminus  V({\cal K})$,
$Y = N_y(G) \setminus  V({\cal K})$,
\\[0.5ex]
and so $[x,X] \subseteq E(G)$ and 
$[y,X] \cap E(G) = \emptyset $.
\\[1ex]
\indent
Since the ${\cal K}_{xy}$-operation in $G$ is symmetric, 
$G \setminus ([x,X] \cup [y,Y])$ has an 
 automorphism 
$\alpha : V(G) \to V(G)$ such that
$\alpha (x) = y$, $\alpha (y) = x$,
$\alpha [{\cal K}] = {\cal K}$, and 
$\alpha (v) = v$ for every $v \in X \cup Y$, and so 
$\alpha [G - ( {\cal K} \cup X \cup Y)] = G - ( {\cal K} \cup X \cup Y)$ and  the ${\cal K}_{xy}$-operation is $\alpha $-symmetric. Let $E^-(G) = E(G \setminus ([x,X] \cup [y,Y]))$.
\\[1ex]
\indent
Let $A(G) = \{e \in E^-(G): \alpha [e] = e\}$ and
$B(G) = \{e \in E^-(G): \alpha [e] \ne e\}$, 
\\[1ex]
and so 
$A(G) \cap B(G) = \emptyset $ and 
$E(G) = A(G) \cup B(G) \cup [x,X] \cup [y,Y]$.
\\[1ex]
\indent
By {\bf \ref{somecoefficients}}, we can assume that 
$s \in \{2, \ldots , v(G) - 1\}$.
We prove our claim by induction on $e(G) = m$.
If $m = 0$, then our claim is obviously true.
Suppose that our claim is true for every graph $G$ with 
$e(G) = m-1$. We prove that it is also true for every graph 
$G$ with $e(G) = m \ge 1$.
Let $e_v = xv$ and $e'_v = yv$ for $v \in X$ and let
\\[1ex]
\indent
$  E[X] = [x,X] = \{e_v: v \in X\}$ and 
$  E'[X] = [y,X] = \{e'_v: v \in X\}$.
\\[1ex]
Then $G' = (G - E[X]) \cup E'[X]$.
Obviously, $E(G) - E[X] = E(G') - E'[X]$.
Let $\varepsilon (u) = u$ if $u \in E(G) - E[X]$ and
$\varepsilon (e_v) = e'_v$ if $e_v \in E[X]$ 
(and so $v \in X$).
Then $\varepsilon : E(G) \to E(G')$ is a bijection.
By {\bf \ref{recursion}}, we have the following recursions for 
$s \le m - 1$:
\bqn
\label{recarsionc(G)}
(m - s) c_s(G) =  \sum \{c_s(G - u): u \in E(G)\}
\eqn
and
\bqn
\label{recarsionc(G')}
(m - s) c_s(G') = \sum \{c_s(G' - \varepsilon (u): u \in E(G)\}
\eqn

Suppose first that $s = m$.
If $G$ is not a forest, then $G'$ is also not a forest.
Therefore we have: $c_m(G) = c_m(G')  = 0$.
So let $G$ be a forest. Then 
$c_m(G) = \gamma (G)$ and $c_m(G') = \gamma (G')$.
In this case it is easy to show that $c_m(G) \ge c_m(G')$.
\\[1ex]
\indent 
Now suppose that $s \le m -1$.
\\[1.5ex]
\indent
{\sc Claim 1.}
{\em If $u \in A(G)$, 
$($i.e.  $\alpha [u] = u$$)$, 
then $c_s(G - u) \ge c_s(G' - u)$.}
\\[1ex]
\indent
{\em Proof.}
If $u \in E(G)  - E({\cal K})$, then
$G' - \varepsilon (u) = (G - u) - E[X]) \cup E'[X] = {\cal K}_{xy}(G - u)$.
\\[1ex]
If $u \in E({\cal K})$, then 
$G' - u = (G - u) - E[X]) \cup E'[X] = {\cal K}'_{xy}(G - u)$,
\\[1ex] 
where ${\cal K}' = {\cal K} - u$. Since $\alpha [u] = u$,
clearly  ${\cal K}'_{xy}$-operation in $G - u$ is 
$\alpha $-symmetric.
In both cases,
since $e(G - u) < e(G)$, we have by the inductive hypothesis:
\\[0.7ex]
$c_s(G - u) \ge c_s(G' - u)$.
\epcl
\\[1.5ex]
\indent
{\sc Claim 2.}
{\em If $u \in [x,X] \cup [y,Y]$, 
then $c_s(G - u) \ge c_s(G' - \varepsilon (u))$.}
\\[1ex]
\indent
{\em Proof.} 
Suppose that  $u = [x,X])$. 
Then $u = xv$ for some $v \in X$.
Obviously, 
$G' - \varepsilon [u] = (G - u) - E[X-v]) \cup E'[X-v] = 
{\cal K'}_{xy}(G - u)$, where ${\cal K'}$ is the subgraph of $G$ induced by ${\cal K} \cup v$. 
Since $\alpha (z) = z$ for every $z \in X \cup Y$and ${\cal K}_{xy}$-operation in 
$G$ is $\alpha $-symmetric,
clearly  ${\cal K}'_{xy}$-operation in $G - u$ is also 
$\alpha $-symmetric. 
Since  $e(G - u) < e(G)$, we have by the inductive hypothesis: $c_s(G - u) \ge c_s(G' - u)$.
\\[0.7ex]
\indent
Similar arguments show that our claim is also true for 
$u \in [y,Y]$.
\epcl
\\[1.5ex]
\indent
By the recursions (\ref{recarsionc(G)}) and  
(\ref{recarsionc(G')})
and by {\sc Claims} 1 and 2,
it is sufficient to prove 
the following inequality for every $u \in B(G)$   
(and so $\alpha [u] \ne u$):
\bqn
\label{c(G-u)+c(G-alpha[u])>}
c_s(G - u) + c_s(G - \alpha [u]) \ge  
c_s(G' - u) + c_s(G' - \alpha [u]).
\eqn

This inequality is a particular case of the following claim.
\\[1.5ex]
\indent
{\sc Claim 3.}
If $Z \subseteq B(G)$, then 
{\em  
$c_s(G - Z) + c_s(G - \alpha [Z]) \ge  
c_s(G' - Z) + c'_s(G - \alpha [Z])$.}
\\[1ex]
\indent
{\em Proof.}
Suppose, on the contrary, that our claim is not true.
Let $(G, Z)$ be an $(|E|, |B - Z|)$-lexicographically 
smallest (or lg-smallest) counterexample to our claim, where 
$E = E(G)$ and $B = B(G)$.  If $Z = B$, then our claim is obviously, true. Therefore such a counterexample exists. Let
\\[1.5ex]
$\sigma (G, Z) = c_s(G - Z) + c_s(G - \alpha [Z])$ and 
$\sigma (G', Z) = c_s(G' - Z) + c_s(G' -\alpha [Z])$.
\\[1.5ex]
By the recursions (\ref{recarsionc(G)}) and  
(\ref{recarsionc(G')})
for $c_s(G)$ and $c_s(G')$, we have:
\\[1.5ex]
\indent
$(m - s) \sigma (G, Z) =
\sum \{ \sigma (G - u, Z): u \in [x,X] \cup [y,Y]\} +
\\[1ex]
\indent
\sum \{ \sigma (G - u, Z): u \in A(G)\} +
\sum \{ \sigma (G - u,Z): u \in B - Z\}.$
\\[1.5ex]
For $\sigma (G', Z)$ we have a similar formula obtained from the above formula by replacing $G$ by $G'$.
Let 
$\Delta (G, Z) = 
\sigma (G, Z) - \sigma (G', Z)$.
Since $(G,Z)$ is a counterexample, we have:
$\Delta (G, Z) < 0$.
\\[1ex]
\indent
Our goal is to get a contradiction by showing  that $\Delta (G, Z) \ge 0$.
From the above relations we have:
\\[1ex]
\indent
$(m - s) \Delta (G, Z) =
\sum \{ \Delta (G - u, Z): u \in [x,X] \cup [y,Y\} +
\\[1ex]
\indent
\sum \{ \Delta (G - u, Z): u \in A(G)\} +
\sum \{ \Delta (G - u,Z): u \in B - Z\}.$
\\[1.5ex]
\indent
Suppose that $u = zv \in [x,X] \cup [y,Y]$.
Then arguments similar to those in the proof of 
{\sc Claim 2}, show that  
$G' - \varepsilon (u) = {\cal K}'_{xy}(G - u)$, where
${\cal K}'$ is the subgraph of $G$ induced by 
${\cal K} \cup v$ and ${\cal K}'_{xy}$-operation in 
$G - u$ is $\alpha $-symmetric.
Since $(|E - u|, |B - Z|)$ is  lg-smaller than $(|E|, |B - Z|)$,
clearly $(G - u, Z)$ is not a counterexample. Therefore 
\\[1ex]
\indent
$\sum \{ \Delta (G - u, Z): u \in [x,X] \cup [y,Y]\} \ge 0$. 
\\[1.5ex]
\indent
Suppose that $u \in A(G)$.
If $u \not \in E({\cal K})$, then 
$G' - \varepsilon (u) = {\cal K}_{xy}(G - u)$.
Since $(|E - u|, |B - Z|)$ is  lg-smaller than $(|E|, |B - Z|)$, obviously,
$(G - u, Z)$ is not a counterexample. Therefore 
$ \Delta (G - u, Z) \ge 0$. 
If $u  \in E({\cal K})$, then 
$G' - \varepsilon (u) = {\cal K}'_{xy}(G - u)$, where 
${\cal K}' = {\cal K} - u$, and
the ${\cal K}'_{xy}$-operation is $\alpha $-symmetric in $G - u$. Since $(|E - u|, |Z|)$ 
is lg-smaller than $(|E|, |B - Z|)$, again
$(G - u, Z)$ is not a counterexample. Therefore 
$ \Delta (G - u, Z) \ge 0$. 
Thus, 
\\[1ex]
\indent
$\sum \{ \Delta (G - u, Z): u \in A(G)\} \ge 0$.
\\[1.5ex]
\indent
Finally, suppose that $u \in B - Z$. 
Since $(|E|, |B - u - Z|)$ is  lg-smaller than 
\\
$(|E|, |B - Z|)$,
again $(G, Z \cup u)$ is not a counterexample. Therefore
\\[0.7ex] 
\indent
$ \Delta (G - u, Z) =  \Delta (G , Z \cup u) \ge 0$, and so  
 $\sum \{ \Delta (G - u, Z): u \in B - Z\} \ge 0$.
 \\[1ex]
\indent
Thus, from the above recursion for $\Delta (G, Z)$ we have: $\Delta (G, Z) \ge 0$.
\epcl
\\[1.5ex]
\indent
Obviously, inequality (\ref{c(G-u)+c(G-alpha[u])>}) is 
a particular case of {\sc Claim 3}, when  
$|B| = 1$.
\ep
\\[1.5ex]
\indent 
Using {\bf \ref{properties-of-Kxy}}, we have  from {\bf \ref{coef<K-operation}}: 
\bs 
\label{t(G),K-operation}
Let $G, F \in {\cal G}^m_n$ and 
$F$ be obtained from $G$ by a symmetric 
${\cal K}_{xy}$-operation. 
Then $t(G) \ge t(F)$
and $t([G]^c) \ge t([F]^c)$.
\es 

Now from {\bf \ref{t(G),K-operation}} we have the following strengthening of {\bf \ref{>t,D-operation}}:
\bs 
\label{>t,K-operation}
Let $G, F \in {\cal G}_n^m$ and 
$F$ be obtained from $G$ by a symmetric 
${\cal K}_{xy}$-operation.
Then 
$G\succeq ^{\tau } F$ and $[G]^c \succeq ^{\tau } [F]^c$.
\es

The following is a generalization of 
{\bf \ref{Gsucceq^pH(G)}} for 
a ${\cal K}_{xy}$-operation.
\bs
\label{Gsucceq^pKxy(G)} 
Let $G$ be a graph and $G'$ be obtained from $G$ by a symmetric  ${\cal K}_{xy}$-operation.
Then $G \succeq ^p G'$. 
\es

{\bf Proof}~ (uses {\bf 
\ref{recursionPhi}}, {\bf \ref{coef<K-operation}}, and
{\bf \ref{t(G),K-operation}}).
Let $v(G) = n$.
Our goal is to prove that 
$\Phi (\lambda , G) \ge  \Phi (\lambda , G')$ for
$\lambda \ge n$.
We prove by induction on $e(G) = m$. Our claim is obviously true for $e(G) = 1$. Suppose that our claim is true for $e(G) = m -1$. We will prove that it is also true for 
$e(G) = m \ge 2$.
By {\bf \ref{properties-of-Kxy}} $(a3)$, we can assume that
$N_x(G) \cap N_y(G) \subseteq V({\cal K})$.
Now, since $G'$ is obtained  from $G$ by 
the ${\cal K}_{xy}$-operation, we have: 
$G' =  G - [x,X] \cup [y,X]$, where 
$X = N_x(G) \setminus  V({\cal K})$,
$Y = N_y(G) \setminus  V({\cal K})$,
$[x,X] = \{xv: v \in X\}$ and $[y,X] = \{yv: v \in X\}$, and so $[x,X] \subseteq E(G)$ and 
$[y,X] \cap E(G) = \emptyset $.
Let $e_v = xv$ and $e'_v = yv$ for $v \in X$ and let
\\[1ex]
\indent
$  E[X] = [x,X] = \{e_v: v \in X\}$ and 
$  E'[X] = [y,X] = \{e'_v: v \in X\}$.
\\[1.5ex]
Then $G' = (G - E[X]) \cup E'[X]$.
Obviously, $E(G) - E[X] = E(G') - E'[X]$.
Let $\varepsilon (u) = u$ 
\\[0.5ex]
if $u \in E(G) - E[X]$ and
$\varepsilon (e_v) = e'_v$ if $e_v \in E[X]$ (and so 
$v \in X$).
Then $\varepsilon : E(G) \to E(G')$ 
\\[0.5ex]
is a bijection.
By {\bf \ref{recursionPhi}}, 
\\[1.5ex]
\indent
$\Phi (\lambda , G) = \Phi (n , G) +
\sum \{\int_n^\lambda \Phi (t, G - u)dt : u \in E(G)\}$.
\\[1.5ex]
Let 
$\Delta (\lambda , G)  = 
\Phi (\lambda , G) - \Phi (\lambda , G')$ and 
$\Delta (\lambda , G - u)  = \Phi (\lambda , G - u) - \Phi (\lambda , G' - \varepsilon (u)).$
\\[1ex]
Then
\bqn
\label{difference}
\Delta (\lambda , G) = 
\Delta (n , G) +
\sum \{\int_n^\lambda (\Delta  (t, G - u))dt : u \in E(G)  \}.
\eqn
By  {\bf \ref{P(n+r)=}}, 
$\Phi (n , G) =   n ^ {m - n +2} t(K_n - E(G))$.
Therefore 
\\[1.5ex]
\indent 
$\Delta (n , G) =  \Phi (n , G) - \Phi (n , G') = 
n^ {m - n + 2} (t(K_n - E(G)) - t(K_n - E(G'))$.
\\[1.5ex]
By {\bf \ref{t(G),K-operation}}, 
$t(K_n - E(G)) - t(K_n - E(G')) \ge 0$.
Therefore 
\bqn
\label{difference-n}
\Delta (n , G) = \Phi (n , G) - \Phi (n , G') \ge 0.
\eqn
Now, using the induction hypothesis and the arguments similar to those in the proof of {\bf \ref{coef<K-operation}}, 
it can be shown that 
\bqn
\label{difference-t}
\sum \{\int_n^\lambda (\Delta (t, G - u))dt : u \in E(G)\} \ge 0.
\eqn
Thus,  our claim follows from  (\ref{difference}),
(\ref{difference-n}), and (\ref{difference-t}). 
\ep
\\[2.5ex]
\indent
From {\bf \ref{properties-of-Kxy}} $(a2)$ and 
{\bf \ref{Gsucceq^pKxy(G)}}, we have  
the following strengthening of {\bf \ref{>t,K-operation}}: 
\bs 
\label{succeq^p,K-operation}
Let $G, F \in {\cal G}^m_n$ and 
$F$ be obtained from $G$ by a symmetric 
${\cal K}_{xy}$-operation. 
Then $t(G) \ge t(F)$
and $[G]^c \succeq ^p [F]^c$.
\es 

From {\bf \ref{Gsucceq^pKxy(G)}}  we have, in particular,  the following strengthening of {\bf \ref{T>^bG}}.
\bs
\label{T>^pG}
Let $G$ be a simple connected graph and $\ddot{G}$ be the graph obtained from $G$ by adding $e(G) - v(G) +1$ isolated vertices. Then for every spanning tree $T$ of $G$ there exists a tree $D$ with $e(G)$ edges such that $T$ is a subgraph of $D$ and  $D\succeq ^p \ddot{G}$.
\es

We can also prove  that  under the assumption in 
{\bf \ref{T>^pG}}, if $e(G) - v(G) +1 > 0$ 
(i.e. if $\ddot{G} \ne G$), then $D\succ\succ ^p G$.
\\[1.5ex]
\indent
Another way to prove {\bf \ref{coef<K-operation}} is by fixing a spanning forest in $G$ and analyzing how it is transformed by the operation that brings $G$ to $G'$.
We demonstrate this approach by giving another proof of 
a particular case of  {\bf \ref{coef<K-operation}} when 
$G'$ is obtained from $G$ by 
a $\diamondsuit $-operation.
\bs
\label{coef<}
Let $G$ be a connected graph with $n$ vertices, $x$ and $y$ two distinct vertices in $G$.
Let 
$G'$ be obtained from $G$ by the $\diamondsuit _{xy}$-operation.
Then 
\\[1ex]
$(a1)$
$c_s(G) > c_s(G')$ for every $s \in \{2, \ldots , n-2\}$ and
\\[1ex]
$(a2)$
$ c_{n-1}(G) = c_{n-1}(G')$ if and only if $x$ is a cut vertex in $G$.
\es

{\bf Proof}~(uses 
{\bf \ref{coefficients}} and {\bf \ref{somecoefficients}}).
By {\bf \ref{somecoefficients}},
we can assume that $s \in \{2, \ldots , n-1\}$.
Let $Z = N_x(G) \cap N_y(G)$, $X = N_x(G) \setminus (Z \cup \{y\})$, and $Y = N_y(G) \setminus (Z \cup \{x\})$.
Let $e_v = xv$ and $e'_v = yv$ for $v \in X$ and let
\\[1ex]
\indent
$E_x = [x,X] = \{e_v: v \in X\}$ and 
$  E_y = [y,X] = \{e'_v: v \in X\}$.
%
\\[1ex]
Since $G'$ is obtained from $G$ by the 
$\diamondsuit _{xy}$-operation, 
\\[1ex]
\indent
$G' = (G - E_x) \cup E_y$
and $[\{x,y\}, Z] \subseteq E(G)$.
Obviously, $E(G) - E_x = E(G') - E_y$.
Let $\varepsilon (u) = u$ if $u \in E(G) - E_x$ and
$\varepsilon (e_v) = e'_v$ if $e_v \in E_x$ 
(and so $v \in X$).
Then $\varepsilon : E(G) \to E(G')$ is a bijection.
For $U \subseteq E(G)$, let
$\varepsilon [U] = \{\varepsilon (u): u \in U\}$.
\\[1ex]
\indent
For a subgraph $S$ of $G$, let
$\vartheta (S)$ be the subgraph of $G$ such that
$V(\vartheta (S)) = V(S)$ and 
$E(\vartheta (S)) = (E(S) \setminus  E_x)\cup
\varepsilon [E(S) \cap E_x]$.
Let ${\cal S}(H)$ denote the set of subgraphs of a graph $H$.
Obviously, $\vartheta $ is a bijection from 
${\cal S}(G)$ to ${\cal S}(G')$.
For  ${\cal A} \subseteq {\cal S}(G)$, let
$\vartheta [{\cal A}] = \{\vartheta (A): A \in {\cal A}\}$. 
\\[1ex]
\indent
Let $P'$ be a forest in $G'$ with at most two components each meeting $\{x,y\}$. Let $P'_x$ and $P'_y$ be the components of $P'$ containing $x$ and $y$, respectively, 
and so if $P'$ has one component, then $P'_x = P'_y$.
Let $P$ be the subgraph in $G$ such that 
$P' = \vartheta(P)$.
Clearly, $\{x,y\} \subseteq  V(P') = V(P)$ and  
$e(P') = e(P)$. 
Let
$\sigma (xz) = yz$ and $\sigma (yz) = xz$
for every $z \in Z$. 
Obviously,
$\sigma : [\{x,y\}, Z] \to [\{x,y\},Z]$ is a bijection.
For $S \subseteq [\{x,y\}, Z]$, let 
$\sigma [S] = \{\sigma (s): s \in S\}$.
Let $A = A(P) =  [\{x,y\}, Z] \cap E(P)$. 
Obviously, $A(P) = A(P')$.
Let 
$\bar{P} = (P -  A) \cup \sigma [A]$
and 
$\bar{P}' = (P' -  A) \cup \sigma [A]$.
Obviously, $\bar{P} = P$ and $\bar{P'} = P'$ if and only if 
$\sigma [A] = A$.
\\[1.5ex]
\indent
We need the following facts.
\\[1ex]
\indent
{\sc Claim 1.}
{\em Suppose that $P'$ is a tree and $P$ is not a tree.
Then $\bar{P}$ is a tree and $\bar{P}'$ is not a tree.
}
\\[1ex]
\indent
{\em Proof.}
Since $v(P) = v(P')$,  $e(P) = e(P')$, $P'$ is a tree, and
$P$ is not a tree, we have: $P$ has a cycle $C$.
If $C$ does not contain vertex $x$ or $C$ contains two edges from $E_x$
then  
$P'$ has a cycle, a contradiction. 
Therefore $C$ is the only cycle in $P$  and $C$ contains exactly one edge $xc$ with $c \in X$ and exactly one edge $xz$ with $z \in Z$.
Thus, $P$ has exactly  two components.
Now our claim follows.
\epcl
\\[1.5ex]
\indent
It is easy to see that the converse of {\sc Claim 1} may be not true. Namely, we have:
\\[1ex]
\indent
{\sc Claim 1'.}
{\em The following are equivalent:
\\[1ex]
$(a1)$ $x$ is not a cut vertex in $G$ and 
\\[1ex]
$(a2)$
there exists a spanning tree $P$ of $G$ such that $\bar{P}$ is a spanning tree and both $P'$ and $\bar{P}'$ are not trees.
}
\\[1ex]
\indent
{\sc Claim 2.}
{\em If $P'$ is a forest with exactly two components, then
$P$ is also a forest with exactly two components.
}
\\[1ex]
\indent
{\em Proof.}
Suppose, on the contrary, $P'$ is a forest with exactly two components but $P$ is not a forest with exactly two components. Then $P$ has at least three components. Then there is a component $Q$ of $P$ that avoids $\{x,y\}$.
Then $Q$ is also a component of $P'$.
However, $P'$ has exactly two components each containing exactly one vertex from $\{x,y\}$, 
a contradiction.
\epcl
\\[1.5ex]
\indent
Actually, the converse of {\sc Claim 2} is also true.
\\[1ex]
\indent
{\sc Claim 2'.}
{\em $P'$ is a forest with exactly two components if and only if $P$ is also a forest with exactly two components.
}
\\[1.5ex]
\indent
Suppose that both $P'$ and $P$ are forest with exactly two components.
Let 
$X^* = \{v \in X: xv \in E(P)\}$ and 
$Y^* = \{v \in Y: yv \in E(P)\}$.
Let, as above,  $P'_x$ and $P'_y$ be the two components of $P'$ such that
$x \in V(P'_x)$ and $y \in V(P'_y)$.
Let us remove from $P'_y$ the edges from $y$ to 
$X^*\cup Y^*$, denote by  
$Q'_y$ the component of the resulting forest containing vertex $y$, and put $R = P'_y - Q'_y$. 
Then  each component of the forest $R$ has exactly one vertex in $X^* \cup Y^*$. Let $R_x = R_x(P) = R_x(P')$ be the union of the components meeting $X^*$ and 
$R_y = R_y(P) = R_y(P')$ be the union of the components meeting $Y^*$.
Let 
\\[2ex]
\indent
$\Delta (P, P') =
(c_r(P) + c_r(\bar{P})) - (c_r(P') + c_r(\bar{P'}))$,
where $r = e(P)$.
\\[2.5ex]
\indent
{\sc Claim 3.}
{\em 
Suppose that both $P$ and $P'$ are forest with exactly two components. Then 
\\[1ex]
$\Delta (P, P') =  2 v(R_x)v(R_y) \ge 0$, and so
$\Delta (P, P') > 0 \Leftrightarrow  v(R_x) > 0~and~v(R_y) > 0$.
} 
\\[1ex]
\indent
{\em Proof.}
Obviously,
\\[1ex]
$c_r(P') = v(P'_x) (v(Q'_y) + v(R_x) + v(R_y))$,
$c_r(\bar{P'}) = v(Q'_y) (v(P'_x) + v(R_x) + v(R_y))$,
\\[1ex]
$c_r(P) = (v(P'_x) + v(R_x) ) (v(Q'_y) +  v(R_y))$, and
$c_r(\bar{P}) = (v(Q'_y) + v(R_x) ) (v(P'_x)+  v(R_y))$.
\\[1ex]
\indent
From the above formulas we have:
\\[1.5ex]
$\Delta (P, P') = 
([(v(P'_x) + v(R_x) ) (v(Q'_y) +  v(R_y))] +
[(v(Q'_y) + v(R_x) ) (v(P'_x)+  v(R_y))]) -
\\[1.5ex]
([v(P'_x) (v(Q'_y) + v(R_x) + v(R_y))] +
[v(Q'_y) (v(P'_x) + v(R_x) + v(R_y))] =
\\[1.5ex]
[v(P'_x) (v(Q'_y) +  v(R_y)+ v(R_x) (v(Q'_y) +  v(R_y))] -
[v(P'_x)(v(Q'_y) + v(R_y)) + v(P'_x)v(R_x)] +
\\[1.5ex]
[(v(Q'_y)(v(P'_x)+  v(R_y)) + v(R_x) (v(P'_x)+  v(R_y))] -
[v(Q'_y) (v(P'_x) + v(R_y))] + v(Q'_y)v(R_x)] =
\\[1.5ex]
 v(R_x) v(Q'_y) + v(R_x) v(R_y) - v(P'_x)v(R_x) +
 v(R_x) v(P'_x)+  v(R_x)v(R_y) - v(Q'_y)v(R_x) =
 \\[1.5ex]
 2 v(R_x)v(R_y).$
 Our claim follows.
\epcl
\\[1.5ex]
\indent
Given a spanning forest $F$ of $G'$, let 
$F_{xy}$ be the minimal subforest of $F$ containing $x$ and $y$ and such that each component of $F_{xy}$ is a component of $F$. Obviously, $F_{xy}$ has at most two components each meeting $\{x,y\}$.
Let ${\cal F}_s(G')$ denote the set of  spanning forests $F$  with $s$ edges in $G'$. 
 Let ${\cal P}(G')$ denote the set of all forests $P'$  in $G'$ having  at most two components each meeting 
 $\{x,y\}$ and let  
${\cal P}_r(G')$ denote the set of all forests in ${\cal P}(G')$ having  $r$ edges. 
Now from {\bf \ref{coefficients}} we have:
\\[0.7ex]
\indent
$c_s(G') = \{ \sum \{\sum \{ \gamma (F): F \in  {\cal F}_s(G'), F_{xy} = P'\}: P' \in {\cal P}(G')\}$.
\\[0.7ex]
Therefore
\bqn
\label{c(G')formula}
c_s(G') =
\sum \{\sum \{c_r(P') c_{s - r}(G - V(P')): P' \in {\cal P}_r(G')\}: 0 \le r = e(P') \le s\}.
\eqn
Similarly,
\\[1.5ex]
\indent
$c_s(G) =
\sum \{\sum \{ \gamma (F): F \in  {\cal F}_s(G), F_{xy} = Q\}: Q \in {\cal P}(G)\} =
\\[1.5ex]
\indent
\{ \sum \{\sum \{c_r(Q) c_{s - r} (G - V(Q)): Q \in {\cal P}_r(G)\}: 0 \le r = e(Q) \le s\}$.
\\[1.5ex]
Recall that  $P$ is the subgraph in $G$ such that 
$P'$ is obtained from $P$ by the $\diamondsuit _{xy}$-operation. 
\\[0.5ex]
\indent
Let $\widetilde{P}' = P$   if $P$ is a forest and 
$\widetilde{P}' = \bar{P}$ if $P$ is not a forest.
By {\sc Claim} 1 and {\sc Claim} 2, if $P$ is not a forest, then  both $P'$ and $\bar{P}$ are  trees.
 
Let
\bqn
\label{c'(G)formula}
c'_s(G) = 
\sum \{\sum \{c_r(\widetilde{P}') c_{s - r}(G - V(P'): P' \in {\cal P}_r(G')\}: 0 \le r = e(P') \le s \}.
\eqn
Obviously, $c_s(G) \ge c'_s(G)$. 
If $P'$ is a tree with $r$ edges , then $\widetilde{P}'$ is also a tree with $r$ edges and $c_r(P') = c_r(\widetilde{P}') = r+1$.
Therefore by {\sc Claim} 3, we have from (\ref{c(G')formula}) and (\ref{c'(G)formula}): 
$c'_s(G) \ge c_s(G')$.
Thus, $c_s(G) \ge c'_s(G) \ge c_s(G')$, and so claim $(a1)$ of our  theorem is true. 

For $2 \le s \le n-2$, there exists a forest $F$ in $G$ with the properties:  both $F_{xy}$ and $F'_{xy}$ are forests with exactly two components, $R_x(F_{xy}) > 0$, and $R_y(F_{xy}) > 0$. Therefore by {\sc Claim 3}, $c_s(G) > c_s(G')$.
Obviously, if $x$ is  a cut vertex of $G$, 
then $c_{n-1}(G) = c_{n-1}(G')$.
If $x$ is not a cut vertex of $G$, then by {\sc Claim 1'},
$c_{n-1}(G) > c'_{n-1}(G)$, and so  
$c_{n-1}(G) > c'_{n-1}(G) \ge c_{n-1}(G')$. Thus,
claim $(a2)$ is true. 
\ep
\\[1.5ex]
\indent
The arguments similar to those in the proof of 
{\bf \ref{coef<}} provide one of possible proofs of 
{\bf \ref{coef<=}}.
\bs
\label{A,B,aHb}
Let  $A$, $B$,  $D$, and $H$ be disjoint graphs, 
$x, y \in V(H)$ and $x \ne y$, and
$d$ is a vertex in  $V(D)$ incident to an edge.
Let $A$ be a path with an end-vertex $a$ and $B$  a 
path with an end-vertex $b$.
Let graph $R$  be obtained from $A$, $B$, and $H$ by identifying $x$ with $a$ and $y$ with $b$.
Let $G_a$ and $G_b$ be obtained from 
$R$ and $D$ by identifying $d$ with $a$ and $b$, respectively.
Suppose that 
\\[1ex]
$(h1)$ there exist an automorphism 
$\eta : V(H) \to V(H)$ such that 
$\eta (a) = b$ and $\eta (b) = a$
{\em (and so two-pole $aHb$ is symmetric)} and
\\[1ex]
$(h2)$ $v(A) \le v(B)$.
\\[0.7ex]
\indent
Then 
\\[1ex]
$(a1)$
$G_a \succeq ^c G_b$ and, moreover,
\\[1ex]
$(a2)$ $G_a \succ\succ ^c  G_b \Leftrightarrow v(A) < v(B)$.
\es

{\bf Proof.}
We prove our claim by induction on $v(A \cup B) = n$. 
If $n = 2$, then our claim is true by 
{\bf \ref{A,B,xHy-operation}}. 
So we assume that $n \ge 3$ and that our claim is true if $v(A \cup B) < n$. Our goal is to prove that the claim  is also true for $v(A \cup B) = n$.

Let $s \in \{2 \ldots v(G) - 1\}$.
Let ${\cal F}_s(G)$ denote the set of spanning forests of $G$ with $s$ edges.
Let $T$ be a tree in $D$ such that $d \in V(T)$.
Put $\{v, z\} = \{a,b\}$ and let
\\[1ex]
\indent
$\sigma _s(G_z) = 
\sum \{\gamma (F): F \in {\cal F}_s(G_z), T \subseteq F\}$.
\\[1ex]
\indent
{\bf It is sufficient to show that 
$\sigma _s(G_a) \ge \sigma _s(G_b)$.} 
\\[1ex]
Let $e(T) = t$ and ${\cal P}_1$ denote the set of trees $P$ in $R$ containing $a$ and $b$ and such that $e(P) \le s - t$.
Let ${\cal P}_2$ denote the set of pairs $(P_a, P_b)$
such that $V(P_a) \cap V(P_b) = \emptyset $,
$P_a$ and $P_b$ are trees, $a \in V(P_a)$,
$b \in V(P_b)$, and $e(P_a \cup P_b) \le s- t$.
Given a spanning forest $F$ of $G$ and $z \in V(G)$, let 
$F_z$ denote the component of $F$ containing $z$.
\\[1ex]
For $P \in {\cal P}_1$, let 
 ${\cal F}_s(G,  P) =
\{F \in {\cal F}_s(G): F_a =  F_b = P \}$.
\\[1ex]
For $(P_a, P_b) \in {\cal P}_2$, let
 ${\cal F}_s(G, (P_a, P_b)) =
\{F \in {\cal F}_s(G): F_a = P_a, F_b = P_b \}$.
\\[1ex]
For $z \in \{a,b\}$, let 
$\sigma _s (G_z, (P_a, P_b)) = \sum \{ \gamma (F): 
F \in {\cal F}_s(G_z, (P_a, P_b)\}$ and
\\[1ex]
$\sigma _s (G_z, P) = \sum \{ \gamma (F): 
F \in{\cal F}_s(G_z, (P)$.
%
Then
\\[1ex]
$\sigma _s (G_z, (P_a, P_b)) =
(v(P_z) + t) v(P_v) \sum \{\gamma (F - (P_a \cup P_b)): F \in {\cal F}_s(G_z, (P_a, P_b))\}$ and
\\[1ex]
$\sigma _s (G_z,  P) =
(v(P) + t)  \sum \{\gamma (F - P): F \in {\cal F}_s(G_z, P)\}$.
\\[1ex]
Let 
\\[1ex]
$\sigma _s^1 (G_z) = \sum \{ \sigma _s (G_z,  P): 
P \in {\cal P}_1\}$ and
$\sigma _s^2 (G_z) = \sum \{ \sigma _s (G_z,  (P_a, P_b)): (P_a, P_b) \in {\cal P}_2\}$.
\\[1ex]
Then
$\sigma _s(G_z) = \sigma _s^1 (G_z)  + 
\sigma _s^2 (G_z)$.
Let $\Delta  _s(R) = \sigma _s(G_a) - \sigma _s(G_b)$.
\\[2ex]
\indent
{\bf We need to  prove that $\Delta  _s(R) \ge 0$.}
\\[1ex]
Clearly, $\sigma _s (G_a,  P) = \sigma _s (G_b,  P)$.
Therefore
\\[1ex]
$\sigma _s^1 (G_a) = \sigma _s^1 (G_b)$, and so
$\Delta  _s(R) = \sigma _s^2(G_a) - \sigma _s^2(G_b)$.
Given $(P_a, P_b) \in {\cal P}_2$, let 
\\[1ex]
$\Delta  _s(R, (P_a, P_b)) = 
\sigma _s (G_a,  (P_a, P_b)) - \sigma _s (G_b,  (P_a, P_b))$. 
\\[1ex]
Let $\sigma (R, (P_a, P_b)) = \sum \{\gamma (F - (P_a \cup P_b)): F \in {\cal F}_s(G, (P_a, P_b))\}$.
Then 
\\[1ex]
$\Delta  _s(R, (P_a, P_b)) =
[(v(P_a) + t) v(P_b) - (v(P_b) + t) v(P_a)]
\sigma (R, (P_a, P_b)) = 
\\[1ex]
t [v(P_b) - v(P_a)] \sigma (R, (P_a, P_b))$.
\\[1ex]
Let 
${\cal P}' = \{(P_a, P_b) \in {\cal P}\}: 
v(P_b) < v(A \cup H - b),  v(P_a) \ne v(P_b)\}$ 
\\[1ex]
and 
${\cal P}'' = \{(P_a, P_b) \in {\cal P}: 
v(P_b) \ge v(A \cup H - b)\}$.
\\[1ex]
Let, accordingly,
\\[1ex]
$\Delta ' _s(R) = \sum \{\Delta  _s(R, (P_a, P_b)):
(P_a, P_b) \in {\cal P}'\}$ and
\\[1ex]
$\Delta  ''_s(R) = \sum \{\Delta  _s(R, (P_a, P_b)):
(P_a, P_b) \in {\cal P}''\}$.
\\[1ex]
Then 
$\Delta  _s(R) = \Delta ' _s(R) + \Delta ''_s(R)$.
\\[1ex]
Since $P_a \subseteq A \cup H - b$,
we  have:
$ v(P_b) \ge v(A \cup H - b)\} 
\Rightarrow v(P_b) \ge  v(P_a)$.
\\[1ex]
Since $t \ge 0$, clearly
$v(P_b) \ge  v(P_a) \Rightarrow 
\Delta  _s(R, (P_a, P_b)) \ge 0$.
\\[1ex]
Now by the previous inequality,
$ v(P_b) \ge v(A \cup H - b)\} 
 \Rightarrow  \Delta  _s(R, (P_a, P_b)) \ge 0$.
\\[1ex]
Therefore $\Delta  ''_s(R) \ge 0$.
Moreover, if $v(B) > v(A)$, then there exist 
$P_a A$ and $P_b $ such that
$v(P_b) > v(P_a)$. Therefore, 
if $v(B) > v(A)$ and $t > 0$, then $\Delta  ''_s(R > 0$.
Since  vertex $d$ in $D$ is not an isolated vertex, there exists a tree  $T$ in $D$ containing $d$ with
$t = e(T) > 0$.
\\[2ex]
\indent
{\bf Thus, it is sufficient to show that} $\Delta '_s(R) \ge 0$.
\\[1ex]
Let ${\cal T}$ denote the set of pairs $(T_a, T_b)$
such that $T_a$ and $T_b$ are trees in $H$,
\\[1ex]
$V(T_a) \cap V(T_b) = \emptyset $,
$a \in V(T_a)$, and
$b \in V(T_b)$.
Let $T'_a = \eta [T_a]$ and $T'_b = \eta [T_b]$
and put 
\\[1ex]
$\pi (T_a, T_b)  = (T'_a, T'_b)$.
Then $\pi : {\cal T} \to {\cal T}$ is a bijection.
\\[1ex]
Let 
\\[1ex]
${\cal T}' = \{(T_a, T_b) \in {\cal T}:  \pi (T_a, T_b) =  
(T_a, T_b)\}$ and
\\[1ex]
${\cal T}'' = \{(T_a, T_b) \in {\cal T}:  \pi (T_a, T_b) \ne  (T_a, T_b)\}$.
\\[1ex]
\indent
Let  ${\cal L}$ denote the set of pairs $(L_a, L_b)$
such that $L_a$ is a path in $A$ containing $a$,  $L_b$ is 
a path in $B$ containing $b$, 
and $v(L_a) \le  v(L_b)$.
Let $\mu (L_a, L_b) = (L'_a, L'_b)$, 
where $L'_a$ is the path in $A$  such that $a \in V(L'_a)$ and $v(L'_a) = v(L_b)$ and $L'_b$ is the path in $B$  such that $b \in V(L'_b)$ and $v(L'_b) = v(L_a)$.
Let ${\cal L}' = \{(L_a, L_b) \in {\cal L}:  \mu (L_a, L_b) =  (L_a, L_b)\}$ and
${\cal L}'' = \{(L_a, L_b) \in {\cal L}:  \mu (L_a, L_b) \ne  (L_a, L_b)\}$.
\\[1ex]
\indent
For $(T_a, T_b) \in {\cal T}$ and $(L_a, L_b) \in {\cal L}$,
let 
$\sigma _s (G_z, (T_a, T_b), (L_a, L_b)) = 
\sigma _s (G_z, (P_a, P_b))$, where
$P_a = T_a \cup L_a$ and $P_b = T_b \cup L_b$.
\\[1ex]
\indent
Let $r = r(T_a, T_b), (L_a, L_b)) = 
s - t - e(T_a \cup T_b \cup L_a \cup L_b)\}$.
\\[1ex]
Since $e(T_z) = e(T'_z)$ and $e(L_z) = e(L'_z)$ for
$z \in\{a,b\}$, we have:
\\[1ex]
$ r(T_a, T_b), (L_a, L_b)) =  r(T'_a, T'_b), (L_a, L_b)) = 
r(T_a, T_b), (L'_a, L'_b)) =  r(T'_a, T'_b), (L'_a, L'_b))$.
\\[1ex]
Therefore
\\[1ex]
$c_r(R - T_a - T_b - L_a - L_b) = 
c_r(R - T'_a - T'_b - L_a - L_b)$ and
\\[1ex]
$c_r(R - T_a - T_b - L'_a - L'_b) = 
c_r(R - T'_a - T'_b - L'_a - L'_b)$.
\\[1ex]
Then 
\\[1ex]
$\sigma _s (G_a, (T_a, T_b), (L_a, L_b)) =
(v(T_a) + e(L_a) + t)( v(T_b) + e(L_b))
c_r(R - T_a - T_b - L_a - L_b)$,
\\[1ex]
$\sigma _s (G_a, (T'_a, T'_b), (L_a, L_b)) =
(v(T'_a) + e(L_a) + t)( v(T'_b) + e(L_b)) 
c_r(R - T'_a - T'_b - L_a - L_b)$,
\\[1ex]
$\sigma _s (G_a, (T_a, T_b), (L'_a, L'_b)) =
(v(T_a) + e(L'_a) + t)( v(T_b) + e(L'_b)) 
c_r(R - T_a - T_b - L'_a - L'_b)$, 
\\[1ex]
$\sigma _s (G_a, (T'_a, T'_b), (L'_a, L'_b)) =
(v(T'_a) + e(L'_a) + t)( v(T'_b) + e(L'_b)) 
c_r(R - T'_a - T'_b - L'_a - L'_b)$.
\\[2ex]
\indent
The similar formulas for $\sigma _s (G_b, ...)$'s are obtained from the above formulas for 
$\sigma _s (G_a, ...)$'s by moving $t$ from the right bracket to the left one. 
Let
\\[2ex]
$\omega _s (G_z, (T_a, T_b), (L_a, L_b)) =
\sigma _s (G_z, (T_a, T_b), (L_a, L_b)) +
\sigma _s (G_z, (T'_a, T'_b), (L_a, L_b)) +
\\[1ex]
\sigma _s (G_z, (T_a, T_b), (L'_a, L'_b)) +
\sigma _s (G_z, (T'_a, T'_b), (L'_a, L'_b))$
\\[1ex]
and 
\\[1ex]
$\Delta '_s(R,  (T_a, T_b), (L_a, L_b)) =
\omega _s (G_a, (T_a, T_b), (L_a, L_b)) -
\omega _s (G_b, (T_a, T_b), (L_a, L_b))$.
\\[1ex]
Then
\\[1.3ex]
$\Delta '_s(R) = 
\sum \{\Delta '_s(R,  (T_a, T_b), (L_a, L_b):
(T_a, T_b) \in  {\cal T}'', (L_a, L_b) \in  {\cal L}''\} + 
\\[1.3ex]
\frac{1}{2} \sum \{\Delta '_s(R,  (T_a, T_b), (L_a, L_b):
(T_a, T_b) \in  {\cal T}', (L_a, L_b) \in  {\cal L}''\} +
\\[1.3ex]
\frac{1}{2} \sum \{\Delta '_s(R,  (T_a, T_b), (L_a, L_b):
(T_a, T_b) \in  {\cal T}'', (L_a, L_b) \in  {\cal L}'\} +
\\[1.3ex]
\frac{1}{4} \sum \{\Delta '_s(R,  (T_a, T_b), (L_a, L_b):
(T_a, T_b) \in  {\cal T}', (L_a, L_b) \in  {\cal L}'\}$.
\\[2ex]
\indent
{\bf Finally, it is sufficient to show that 
each $\Delta '_s(R,  (T_a, T_b), (L_a, L_b)) \ge 0$.} 
\\[1ex]
From the above formula  we have:
\\[1ex]
$~~~~~~~~~~~~~~~~~~~~~~~~~~~~~~~~~~~~~~~
\Delta '_s(R,  (T_a, T_b), (L_a, L_b)) =$
\bqn
\label{Delta's}
t (v(L_b) - v(L_a)) 
[c_r(R - T_a - T_b - L_a - L_b) - c_r(R - T_a - T_b - L'_a - L'_b)]
\eqn
Let $R_{ab} = R - T_a - T_b$.
Then $R_{ab}$ is the disjoint union of three graphs:
\\[1ex]
$R_{ab} = (A - a) \cup (B - b) \cup (H - T_a - T_b)$, and so  
\\[1ex]
%
$R - T_a - T_b - L_a - L_b =
(A - L_a) \cup (B - L_b) \cup (H - T_a - T_b)$.
\\[1ex]
\indent
Therefore
\\[1ex]
$c_r(R - T_a - T_b - L_a - L_b) =
\sum \{c_p((A - L_a) \cup (B - L_b)) c_q(H - T_a - T_b):
p + q = r\}$ and
\\[1ex]
$c_r(R - T_a - T_b - L'_a - L'_b) =
\sum \{(c_p(A - L'_a) \cup (B - L'_b)) c_q(H - T_a - T_b):
p + q = r\}$.
\\[1ex]
\indent
Thus 
\\[1.3ex]
$[c_r(R - T_a - T_b - L_a - L_b) - c_r(R - T_a - T_b - L'_a - L'_b)] =
\\[1ex]
\sum \{[c_p((A - L_a) \cup (B - L_b)) - c_p((A - L'_a) \cup (B - L'_b))] c_q(H - T_a - T_b): p + q = r\}$.
\\[1.3ex]
\indent
Let 
$A - L_a = A_1$, $B - L_b = B_1$, 
$A - L'_a = A_2$, and $B - L'_b = B_2$.
\\[1ex]
Then $v(A_1 \cup B_1) = v(A_2 \cup B_2)$.
Since $v(L'_a) = v(L_b)$, $v(L'_b) = v(L_a)$, and 
$v(L_a) \le v(L_b)$, 
\\[0.7ex]
clearly $v(A_1) \ge v(A_2)$.
\\[1.5ex]
\indent
Let $\bar{A} = A_2$, $\bar{B} = B_2$, and
$\bar{D}$ be the  path with $v(A_1) - v(A_2) +1$ vertices disjoint from $\bar{A} \cup \bar{B}$.
Let $\bar{a}$, $\bar{b}$, and $\bar{d}$ be  end-vertices of 
paths $\bar{A}$, $\bar{B}$, and $\bar{D}$, respectively.
Let $\bar{H}$ be the graph consisting of  exactly two isolated vertices $\bar{x}$ and $\bar{y}$.
Obviously,  
\\[0.7ex]
$v( A ) \le v( B ) \Rightarrow v(\bar{A}) \le v(\bar{B})$.
\\[0.7ex]
\indent
Let us put in the claim (we are proving)
\\[0.7ex]
 $(A, a): = (\bar{A}, \bar{a})$, 
$(B, b): = (\bar{B}, \bar{b})$, $(D, d): = (\bar{D}, \bar{d})$, 
$(H, x, y): = (\bar{H}, \bar{x}, \bar{y})$,
$G_a: = A_1 \cup B_1$,
\\[0.7ex]
 and $G_b: = A_2 \cup B_2$.
 \\[0.7ex]
 \indent
Clearly, 
after this replacement the assumptions $(h1)$ and $(h2)$ are satisfied and
\\[0.7ex] 
$v(\bar{A} \cup \bar{B}) < v(A \cup B) = n$. 
Therefore by the induction hypothesis, 
\\[0.7ex]
$c_p(A_1 \cup B_1) \ge c_p(A_2 \cup B_2)$.
Obviously, 
\\[0.7ex]
$c_p(A_1 \cup B_1) - c_p(A_2 \cup B_2) = 
c_p((A - L_a) \cup (B - L_b)) - c_p((A - L'_a) \cup (B - L'_b))$.
\\[0.7ex]
\indent
Now since $v(L_b) \ge v(L_a)$,
we have from (\ref{Delta's}):
$\Delta '_s(R,  (T_a, T_b), (L_a, L_b)) \ge 0$.
\ep
\\[2ex]
\indent
Using the arguments similar to those in the proof of 
{\bf \ref{A,B,aHb}}, it is not hard to prove the following generalization of {\bf \ref{A,B,aHb}}.
\bs
\label{A,F,aHb}
Let  $A$,  $D$, $F$, and $H$ be disjoint graphs, 
$d$ be a non-solated vertex of $D$, $b \in V(F)$,
 and 
$x, y \in V(H)$, where  $x \ne y$.
Let $A$ be a path with an end-vertex $a$.
Let $R$  be obtained from $A$, $F$,  and $H$ by identifying $x$ with $a$ and $y$ with $b$.
Let graphs $G_a$ and $G_b$ be obtained from 
$R$ and $D$ by identifying $d$ with $a$ and $b$, respectively.
Suppose that 
\\[1ex]
$(h1)$ two-pole $xHy$ is symmetric and
\\[1ex]
$(h2)$ $F$ has a path $bBt$ such that 
$v(A) \le v(B)$.
Then  
\\[1ex]
$(a1)$ $G_a \succeq ^c G_b$ and
\\[1ex]
$(a2)$ $v(A) < v(B) \Rightarrow G_a \succ\succ  ^c G_b$.
\es

The operation described in {\bf \ref{A,F,aHb}}
is shown in Figure \ref{A,F,aHb-operation}.

\begin{figure}
\begin{center}
\includegraphics[scale=0.8]{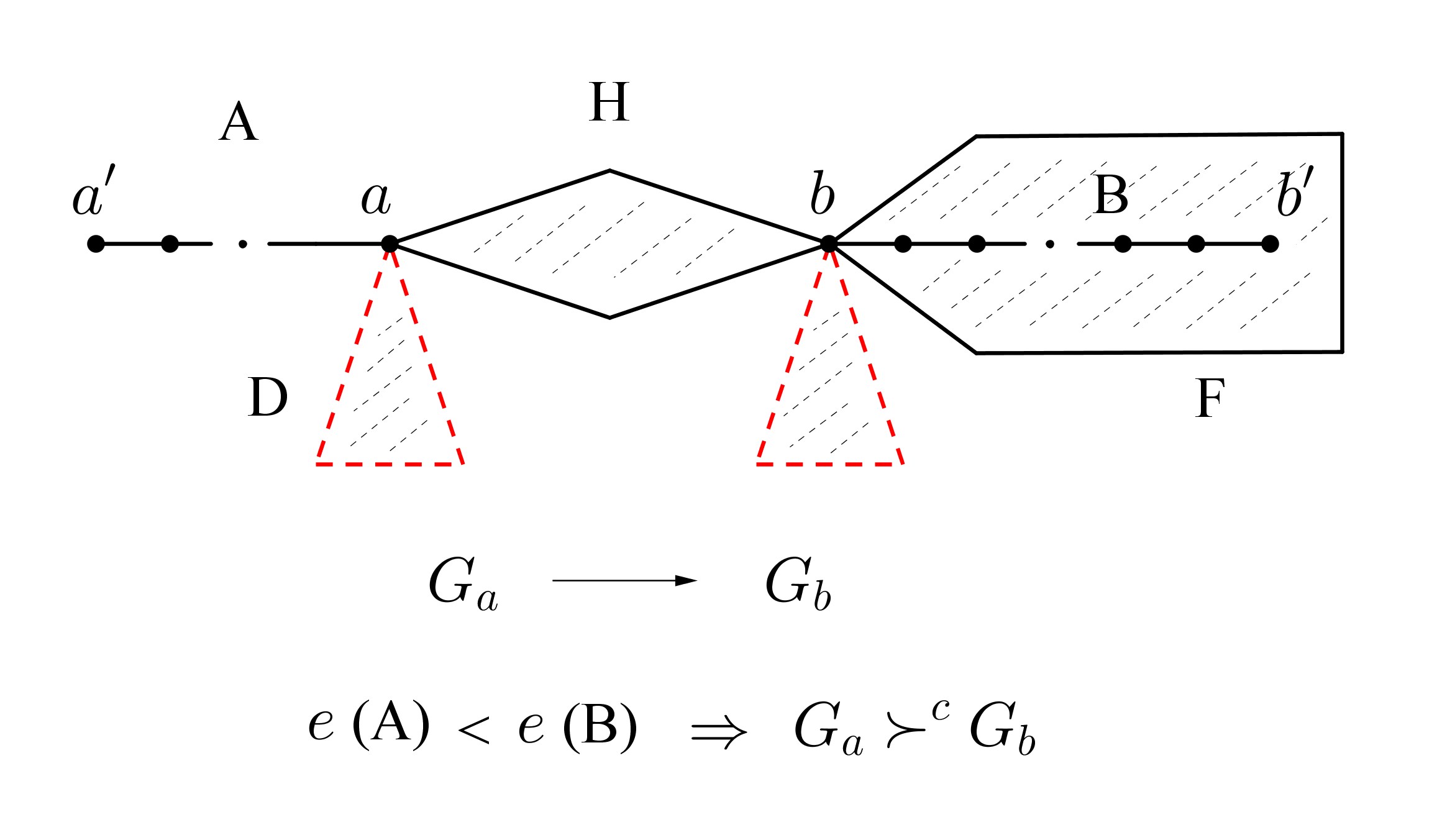}
\end{center}
\caption{The operation in {\bf \ref{A,F,aHb}}}
\label{A,F,aHb-operation}
\end{figure}

\begin{figure}
\begin{center}
\includegraphics[scale=0.48]{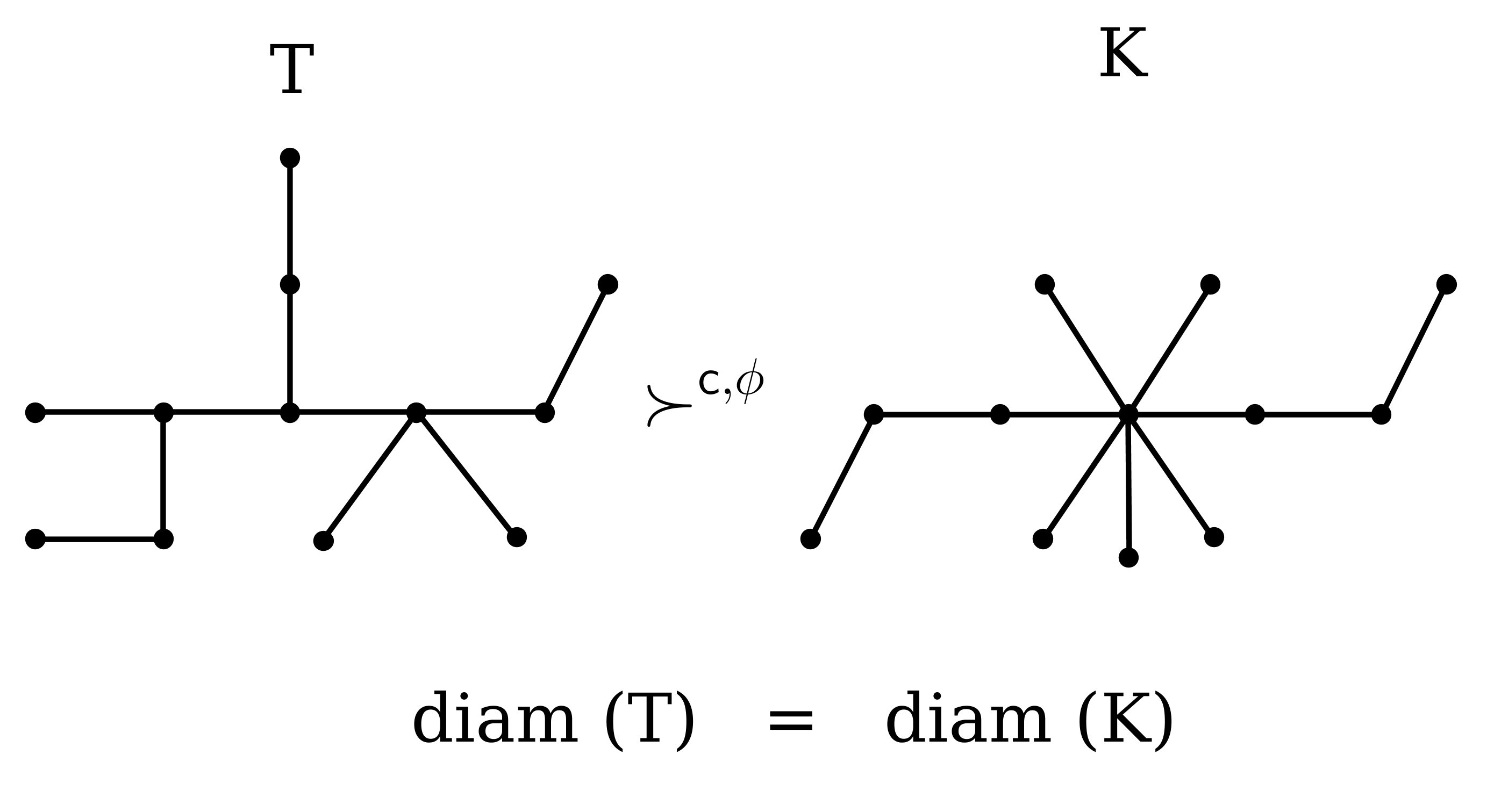}
\end{center}
\caption{$T \succ^{c,\phi} K$}
\label{D>K}
\end{figure}

Now, using the $\succeq ^c$-monotonicity of the $xPy$-operation and the operation described in 
{\bf \ref{A,F,aHb}}, we can obtain the results on the poset 
$({\cal T}_n, \succeq ^c)$ similar to {\bf \ref{Dr<Yr}} and {\bf \ref{Tr<Yr}} on the poset $({\cal T}_n, \succeq ^\phi)$.
\\[1ex]
\indent
We remind that for $n$-vertex trees $T$ and $T'$, 
$T \succ\succ ^c T'$ if and only if $c_s(T) > c_s(T')$ for every $s \in \{2, \ldots, n-2\}$.
Let, as above,  $\delta '(G)$ be the minimum non-leaf vertex degree in $G$.
The following theorem is an analog of  {\bf \ref{Dr<Yr}} with $ \succ\succ  ^\phi $  replaced by $\succ\succ  ^c$
(see the corresponding definitions in 
{\bf \ref{D(r),Y}}). 
\bs {\em \cite{Kcourse}}
\label{Yr>cDr} 
Let $r \ge 3$ and $n \ge r + 2$.
Then 
\\[1ex]
$(a1)$ for every $D \in {\cal D}_n(r) \setminus {\cal K}_n(r)$ there exists $Y \in {\cal K}_n(r)$ such that $D \succ\succ  ^c Y$,  
\\[1ex]
$(a2)$  $D \succ\succ  ^c K_n(r)$ for every 
$D \in {\cal K}_n(r) \setminus \{K_n(r)\}$, and therefore, 
{\em (from $(a1)$ and $(a2)$)}
\\[1ex]
$(a3)$ 
$D  \succ\succ  ^c K_n(r)$ for every 
$D \in {\cal D}_n(r) \setminus \{K_n(r)\}$
{\em (see Figure \ref{D>K})},
\\[1ex]
$(a4)$
$({\cal D}_n(3), \succeq ^c )$ is a linear order poset, namely, for $T, T' \in {\cal D}_n(3)$ we have:
\\[1ex]
 $T \succ\succ ^c T' \Leftrightarrow \delta '(T) > \delta ' (T')$, and
\\[1ex]
$(a5)$
$({\cal D}_n(4,r), \succeq ^c )$ is a linear  poset, namely, for $T, T' \in {\cal D}_n(4)$ we have:
\\[1ex]
 $T \succ\succ ^c T' \Leftrightarrow  T > ^u T' $.
\es

Obviously, claim $(a3)$ in {\bf \ref{Yr>cDr}} follows from
claims $(a1)$ and $(a2)$.
Claim $(a1)$  in {\bf \ref{Yr>cDr}} follows from 
{\bf \ref{coef<=}} and the fact that every  tree in $D \in {\cal D}_n(r) \setminus {\cal K}_n(r)$ can be transformed into a tree  in
${\cal K}_n(r)$ by a series of $xPy$-operations (so that every intermediate tree is also in ${\cal D}_n(r)$).
Claim $(a2)$ in {\bf \ref{Yr>cDr}} follows from 
{\bf \ref{A,F,aHb}} and the fact that every  tree in 
$D \in {\cal K}_n(r)\setminus \{K_n(r)\} $ can be transformed into $K_n(r)$
by a series of operations described in {\bf \ref{A,F,aHb}} 
(so that every intermediate tree is also in ${\cal K}_n(r)$).
Claims $(a4)$ and $(a5)$  can be proven in the same way using {\bf \ref{A,F,aHb}}.
\\[2ex]
\indent
Similarly, we can prove the following analog of {\bf \ref{Tr<Yr}} for $\succ\succ ^c$
(see the corresponding definitions in {\bf \ref{L(r),S(r)}}):
 \bs {\em \cite{Kcourse}}
\label{Yr>^cLr} 
Let $r \ge 3$, $n \ge r + 2$,   
and   
$L \in {\cal L}_n(r)$. 
Then 
\\[1ex]
$(a0)$ $L_n(r) \succ\succ  ^c  L_n(r + 1)$ for every 
$r \in \{2, \dots , n - 2\}$, 
\\[1ex]
$(a1)$
$({\cal S}_n(r), \succeq ^c )$ is a linear  poset, namely, for $T, T' \in {\cal S}_n(r)$ we have:
\\[1ex]
 $T \succ\succ ^c T' \Leftrightarrow  T > ^w T' $,
\\[1ex]
$(a2)$
$M_n(r) \succ\succ  ^c L $ for every 
$L \in {\cal S}_n(r) \setminus \{M_n(r)\}$, 
\\[1ex]
$(a3)$ for every 
$L \in {\cal L}_n(r) \setminus {\cal S}_n(r)$ 
there exists $Z \in {\cal S}_n(r)$ such that 
$L \succ\succ  ^c Z$,
\\[1ex]
$(a4)$ 
$L  \succ\succ  ^c L_n(r)$ for every 
$L \in {\cal S}_n(r) \setminus \{ L_n(r)\}$,  and
therefore 
{\em (from $(a3)$ and $(a4)$)}
\\[1ex]
$(a5)$
$L  \succ\succ  ^c L_n(r)$ for every 
$L \in {\cal L}_n(r) \setminus \{L_n(r)\}$, i.e.  $L_n(r)$ is the  $\succ\succ  ^c$-minimum graph in ${\cal L}_n(r)$.
\es

\begin{figure}
\begin{center}
\includegraphics[scale=0.40]{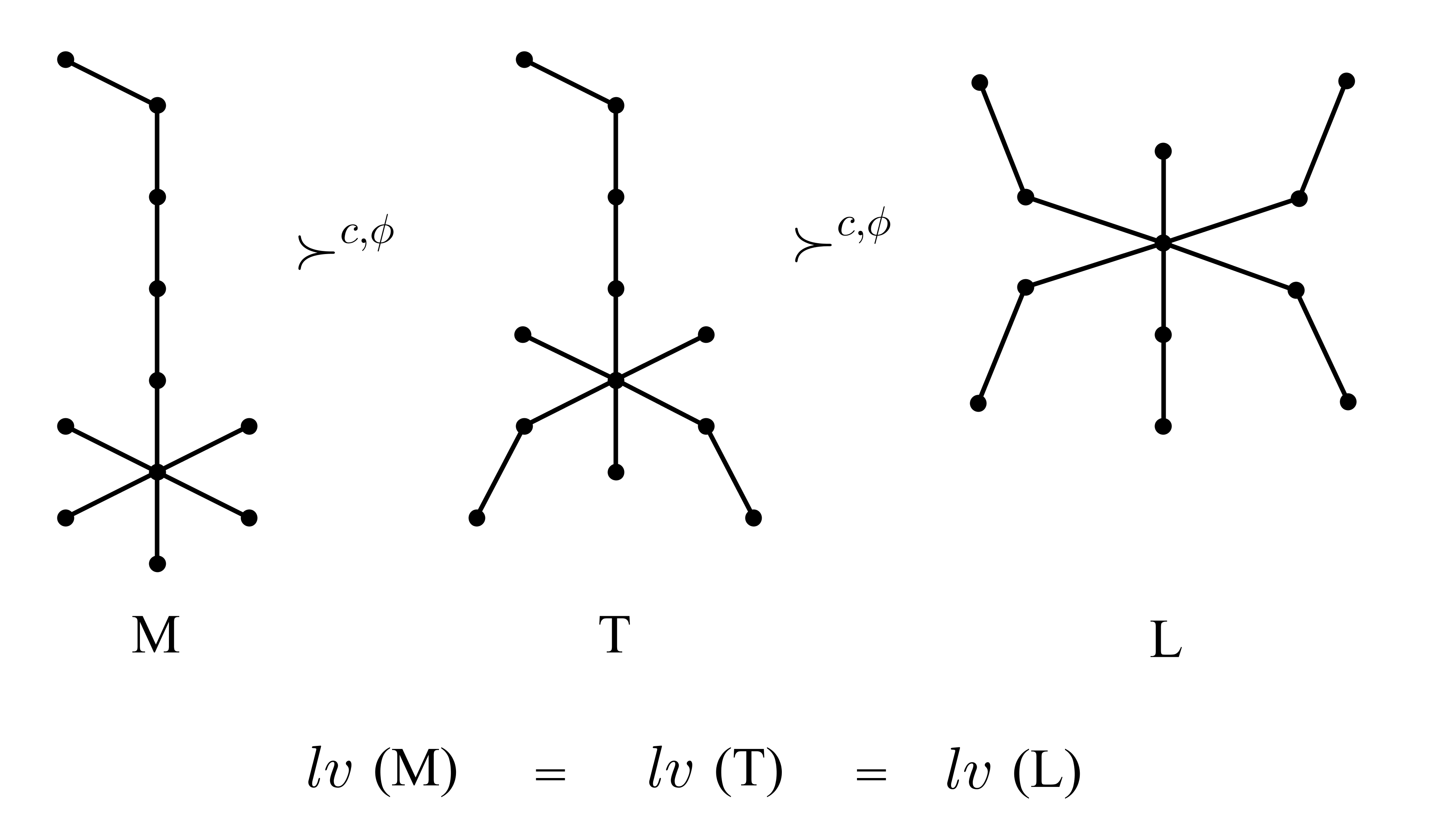}
\end{center}
\caption{$M \succ ^{c, \phi }T \succ ^{c, \phi } L$}
\label{L>L(r)}
\end{figure}


Figure \ref{L>L(r)} illustrates claims $(a2)$ and $(a5)$ in {\bf \ref{Yr>^cLr}}.
\\[1ex]
\indent
Let ${\cal L}^1_n(r)$ denote the set of graphs with $n$ vertices, $r$ leaves, and exactly one cycle (and so $v(G) = e(G) = n$). 
It can be shown that unlike ${\cal L}_n(r)$ the set ${\cal L}^1_n(r)$ does not have in general the  $\succ\succ  ^c$-minimum graph.
Here are some results illustrating this situation.
%
\\[1ex]
\indent
Let $Y = Y^s(r)$ be the tree obtained from $r \ge 1$ disjoint paths
$P_i$,  $1\le i \le r$, of $s \ge 1$ edges  by specifying one end-vertex of each path and identifying these specified  end-vertices of all paths. Then 
$Y$ has $n = rs +1$ vertices, exactly one vertex (say, $y$) of degree $r$,  exactly $r$ vertices of degree 1, and  the other vertices of degree two, and so $Y \in {\cal S}_n(r)$.

 Let  $Z$ be a path with $s \ge 2$ edges,  $z$ a leave of $Z$, and $z'$ the vertex in $Z$ adjacent to $z$.
 Let $\triangle $ be the triangle and
 $F$ be the graph obtained from disjoint $\triangle $ and $Z$ by identifying a  vertex  in $\triangle $ with vertex $z$ in $Z$.

Now let $A =  A^s(r)$ and $B =  B^s(r)$ be the graphs obtained from disjoint $F$ and $Y = Y^s(r-1)$ by identifying vertex $y$ in $Y$ with $z$ and $z'$, respectively.
Then $A^s(r)$ is isomorphic to $Y^s(r)$.
Both $A$ and $B$ have $n = rs + 3$ vertices and edges, exactly $r$ leaves, and exactly one cycle which is the triangle $\triangle $, and so $A, B \in {\cal L}^1_n(r)$. 

As above, $c_i(G)$ is the $i$-th coefficient of  the Laplacian polynomial of $G$, and so  $c_0(G) = 1$, $c_1(G) = 2e(G)$, and $c_{n-1}(G) = n t(G)$.

\bs {\em \cite{Kcourse}}
\label{Ars,Brs}  
Let $n$ and $r$ be integers such that $n \ge r + 3$ and $r \ge 2$.
Then
\\[1ex]
$(a0)$ 
$c_1 (A) = c_1(B) = 2n$ and $c_{n-1}(A) = c_{n-1}(B) = 3n$,
\\[1ex]
$(a1)$  $c_{n-2}(A) > c_{n-2} (B)$, and
\\[1ex]
$(a2)$ $c_2 (B) > c_2 (A)$.
\es

Here are more details about graphs  $A$ and $B$ for $r = 2$.
\bs 
{\em \cite{Kcourse}}
\label{A,B} 
Let $A$ and $B$ be graphs described above with $r = 2$, and so $v(A) = v(B) = n = 2s + 3$.  Then
\\[1ex]
$(a1)$ if $s = 4$ $($i.e. $n = 11$$)$, then $c_9(A) > c_9 (B)$ and
$c_j (B) > c_j(A)$ for $8 \ge j \ge 2$,
\\[1ex]
$(a2)$ if $s = 5$ or $6$, then $c_i(A) > c_i (B) $ for $2s +1 \ge i \ge 2s $ and
$c_j(B) > c_j(A)$ for 
\\[0.5ex]
$2s -1  \ge j \ge 2$,
\\[1ex]
$(a3)$ if $s = 7$, then $c_i(A) > c_i (B) $  for $2s +1 \ge i \ge 2s -1$ and
$c_j(B) > c_j(A)$ for $2s - 2 \ge j \ge 2$, and
\\[1ex]
$(a4)$ if $s \ge 7$, then $c_i(A) > c_i (B) $ for $2s +1 \ge i \ge 2s -1$ and 
$c_2(B) > c_2(A)$.
\es







Let $T$ be an $n$-vertex tree of maximum degree $r$.
Then $T$ can be transformed  to an $n$-vertex star-tree $S$ with $r$ leaves by a series of the inverse 
$P_{xy}$-operations. Then by {\bf \ref{coef<=}},
$S \succeq ^c T$. Therefore by {\bf  \ref{Yr>^cLr}} 
$(a2)$, we have:
 \bs {\em \cite{Kcourse}}
 \label{TwithMaxDg=r}
 Let $T$ be an $n$-vertex tree of maximum degree $r$ and $T$ is not isomorphic to $M_n(r)$.
 Then $M_n(r)  \succ\succ  ^c T$.
\es

Using the $\succeq^c$-monotonicity of the operation described in {\bf \ref{A,B,aHb}}, we can prove the following results on trees with exactly three leaves (i.e. on trees in ${\cal S}(3)$). 
\bs
 \label{T<Y} {\em \cite{Kcourse}}
For every tree $F$ non-isomorphic to a path there exists $T \in {\cal S}(3)$ with
 $v(T) = v(F)$ such that $T \succ\succ  ^c F$.
 \es

From {\bf \ref{bestTr}}, 
{\bf \ref{T>^cG}},
and {\bf \ref{T<Y}} we have:
\bs {\em \cite{Kcourse}}
\label{conG<Y}
Let $G$ be  connected graph and
$n =  e(G) + 1$ {\em (and so $v(G) \le n$)}. 
Suppose that $G$ is not a path and not in ${\cal S}_n(3)$. 
Then there exists  $T \in {\cal S}_n(3)$ 
such that $P_n \succ\succ  ^{c,\phi} M_n(3) \succ\succ  ^{c,\phi} T \succ\succ  ^{c,\phi} G$.
\es
 
From {\bf \ref{GtoTreshGraph}}, 
{\bf \ref{T>^cG}},
{\bf \ref{Gsucceq^pKxy(G)}},
and
{\bf \ref{coef<}} we have:
\bs {\em \cite{Kcourse}}
\label{G,threshold}
Let $G$ be a connected graph in ${\cal G}^m$. Then for every spanning tree $T$ of $G$ and every clique $K$ in $G$ 
there exist a tree $D$ and  
a connected threshold graph $H$ in ${\cal G}^m$ such that $D$ contains $T$, $H$ contains $K$, and  $D \succeq ^{c,p} G \succeq ^{c,p} H$. 
\es

Given a symmetric function $g$ on $k$ variables and a graph $F$ with $k$ components, let 
$g[F] = g\{v(C): C \in Cmp(F)\}$.
\bs {\bf Remark.}
 \label{remark}
 Let $G$ be a graph with $n$ vertices and let 
\\[0.7ex]
\indent
$q_s(G) = \sum \{g[F]: F \in {\cal F}(G), e(F) = s\}$, 
\\[0.7ex]
where  $g$ is  a symmetric concave function on 
$n - s$ variables.
Then theorems 
{\bf \ref{coef<K-operation}} and {\bf \ref{A,F,aHb}} remain true if $c_s(G)$ is replaced by $q_s(G)$.
 \es

\section{On reliability  posets of graphs}
\label{Rposet}

\indent

In this section we will describe some results on the monotonicity of the operations  defined in Section 
\ref{operations} with respect to some reliability posets  
and on the problems 
${\cal R}_{max}$ and  ${\cal R}_{min}$ of finding  maximum and minimum reliable graphs among the graphs of the same size.
\\[1.5ex]
\indent
As  above, $R^k(p,G)$
denotes the probability that
the random graph $(G,p)$ has at most $k$ components and $R^1(p,G) = R(p,G)$. 
Let $f_k(G)$ denote the number of spanning forests of $G$ with $k$ components.
Obviously, if $G \in {\cal G}_n^m$ and $k \ge n-2$, then
$R^k(p,G)$ depends only on $p$ and $m$ (and is easy to find). Therefore we will always assume  that if $G \in {\cal G}_n^m$, then $k \le n-3$.
For a graph $G \in \bar{{\cal G}}_n^m$, we have: 
\\[1.5ex]
\indent
$R^k(p,G) = \sum \{a^k_s(G) p^sq^{m-s}: 
s \in \{n - k, \ldots , m\}\}$,
\\[1.5ex]
where $q = 1 - p$ and $a^k_s(G)$ is the number of spanning subgraphs of $G$ with $s$ edges and at most 
$k$ components, and so  $a^k_{n-k}(G) = f_k(G)$.  
Let 
\\[1.5ex]
$Q^k_1(x, G) = \sum \{a^k_s(G) x^{m-s}: s \in \{n - k, \ldots , m\}\}$ 
and
\\[1.5ex]
$Q^k_0(x, G) = \sum \{a^k_s(G) x^{s}: s \in \{n - k, \ldots , m\}\}$. Then
\\[1.5ex]
$R^k(p,G) = 
p^m Q^k_1(q/p, G) = q^m Q^k_0(p/q, G)$.
\\[1.5ex]
\indent
It turns out that the coefficients of $R^k(x,G)$  (and accordingly, of $Q^k_0(x,G)$ and $Q^k_1(x,G))$ satisfy the following recursions similar to those for
$\Phi (\lambda ,G)$
(see {\bf \ref{recursion}} and {\bf \ref{recursionPhi}}).
\bs {\em \cite{Kcourse,Krandgraph}}
\label{recursionR(G)}
Let $G \in \bar{{\cal G}}_n^m$. Then 
\\[1.5ex]
$(a1)$ 
$(m-s) a^k_s(G) = \sum \{a^k_s(G - u); u \in E(G)\}\}$ for 
$s \in \{n - k, \ldots , m\}$, 
and 
\\[1.5ex]
$(a2)$
$s a^k_s(G) = \sum \{a^k_s(G /  u); u \in E(G)\}\}$ for 
$s \in \{n - k, \ldots , m\}$,
and so
\\[1.5ex]
$(a3)$
$Q^k_0(x,G) = \sum \{ \int _{t = 0}^x Q^k_0(t, G - u)dt:
 u \in E(G)\}$ and
\\[1.5ex]
$(a4)$
 $Q^k_1(x,G)  = \sum \{ \int _{t = 0}^x Q^k_1(t, G / u)dt:
 u \in E(G)\} $. 
\es

Function $R^k(p,G)$ and its coefficients $a^k_s(G)$ satisfy the following
useful ``{\em deletion-contraction}''  formulas: 
\bs {\em \cite{Kcourse,KoperProb,Kprob-network}}
\label{-,/recursionR(G)}
Let $G \in \bar{\cal G}_n^m$ be a graph with at least one edge and $u \in E(G)$.  Then 
\\[1.5ex]
$(a)$ 
$a^k_s(G) =a^k_s(G - u) + a^k_s(G/u)$ for every 
 $s \in \{n - k, \ldots , m\}$, 
 and so
\\[1.5ex]
$(r)$
$R^k(p,G) = R^k(p,G - u) + R^k(p, G/u)$.
\es

The relations in 
{\bf \ref{recursionR(G)}} and {\bf \ref{-,/recursionR(G)}} are pretty useful for proving some claims on 
$a^k_s(G)$ and $R^k(p,G)$ by induction.
For example, it is very easy to prove 
{{\bf \ref{maxR(G),cycl(G)<4}} $(a0)$ below 
by induction using  {\bf \ref{recursionR(G)}} $(a4)$.
\\[1ex]
\indent
The following result is a generalization of 
{\bf \ref{t(G),H-operation}}.
\bs {\em \cite{Kcourse,KoperProb}} 
\label{R(G),H-operation}
Let $G \in {\cal G}_n^m$ and let graph $H$ be obtained from $G$ by a symmetric hammock-operation.
Then 
\\[1ex]
\indent
$a_s^k(G) \ge a_s^k(H)$ for every $s \in \{n-k, \ldots , m\}$, i.e.
$G \succeq ^{a(k)} H$, and so
$G \succeq ^{r(k)} H$.
\es

In {\bf \cite{KoperProb}} this theorem 
(and in \cite{KoperTrees} theorem 
{\bf \ref{t(G),H-operation}}) was proved using
the deletion-contraction  formula in {\bf \ref{-,/recursionR(G)}}.
Theorem {\bf \ref{R(G),H-operation}}
can also be proved using the approach in the proof of {\bf \ref{coef<}} and in Section \ref{other-parameters}.
\\[1.5ex]
\indent
From {\bf \ref{GtoTreshGraph}} and {\bf \ref{R(G),H-operation}} we have:
\bs {\em \cite{Kcourse}}
\label{threshR^k<=}
For every graph  $G$ in ${\cal C}_n^m$ there exists a threshold graph $F$ in ${\cal C}_n^m$ such  that
\\[1ex]
$a_s^k(G) \ge a_s^k(F)$ for every 
$s \in \{n-k, \ldots , m\}$, i.e.
$G \succeq ^{a(k)} F$, and so $G \succeq ^{r(k)} F$.
\es

From {\bf \ref{properties-of-Kxy}} $(a2)$ and 
{\bf \ref{R(G),H-operation}} we have:

\bs {\em \cite{Kcourse,KoperTrees}}
\label{>R,D-operation}
Let $G \in {\cal G}_n^m$ and let graph $H$ be obtained from $G$ by a symmetric hammock-operation.
Then 
\\[1.5ex]
$(c1)$
$a_s^k([G]^c) \ge a_s^k([H]^c)$ for every $s \in \{n-k, \ldots , m\}$, i.e. 
$[G]^c \succeq ^{a(k)} [H]^c$, and so
\\[1.5ex]
$(c2)$
$[G]^c \succeq ^{r(k)} [H]^c$.
\es

The arguments in the proofs of  {\bf \ref{coef<K-operation}} can be used to prove the following.
\bs {\em \cite{Kcourse}}
\label{R(G),Kxy-operation}
Let $G \in {\cal G}_n^m$,  
$F$ be the graph obtained from $G$ by a symmetric  
${\cal K}_{xy}$-operation, and $k $ be a positive  integer at most $n$. 
Then $G  \succeq ^{r(k)} F$.
\es

From {\bf \ref{properties-of-Kxy}} $(a2)$ and 
{\bf \ref{R(G),Kxy-operation}} we have:
\bs {\em \cite{Kcourse}}
\label{R(G),Kxy-operation,cmpl}
Let $G \in {\cal G}_n^m$,  
$F$ be the graph obtained from $G$ by a symmetric  
${\cal K}_{xy}$-operation, and $k $ be a positive  integer at most $n$. 
Then $[G]^c  \succeq ^{r(k)} [F]^c$.
\es

Obviously, {\bf \ref{R(G),Kxy-operation}} and 
{\bf \ref{R(G),Kxy-operation,cmpl}} are 
generalizations of  {\bf \ref{R(G),H-operation}} and 
{\bf \ref{>R,D-operation}}, respectively.
\\[1.5ex]
\indent
From 
{\bf \ref{T>G}} and  
{\bf \ref{>R,D-operation}} we have the following generalization of {\bf \ref{T>^bG}}.
\bs {\em \cite{Kcourse,KoperTrees}} 
\label{T>^rG}
Suppose that $H$ is a connected graph and $n \ge e(H) +1$. Then for every spanning tree $T$ of $F$ there exists a tree $D$ with $e(H)$ edges such that $D$ contains $T$ and 
\\[1.5ex]
\indent
$a_s^k(K_n - E(D)) \ge a_s^k(K_n - E(H))$ for every 
$s \in \{n-k, \ldots , m\}$, i.e. 
\\[1.5ex]
\indent
$K_n - E(D) \succeq ^{a(k)} K_n - E(H)$, and so  
$G \succeq ^{r(k)} H$.
\es

We remind that
the numbers $r(G) = v(G) - cmp(G)$ and 
$r^*(G) = e(G) - v(G) + cmp(G)$ are called the {\em rank} and the {\em corank} (or the {\em cyclomatic number}) of a graph $G$.

\subsection{On the $\succ\succ ^{a(k)}$-maximization problem for  graphs with  ``small'' corank} 
\label{max-small}

\indent

Let us consider the following generalization of problem of 
problem ${\cal T}_{max}$ in Section \ref{introduction}:
find a graph $M \in {\cal C}_n^m$  such that 
\\[1ex]
\indent
$a_s^k(M) = \max \{a_s^k(G): G \in  \bar{{\cal G}}_n^m\}$, 
where $1 \le k < n$ and $s \in \{n - k, \ldots, m\}$.
\\[1ex]
\indent
Does a solution of this problem depends on $s$ or on $k$ ?
It turns out that a solution may depend on $s$ 
(see {\bf \ref{non-uniform}}).
From  claim {\bf \ref{maxR(G),cycl(G)<4}} it follows that if $m \le n+2$, then a solution does not depend on $s$, and therefore it is a solution for the $\succ\succ ^{a(k)}$- and $\succ\succ ^{r(k)}$-maximization problems as well.
\\[0.7ex]
\indent
Let $F$ be a graph and $\sigma $ a positive integer.
Let   ${\cal F}_{\sigma }$ denote the set of graphs $H$ that are obtained from $F$ by subdividing each edge $u$ by $s(u)$ vertices, where  
$\sum \{s(u): u \in E(F)\} = \sigma $, and so 
$v(H) = v(F) + \sigma $ and $e(H) = e(F) + \sigma$).
We call a graph $H$ from ${\cal F}_{\sigma }$
$\sigma $-{\em uniform} if $|s(u) - s(u')| \le 1$ for every $u, u' \in E(G$.
Let $z =  \lfloor \sigma /e(F)\rfloor$.
If $H$ is a $\sigma $-uniform graph in ${\cal F}_{\sigma }$, then
$s(u) \in \{z, z +1\}$
for every $u \in E(F)$. 
Let $E_0 = \{u \in E(F): s(u) =z\}$ and 
$E_1 = \{u \in E(F): s(u) = z + 1\}$. Let $F_i$ be the subgraphs  induced by $E_i$ in 
$F$, $i \in \{0,1\}$. 
Obviously,  if $\sigma = 0 \bmod e(F)$,
then $s(u) = z$ for every $u \in E(F)$, and 
therefore there is exactly one $\sigma $-uniform graph in ${\cal F}_{\sigma }$
up to isomorphism; we denote his graph $F^*$. 
\\[0.7ex]
\indent
Let $B$ denote  the graph with two vertices and three parallel edges. Then all $\sigma $-uniform graph in 
${\cal B}_\sigma $ are isomorphic to the same graph which we denote by  
$\Theta ^*_n$, where $n = \sigma + 2  = v(\Theta ^*_n)$. 
\\[0.7ex]
\indent
Let $Q = K_4$ and $n = \sigma + 4$.
\\[1ex]
\indent 
If $\sigma \bmod 6  \in \{0, 1, 5\}$, then
all $\sigma $-uniform graph in  ${\cal Q}_{\sigma }$ are isomorphic to the same graph which we denote by   
$Q ^*_n$.
\\[0.7ex]
\indent
Now we will define the $\sigma $-uniform graph 
$Q_n^*$ in case when 
$\sigma \bmod 6 \not \in \{0, 1, 5\}$, i.e. when $\sigma \bmod 6 = r  \in \{2,3,4\}$, and so $e(Q_1) = r$.
Let $H$ be a  $\sigma $-uniform graph in ${\cal Q}_{\sigma }$.
Put $Q_n^* = H$ if one of the following holds:
\\[0.7ex]
$(h1)$ $r = 2$ and $Q_1$ is a matching, 
\\[1ex]
$(h2)$
$r = 4$ and $Q_0$ is a matching, and 
\\[1ex]
$(h3)$
$r = 3$ and $Q_1$ (as well as $Q_0$) is a 3-edge path.
\\[1ex]
\indent
Let, as above, $O_n$ denote the cycle with $n$ vertices.
\\[1ex]
\indent
Using {\bf \ref{recursionR(G)}} and {\bf \ref{-,/recursionR(G)}}, 
 we were able  to prove  the following.
\bs {\em \cite{Krandgraph,Kcourse}}
\label{maxR(G),cycl(G)<4}
Let $G \in \bar{{\cal G}}_n^m$ and 
$k \in \{1, \ldots , n - 3\}$.
Then the following  holds.
\\[1ex]
$(a0)$  Suppose that $e(G)  = n$.
If  $G$ is not a cycle, then  
\\[1ex]
\indent
$O_n \succ\succ ^{a(k)} G$, and so
$O_n \succ\succ ^{r(k)} G$.
\\[1ex]
$(a1)$ Suppose that $e(G)  = n +1$.
If  $G$ is not 
$\Theta _n^*$, then 
\\[1ex]
\indent
$\Theta _n^* \succ\succ ^{a(k)} G$, and so
$\Theta _n^* \succ\succ ^{r(k)} G$.
\\[1ex]
$(a2)$ Suppose that $e(G)  = n + 2$.
If $G$ is not 
$Q _n^*$, then 
\\[1ex]
\indent
$Q _n^* \succ\succ ^{a(k)} G$, and so
$Q _n^* \succ\succ ^{r(k)} G$.
\es

\subsection{On the $\succ\succ ^{a(k)}$-maximization problem for  graphs with  ``large'' corank}
\label{max-large}

\indent

Let us reformulate the problem in the previous part \ref{max-small} as  follows:
find in a complete graph $K_n$ a set $Z$ of $z$ edges (and the corresponding subgraph of $K_n$ induced by $Z$) such that 
\\[1ex]
\indent
$a_s^k(K_n - Z) = \max \{a_s^k(K_n - A): A \subset E(K_n), 
~and~|A| = z\}$.
\\[1.5ex]
The following result gives a solution to this problem for the graphs with relatively ``large'' corank, i.e. with relatively ``small'' $|A|$.
\\[1ex]
\indent
From {\bf \ref{>R,D-operation}} we have, in particular, the following generalization of {\bf \ref{maxtree-graph}}.
\bs {\em \cite{Kcourse,KoperProb}}
\label{maxR-graph}
Let $H$ be a subgraph of $K_n$  and with $r$ edges with no isolated vertices.
Suppose that
$H$ is not isomorphic to $rP^1$ or $P^2 + (r-2)P^1$.
\\[1.5ex]
$(a1)$
If $n\ge 2r$, then
\\[1.5ex]
\indent
$K_n - E(rP^1) \succ\succ ^{a(k)} K_n - E(P^2 + (r-2)P^1)  
\succ\succ ^{a(k)} K_n - E(H)$, and so
\\[1.5ex]
\indent 
$K_n - E(rP^1) \succ\succ ^{r(k)} K_n - E(P^2 + (r-2)P^1)  
\succ\succ ^{r(k)} K_n - E(H)$.
\\[1.5ex]
$(a2)$
If $n= 2r -1$, then
\\[1.5ex]
\indent
$K_n - E(P^2 + (r-2)P^1)  
\succ\succ ^{a(k)} K_n - E(H)$, and so
\\[1.5ex]
\indent 
$K_n - E(P^2 + (r-2)P^1)  
\succ\succ ^{r(k)} K_n - E(H)$.
\es

\subsection{On the $\succ\succ ^{a(k)}$-minimization problem for  graphs with  ``small'' corank} 
\label{min-small}

\indent

Let us consider the following generalization of problem of 
problem ${\cal T}_{min}$ in Section \ref{introduction}:
find a graph $L \in {\cal C}_n^m$  such that 
\\[1.5ex]
\indent
$a_s^k(L) = \min \{a_s^k(G): G \in  {\cal C}_n^m\}$, 
where $1 \le k < n$ and $s \in \{n - k, \ldots, m\}$.
\\[1.5ex]
\indent
Does a solution of this problem depends on $s$ or on $k$ ?
The following result gives a solution to this problem for the graphs of relatively ``small'' corank. It turns out that in this case the solution does not depend on $s$ and is the solution for the $\succ\succ ^{a(k)}$- and  $\succ\succ ^{r(k)}$-minimization problems as well.
\bs {\em \cite{Kcourse,Kproceed}}
\label{minR(G)forsome(m,n))}
Let $G \in {\cal C}_n^m$ and $G \not \in {\cal H}_n^m$.
\\[1.5ex]
$(a0)$
If $m = n \ge 3$, then 
\\[1.5ex]
\indent
$G \succ\succ ^{a(k)} L$ for every 
$L \in {\cal H}_n^n$, and so 
$G \succ\succ ^{a(k)} F_n^n $ and $G \succ\succ ^{r(k)} F_n^n = W$.
\\[1.5ex]
$(a1)$ If $n \ge 4$ and $m = n+1$, then 
\\[1.5ex]
\indent
$G \succ\succ ^{a(k)}L$ for every 
$L \in {\cal H}_n^{n+1}$, and so
$G \succ\succ ^{a(k)} F_n^{n+1}$ and $G \succ\succ ^{r(k)} F_n^{n+1}$.
\\[1.5ex]
$(a2)$ If $n \ge 5$ and $m = n+2$, then
\\[1ex]
\indent 
$G \succ\succ ^{a(k)} L$ for every 
$L \in {\cal H}_n^{n+2}$, and so 
$G \succ\succ ^{a(k)} F_n^{n+2}$ and
$G \succ\succ ^{r(k)} F_n^{n+2}$.
\\[1.5ex]
$(a3)$ If $n \ge 6$ and $n+2 \le m \le 2n - 2$, 
then 
\\[1.5ex]
\indent
$G \succ\succ ^{a(k)} L$ for every 
$L \in {\cal H}_n^m$, and so 
$G \succ\succ ^{a(k)} F_n^m$
and $G \succ\succ ^{r(k)} F_n^m$.
\es

\subsection{On the $\succ\succ ^{a(k)}$-minimization problem for  graphs with  ``large'' corank} 
\label{min-large}

\indent

Let us reformulate the problem in the previous part \ref{min-small} as  follows:
find in a complete graph $K_n$ a set $Z$ of $z$ edges 
(and the corresponding subgraph of $K_n$ induced by $Z$) such that 
\\[1.5ex]
\indent
$a_s^k(K_n - Z) = \min \{a_s^k(K_n - A): A \subset E(K_n), 
K_n - A~is~connected,~and~|A| = z\}$.
\\[1.5ex]
\indent
A similar question is whether a solution of this problem depends on $s$ or on $k$.

The following result gives a solution to this problem for the graphs with relatively ``large'' corank, i.e. with relatively ``small'' $|A|$.
\bs {\em \cite{Kcourse,Kproceed}}
\label{minR-graph}
Let $H$ be a subgraph of $K_n$ with no isolated vertices and with $z \ge 1$ edges. 
\\[1.5ex]
$(a0)$
If $z \le n - 2$ and $H$ is not isomorphic to 
$F_{z+1}^z = S^z$, then 
\\[1.5ex]
\indent
$K_n - E(H) \succ\succ ^{a(k)} K_n - E(F_{z+1}^z)$, and so
$K_n - E(H) \succ\succ ^{r(k)} K_n - E(F_{z+1}^z)$.
\\[1.5ex]
$(a1)$
If $z = n - 1$, 
then 
\\[1.5ex]
\indent
$K_n - E(H) \succeq ^{a(k)} K_n - E(F_z^z )$, and so
$K_n - E(H) \succeq ^{r(k)} K_n - E(F_z^z )$.
\\[1ex]
$(a2)$
If $z = n$,  then 
\\[1.5ex]
\indent
$K_n - E(H) \succeq ^{a(k)} K_n - E(F_{z-1}^z)$, and so
$K_n - E(H) \succeq ^{r(k)} K_n - E(F_{z-1}^z)$.
\es 

Notice that claim $(a0)$ in {\bf \ref{minR-graph}} is a generalization of {\bf \ref{mintree-graph}} and follows from 
{\bf \ref{GtoTreshGraph}} and 
{\bf \ref{>R,D-operation}}.
\\[2ex]
\indent
In {\bf \ref{extreme-threshold-graph}} 
we defined the so-called extreme threshold graphs $F_n^m$ and the set  ${\cal H}_n^m$.
Obviously, 
$F_n^m \in {\cal H}_n^m \subseteq {\cal C}_n^m$ and
if $G, G' \in {\cal H}_n^m$, then $G \asymp ^r G'$.
\bs {\em \cite{Kcourse,Kproceed}}
\label{minR-graphF_n^m}
Let $n$ and $z$ be positive integers, $n \ge 3$, and 
$n \ge z$. Let $m = {n\choose 2} - z$.
Suppose that $m \ge n-1$
$($i.e. ${\cal C}_n^m \ne \emptyset $$)$ and $1 \le z \le n$. 
Then 
\\[1.5ex]
\indent
$G \in {\cal C}_n^m \setminus {\cal H}_n^m \Rightarrow
G  \succ\succ ^{a(k)} L$ for every $L \in  {\cal H}_n^m$, and in particular,  $G  \succ\succ ^{a(k)} F_n^m$.
\es 

\bp The following claim is obviously true.
\\[1.5ex]
\indent
{\sc Claim.}  
\\[1ex]
$(c1)$
{\em If $1 \le z \le n-2$, then
$F_n^m = K_n - E(F_{z+1}^z)$.
\\[1ex]
$(c2)$
If $z = n - 1$, then  $F_n^m = K_n - E(F_z^z)$.
\\[1ex]
$(c3)$
If $z = n$, then  $F_n^m = K_n - E(F_{z-1}^z)$.}
\\[1.5ex]
\indent
Now the claim of the theorem follows from the above {\sc Claim} and 
{\bf \ref{minR-graph}}.
\ep

\subsection{Some problems on the reliability poset of graphs}
\label{problems}

\indent

Using {\bf \ref{maxtree-graph}}, we were able to prove the following about the existence of 
$\succeq ^r$-maximum graphs in ${\cal G}_n^m$.
\bs
\label{non-uniform} 
{\em \cite{Kcourse,Krandgraph,Knotuniform,Krandgr}}
Let,  as above, $G \succeq ^r F$ if $R(p,G) \ge R(p,F)$ for every $p \in [0,1]$.
Then
\\[1.5ex]
$(a1)$ for every $n \le 5$ and $n-1 \le m \le e(K_n)$, 
${\cal G}_n^m$
has an $\succeq ^r$-maximum graph and
\\[1.5ex]
$(a2)$
for every $n \ge 6$ there exists $m = m(n)$ such that
${\cal G}_n^m$ does not 
have  an $\succeq ^r$-maximum graph.
\es

From {\bf \ref{non-uniform}} we have, in particular:
\bs
\label{non-uniform,succ^a} 
{\em \cite{Kcourse,Krandgraph,Knotuniform,Krandgr}}
For every $n \ge 6$ there exists $m = m(n)$ such that
${\cal G}_n^m$ does not 
have  an $\succeq ^a$-maximum graph.
\es

Paper \cite{Knotuniform} containing theorem 
{\bf \ref{non-uniform}} was also mentioned in a survey paper \cite{Kar} published in Journal of Graph Theory in 1982.
Accidentally, in 1986
the author of \cite{B} put forward a conjecture 
contradicting {\bf \ref{non-uniform}} and {\bf \ref{non-uniform,succ^a}}, namely, 
saying  that ${\cal G}_n^m$
has a $\succeq ^r$-maximum graph for every $(n, m)$ with $n-1 \le m \le e(K_n)$ and claiming in addition that his conjecture is true for every $n \le 6$. 
\\[1.5ex]
\indent
Here are some interesting problems related with
{\bf \ref{non-uniform}}.
As above, let  $M(p)$ be a graph in $\bar{{\cal G}}_n^m$, namely,
\\[1.5ex]
\indent
$~~~~~~~~~~~~~~~~~~~~~~~
R(p,M(p)) = \max \{R(p, G): G \in \bar{{\cal G}}_n^m\}$, 
\\[1.5ex]
i.e. $M(p)$ is a  most reliable  $p$-random graph in $\bar{{\cal G}}_n^m $.
\\[1ex]
\indent
Let ${\cal M}_n^m = \{M(p) \in {\cal G}_n^m: p \in [0,1]\}$.
Obviously, ${\cal M}_n^m$ is a finite set.
By {\bf \ref{non-uniform}}, there are infinitely many pairs $(n,m)$ such that $|{\cal M}|_n^m \ge 2$.
\bs {\bf Problem.} {\em \cite{Kcourse}}
\label{calMnmIsBounded}
Is there a number $N$ such that $|{\cal M}_n^m| \le N$ for every  $n \ge 2$ and $m \le n(n-1)/2$ $?$ 
\es

Given $G, F \in {\cal G}_n^m$ such that 
$R(p,G) \not \equiv R(p,F)$, let
$crs\{R(G,F)\}$ denote the number of zeros of 
$R(p,G) - R(p,F)$ in $(0,1)$(with their multiplicities). In \cite{Kprob-network} we observed that there are pairs  
$(G,F)$  of equi-size graphs  such that
$crs\{R(G,F)\}\ge 1$. Are there such pairs $(G,F)$ with 
$crs\{R(G,F)\}\ge 2$ ? This was a long standing question until
we gave in \cite{Kcross} a construction providing for every integer $k \ge 1$ an equi-size  pair $(G,F)$ such that
$crs\{R(G,F)\} = k$. A similar question arises for equi-size pairs $(G,F)$ of $R(p)$-maximum graphs.
Let $crs\{R_n^m\} = \max\{crs\{R(G,F)\}: G, F \in {\cal M}_n^m,
R(p,G) \not \equiv R(p,F)\}$.
As we have mentioned in {\bf \ref{non-uniform}}, for every $n \ge 6$ there exists $m = m(n)$ such that 
$crs\{R_n^m\} \ge 1$. 
\bs {\bf Problem.} {\em \cite{Kcourse}}
\label{cross}
Are there pairs $(n,m)$ such that $crs\{R_n^m\} \ge 2$ $?$
Furthermore, is it true that for every integer $k \ge 1$ there exists $(n,m)$ such that $crs\{R_n^m\} \ge  k$ $($or, moreover, $crs\{R_n^m\} =  k$$)$  $?$ 
\es

In Section \ref{introduction} we mentioned the  problem on the minimal elements of the $\succeq ^r$-poset on ${\cal C}_n^m$. Here is a more general version of this problem.
\bs {\bf Problem.} {\em \cite{Kcourse,Kproceed}}
\label{problem-succ^r}
Is the following claim true $?$
\\[1.5ex]
\indent
{\sc Claim}. Let $m \ge n-1$, and so ${\cal C}_n^m \ne \emptyset $. Then there exists $L \in {\cal C}_n^m$ such that 
$G  \succeq ^{r(k)} L$  for every $G \in {\cal C}_n^m$.
\es

A similar problem concerns the $\succeq ^{a(k)}$-minimal graphs in ${\cal C}_n^m$.
\bs {\bf Problem.} {\em \cite{Kcourse,Kproceed}}
\label{problem-succ^a}
Is the following claim true $?$
\\[1.5ex]
\indent
{\sc Claim}. Let $m \ge n-1$, and so ${\cal C}_n^m \ne \emptyset $. Then
there exists $L \in {\cal C}_n^m$ such that 
$G  \succeq ^{a(k)} L$ for every $G \in {\cal C}_n^m$. 
\es

Let, as above,  ${\cal F}_n^m$ denote the set of connected threshold graphs with $n$ vertices and $m$ edges, and so 
${\cal F}_n^m \subseteq {\cal C}_n^m$. 
From {\bf \ref{threshR^k<=}} it follows that Problem  
{\bf \ref{problem-succ^r}} is equivalent to the following problem.
\bs {\bf Problem.} {\em \cite{Kcourse,Kproceed}}
\label{problem-succ^rF}
Is the following claim true $?$
\\[1.5ex]
\indent
{\sc Claim.}  Let $m \ge n-1$.Then
there exists $F \in {\cal F}_n^m$ such that 
$G  \succeq ^{r(k)} F$ for every $G \in {\cal F}_n^m$. 
\es

Similarly, from {\bf \ref{threshR^k<=}} it follows that Problem  
{\bf \ref{problem-succ^a}} is equivalent to the following problem.
\bs {\bf Problem.} {\em \cite{Kcourse,Kproceed}}
\label{problem-succ^aF}
Is the following claim true $?$
\\[1.5ex]
\indent
{\sc Claim.}  Let $m \ge n-1$.Then
there exists $F \in {\cal F}_n^m$ such that 
$G  \succeq ^{a(k)} F$ for every $G \in {\cal F}_n^m$. 
\es

Obviously, {\sc Claim} in {\bf \ref{problem-succ^a}} implies
{\sc Claim} in {\bf \ref{problem-succ^r}} and 
{\sc Claim} in {\bf \ref{problem-succ^aF}} implies
{\sc Claim} in {\bf \ref{problem-succ^rF}}.
\\[2ex]
\indent
There are some results supporting {\sc Claim} in 
{\bf \ref{problem-succ^a}} and indicating that 
${\cal H}_n^m$ is the set of all $\succ\succ ^{a(k)}$-minimum graphs in ${\cal C}_n^m$ (see, for example, 
{\bf \ref{minR(G)forsome(m,n))}} and {\bf \ref{minR-graph}}).

\section{On some  other graph parameters}
\label{other-parameters}

\indent

We call a  component of a graph {\em non-trivial} if it has at least one edge, and {\em trivial} otherwise.
\\[0.5ex]
\indent
Let ${\cal F}(G, v, e)$ denote the set of forests $F$ in 
$G$ with $v$ vertices and $e$ edges such that each component of $F$ is non-trivial, and so $F$ is induced by its set of edges and has 
$ v - e$ components. In particular, 
${\cal F}(G, 2e, e) = M(G,e)$ is the set of matchings in $G$ with $e$ edges.
Given a vertex $x$ of $G$, let ${\cal F}^x(G, v, e)$ denote the set of forests $F$ in ${\cal F}(G, v, e)$ containing vertex $x$.
Let  $|{\cal F}^x(G, v, e)| = f^x(G, v, e)$ and 
$|M(G, f)| = m(G, e)$.
\\[1ex]
\indent
We  remind  some notations and notation 
from  the proof of {\bf \ref{coef<}}.
\\[1ex]
\indent
Let $G = (V, E)$ be a simple graph, $x, y \in V$, $x \ne y$, $Z = N_x(G) \cap N_y(G)$, $X = N_x(G) \setminus (Z \cup \{y\})$, and $Y = N_y(G) \setminus (Z \cup \{x\})$.
\\
Let $e_v = xv$ and $e'_v = yv$ for $v \in X$ and let
\\[1ex]
\indent
$E_x = [x,X] = \{e_v: v \in X\}$ and 
$  E_y = [y,X] = \{e'_v: v \in X\}$.
%
\\[1ex]
Let $G'$ be obtained from $G$ by the 
$\diamondsuit _{xy}$-operation, i.e. 
$G' = (G \setminus E_x) \cup E_y$.
Obviously,
\\[0.5ex]
 $[\{x,y\}, Z] \subseteq E(G)$ and
$E(G) - E_x = E(G') - E_y$.
\\[1ex]
\indent
Let $\varepsilon (u) = u$ if $u \in E(G) - E_x$ and
$\varepsilon (e_v) = e'_v$ if $e_v \in E_x$ 
(and so $v \in X$).
Then 
\\[0.5ex]
$\varepsilon : E(G) \to E(G')$ is a bijection.
For $U \subseteq E(G)$, let
$\varepsilon [U] = \{\varepsilon (u): u \in U\}$.
\\[1ex]
\indent
For a subgraph $S$ of $G$, let
$\vartheta (S)$ be the subgraph of $G$ such that
\\[0.5ex]
$V(\vartheta (S)) = V(S)$ and 
$E(\vartheta (S)) = (E(S) \setminus  E_x)\cup
\varepsilon [E(S) \cap E_x]$.
\\[1ex]
\indent
Let ${\cal S}(H)$ denote the set of subgraphs of a graph $H$.
Obviously, $\vartheta $ is a bijection from 
${\cal S}(G)$ to ${\cal S}(G')$.
For  ${\cal A} \subseteq {\cal S}(G)$, let
$\vartheta [{\cal A}] = \{\vartheta (A): A \in {\cal A}\}$.
\bs {\em \cite{Kcourse}}
\label{match}
Let $G$ be a connected graph with $n$ vertices, $x$ and $y$ two distinct vertices in $G$.
Let $v$ and $e$ be integers such that 
${\cal F}(G, v, e) \ne \emptyset $.
Let  
$G'$ be obtained from $G$ by 
the $\diamondsuit _{xy}$-operation.
Then 
$f^x(G, v, e) \ge f^x(G', v, e)$ and 
$m(G, e) \ge m(G',e)$.
\es

{\bf Proof}.
The proof of this theorem is similar to but 
simpler than the proof of {\bf \ref{coef<}}.
\\[1ex]
\indent
Let $P'$ be a forest in $G'$ with at most two components each being non-trivial and each meeting $\{x,y\}$.
Let $P$ be the subgraph in $G$ such that 
$P' = \vartheta (P)$.
Clearly, $\{x,y\} \subseteq  V(P') = V(P)$ and  
$e(P') = e(P)$. 
Let $\sigma  (xz) = yz$ and $\sigma (yz) = xz$
for every $z \in Z$. 
Obviously,
$\sigma   : [\{x,y\}, Z] \to [\{x,y\},Z]$ is a bijection.
For $S \subseteq [\{x,y\}, Z]$, let 
$\sigma   [S] = \{\sigma   (s): s \in S\}$.
Let $A = A(P) =  [\{x,y\}, Z] \cap E(P)$. 
Obviously, $A(P) = A(P')$.
Let 
$\bar{P} = (P -  A) \cup \sigma   [A]$
and 
$\bar{P}' = (P' -  A) \cup \sigma   [A]$.
Obviously, $\bar{P} = P$ and $\bar{P'} = P'$ if and only if 
$\sigma   [A] = A$.
Clearly, both $P$ and $P' = \vartheta (P)$ have at most two components.
\\[1ex]
\indent
We need  the following simple facts (see similar claims in the proof of {\bf \ref{coef<}}).
\\[1ex]
\indent
{\sc Claim 1.}
{\em Suppose that $P'$ is a tree and $P$ is not a tree.
Then $\bar{P}$ is a tree and $\bar{P}' = \vartheta (P)$ is not a tree.
}
\\[0.7ex]
\indent
{\em Proof.}
Since $v(P) = v(P')$,  $e(P) = e(P')$, $P'$ is a tree, and
$P$ is not a tree, we have: $P$ has a cycle $C$.
If $C$ does not contain vertex $x$ or $C$ contains two edges from $E[X]$, then  
$P'$ has a cycle, a contradiction. 
Therefore $C$ is the only cycle in $P$  and $C$ contains exactly one edge $xc$ with $c \in X$ and exactly one edge $xz$ with $z \in Z$.
Then $\bar{P}$ is a tree, $x \in V(P') = V(\bar{P})$, and $\bar{P}' = \vartheta (\bar{P})$ is not a tree (namely, it has a  cycle).
\epcl
\\[1ex]
\indent
It is also easy to prove the following two claims.
\\[1ex]
\indent
{\sc Claim 2.}
{\em $P'$ has two components if and only if  
$P = \vartheta ^{-1}(P)$ has  two components.}
\\[1ex]
\indent
{\sc Claim 3.}
{\em Suppose that both $P'$ and $P = \vartheta ^{-1}(P')$ have two components. Then
\\[1ex]
\indent
$(a1)$ both  $\bar {P}$ and $\bar {P}' = \vartheta (\bar {P})$ have two components and
\\[1ex]
\indent
$(a2)$
if $P'$ has no trivial components and $P$ has a trivial component $($namely, $y$$)$, then $\bar {P}$  has no  trivial components and $\bar {P}'$ has a trivial component $($namely, $x$$)$.}
\\[1.5ex]
\indent
Given a  forest $F $ of $G$ containing $x$ and $y$, let 
$F_{xy}$ be the minimal subforest of $F$ containing $x$ and $y$ and such that each component of $F_{xy}$ is a component of $F$. Obviously, $F_{xy}$ has at most two components each meeting $\{x,y\}$.
Let $\bar{F} = F - F_{xy} + \bar{F}_{xy}$.
\\[1ex]
\indent
Let ${\cal F} = {\cal F}^x(G, v, e)$, 
${\cal F}' = {\cal F}^x(G', v, e)$, and
$\tilde {{\cal F}} = \vartheta ^{-1}[{\cal F}']$.
\\[1ex]
Let ${\cal A}'$ and ${\cal B}'$ be subsets of ${\cal F}'$
such that 
${\cal F}' = {\cal A}' \cup {\cal B}'$, 
$\vartheta ^{-1}[{\cal A}'] \subseteq {\cal F}$, and
$\vartheta ^{-1}[{\cal B}'] \cap {\cal F} = \emptyset $,
\\[0.7ex]
 and so ${\cal A}' \cap{\cal B}' = \emptyset $.
Let $\tilde {{\cal A}} = \vartheta ^{-1}[{\cal A}']$ and 
$\tilde {{\cal B}} = \vartheta ^{-1}[{\cal B}']$. Then
\\[0.7ex]
\indent
$\tilde {{\cal F}} = \tilde {{\cal A}} \cup \tilde {{\cal B}}$,
$\tilde {{\cal A}} \cap \tilde {{\cal B}} = \emptyset $,
$|{\cal A}| = |{\cal A}'|$, and $|{\cal B}| = |{\cal B}'|$.
\\[1ex]
\indent
Suppose first that ${\cal B}' = \emptyset $. Then 
$ \vartheta ^{-1}$ is an injection from ${\cal F}'$ to 
${\cal F}$ and we are done.
\\[1ex]
\indent
Now suppose that  ${\cal B}' \ne \emptyset $. 
Let $c(H)$ denote the number of non-trivial components of a graph $H$.
Let $B' \in {\cal B}'$, and so $B = \vartheta ^{-1} (F')$ is not an element  of ${\cal F}$.
Obviously, $e(B) = e(B') = e$ and  $x \in V(B')$.
Now  $x \in V(B') \Rightarrow x \in V(B)$. 
Therefore $c(B) \ne c(B')$ or, equivalently, 
$c(B_{xy}) \ne c(B'_{xy})$.
Then by {\sc Claims} 1 and 3, $c(\bar{B}) = c(B')$ and
$c(\bar{B}') \ne  c(B')$, and so $ \bar{B} \in {\cal F}$ and
$\bar{B}' \not \in {\cal F}'$. 
Let  
\\[1ex]
\indent
$ \tilde {\vartheta }(F')  =   \vartheta ^{-1} (F')$ if 
$F' \in {\cal A}'$ and 
$ \tilde {\vartheta }(F')  =  \bar{F}$ if  $F' \in {\cal B}'$.
\\[1.5ex]
\indent
It is easy to see that the following holds. 
\\[1.5ex]
\indent
{\sc Claim 4.}
{\em Let $F'_1, F'_2 \in {\cal F}'$. Then
$F'_1 \ne F'_2 \Leftrightarrow 
 \tilde {\vartheta } (F'_1) \ne  \tilde {\vartheta } (F'_2)$.
}
\\[1ex]
It follows that  $\tilde {\vartheta }$ is an injection from 
${\cal F}'$ to ${\cal F}$.
Thus $|{\cal F}| \ge |{\cal F}'|$.
\\[0.5ex]
\indent
Similar (but much simpler) arguments show that 
$m(G, e) \ge m(G',e)$.
\ep
\\[1.5ex]
\indent
By the above definition, a component of a graph is {\em non-trivial} if it has at least one edge.
Notice that if this  definition  
is replaced by:  ``a component of a graph is 
{\em non-trivial} if it has at least $s$ edges with $s \ge 2$'', then the claim of theorem {\bf \ref{match}} is no longer true. 
\\[1ex]
\indent
Obviously,
$f^x(G, v(G), v(G) - 1) = t(G)$.
Therefore theorem {\bf \ref{t(G),H-operation}} $(a1)$ for the $\diamondsuit$-operation is a particular case of 
{\bf \ref{match}}. 
Also, it  follows from {\bf \ref{match}} that  a similar result holds for the number of  spanning forests  in $G$ with a given number of edges (or, the same, with a given number of components). 
\\[1ex]
\indent
Let ${\cal S}^x(G, v, e, k)$ denote the set of subgraphs of $G$ containing vertex $x$  and having $v$ vertices, $e$ edges and at most $k$ components. 
Let $s^x(G, v, e, k) = |{\cal S}^x(G, v, e, k)|$.
\\[1ex]
\indent
The arguments similar to those in the proof of {\bf \ref{match}}  can be used to prove  the following generalization of  {\bf \ref{match}}.
\bs {\em \cite{Kcourse}}
\label{x-subgraphs}
Let $G$ be a connected graph and$x$ and $y$ two distinct vertices in $G$.
Let $v$ and $e$ be integers such that 
${\cal S}^x(G, v, e, k) \ne \emptyset $.
Let  
$G'$ be obtained from $G$ by 
the $\diamondsuit _{xy}$-operation.
Then 
$s^x(G, v, e, k) \ge s^x(G', v, e, k)$.
\es

Let $h_0(G)$ and $h_1(G)$ denote the number of Hamiltonian cycles and Hamiltonian paths, respectively.
 \bs {\em \cite{Kcourse}}
\label{HamCycle,diamond-operation}
Let $G \in {\cal G}_n^m$ and 
$G'$ be the graph obtained from $G$ by an $\diamondsuit _{xy}$-operation.
Then 
$h_s(G) \ge h_s(G')$ for $s \in \{0,1\}$. 
\es

{\bf Proof}. The proof is similar to but much simpler than the proof of
{\bf \ref{coef<}}. 
\\[1ex]
\indent
For a graph $F$, let ${\cal H}_0(F)$ be the set of Hamiltonian cycles of $F$ and ${\cal H}_1(F)$  the set of Hamiltonian paths of $F$, and so $h_s(F) = |{\cal H}_s(F)|$.
\\[1.5ex]
\indent
{\sc Case 1.}
Consider  a Hamiltonian cycle $C'$ in $G'$.
Let $C'$ be a Hamiltonian cycle in $G'$ and $C = \vartheta ^{-1}(C')$.
If  $E(C') \cap [y, X] = \emptyset $, then $C$ is also 
a Hamiltonian cycle  in $G$.
If $x$ is incident to at most one edge in $G'$, then $G'$ has no Hamiltonian cycle.
Therefore we assume that
$x$ is incident to at least two edges in $G'$ and $E(C') \cap [y, X] \ne \emptyset $, and so $E(C') \cap [y, X]$ has either one or two edges. 
Obviously, $C'$ has either exactly two edges $xz_1$ and $xz_2$ in $[x,Z]$ or exactly one edge $xz$ in $[x,Z]$ and edge $xy$ (provided $xy \in E(G'))$.
\\[1.5ex]
\indent
{\sc Case 1.1.} Suppose that $C'$ has exactly one edge $xz$ in $[x,Z]$ and edge $xy$. Then $C'$ has exactly one edge $yx'$ in  $[y, X]$. Then $C = C' - yx' + xx'$ is not a Hamiltonian cycle in $G$, namely,   $y$ is incident to one edge in $C$ and $C - y$ is a cycle containing $xz$. 
Put $\bar{C} = C - xz + yz$. Then $\bar{C}$ is a 
Hamiltonian cycle in $G$ and $\bar{C}'$ is not 
a Hamiltonian cycle in $G'$.
\\[1,5ex]
\indent
{\sc Case 1.2.}
Now suppose that $C'$ has exactly two edges $xz_1$ and $xz_2$ in $[x,Z]$.
\\[1ex]
\indent
Suppose that $C'$  has exactly one edge $yx'$ in $[y, X]$
(and so $x' \in X$). Then $C = C' - yx' + xx'$ is not a Hamiltonian cycle in $G$, namely, 
$C$ has a unique cycle $D$ and vertex $y$ is of degree one in $C$. Obviously, $D$ has exactly one edge in $\{xz_1,xz_2\}$, say $xz_1$. Put $\bar{C} = C - xz_1 + yz_1$. Then $\bar{C}$ is a 
Hamiltonian cycle in $G$ and $\bar{C}'$ is not a Hamiltonian cycle in $G'$.
\\[1ex]
\indent
Now suppose that  $E(C') \cap [y, X]$ has exactly two edge $yx_1$ and $yx_2$ (and so $x_1,x_2 \in X$). 
Then $C = C' - \{yx_1, yx_2\} + \{xx_1, xx_2\}$ is not a Hamiltonian cycle in $G$, namely, $C$ has exactly two cycles $D_1$ and $D_2$, each $xx_i \in  D_i$, and $y$ is an isolated vertex in $C$.
 Put $\bar{C} = C -  \{xz_1, xz_2\} +  \{yz_1, yz_2\}$. 
 Then $\bar{C}$ is a Hamiltonian cycle in $G$ and 
 $\bar{C}'$ is not a Hamiltonian cycle in $G'$.
 \\[1.5ex]
 \indent
 {\sc Case 2.}
 Now consider a Hamiltonian path $P'$ in $G'$ and assume that $E(P') \cap [y, X] \ne \emptyset $, and so 
 $E(P') \cap [y, X]$ has either one or two edges. 
Obviously, either $P'$ contains $xy$ (provided $xy \in E(G')$) or has one or two edges in $[x,Z]$.
\\[1.5ex]
\indent
 {\sc Case 2.1.} Suppose that $xy \in E(P')$. 
If $x$ is an end of $P'$, then $P'$ has exactly one edge $x'y$ in  $[y, X]$ and $P = P' - x'y + x'x$ is a Hamiltonian path in $G$.
Otherwise,  the situation is similar to  {\sc Case 1.1}  on Hamiltonian cycles.
 \\[1.5ex]
\indent
 {\sc Case 2.2.} Now suppose that $P'$  has one or two edges in $[x,Z]$.
It is easy to see that $P$ is a Hamiltonian path in $G$ if and only if $P'$ has exactly one  edge in $[x,Z]$, exactly one edge $yx'$ in $ [y,X]$, and  in $P'$
$y$ is closer to $x$ than $x'$.
\\[1.5ex]
\indent
Suppose that $x$ is an end-vertex of $P'$, and so
$P'$ has exactly one edge $xz$  in $[x,Z]$. 
If $P'$  has exactly one edge $yx'$ in $[y,X]$ and in $P'$
$x'$ is closer to $x$ than $y$, then $P = P' - yx' + xx'$ is not a Hamiltonian path in $G$, namely, $P$ has exactly two components and one of them a cycle containing $xz$. Put $\bar{P} = P - xz + yz$. Then  
$\bar{P}$ is a Hamiltonian path in $G$ and $\bar{P}'$ is not a Hamiltonian path in $G'$.
If $P'$  has two edges $yx_1$ and $yx_2$ in $ [y,X]$,
then  
$P = P' - \{yx_1, yx_2\} + \{xx_1, xx_2\}$ is not a Hamiltonian path in $G$, namely, $P$ has exactly one  cycle $D$, $xz \in E(D)$, and $y$ is an isolated vertex in $P$. 
Put $\bar{P} = P - xz +  yz$. 
 Then $\bar{P}$ is a Hamiltonian path in $G$ and 
 $\bar{P}'$ is not a Hamiltonian cycle in $G'$.
\\[1.5ex]
\indent
Now suppose that $x$ is not an end-vertex of $P'$, and so
$P'$ has exactly two edges  in $[x,Z]$. Then the 
situation is similar to {\sc Case 2.1}  on Hamiltonian cycles.
\\[2ex]
\indent
Thus, from the above {\sc Cases} it follows that the following claim is true. 
\\[1ex]
{\sc Claim 1.}
{\em Let $s \in \{0,1\}$.
Suppose that $Q' \in {\cal H}_s(G')$ and 
$Q \not \in {\cal H}_s(G)$.
Then $\bar{Q} \in {\cal H}_s(G)$ and $\bar{Q}' \not \in {\cal H}_s(G')$. }
\\[1.5ex]
\indent
Let  ${\cal H} = {\cal H}_s(G)$ and 
${\cal H}' =   {\cal H}_s(G')$.
Let ${\cal A}'$ and ${\cal B}'$ be subsets of ${\cal H}'$
such that 
\\[0.7ex]
${\cal H}' = {\cal A}' \cup {\cal B}'$, 
$\vartheta ^{-1}[{\cal A}'] \subseteq {\cal H}$, and
$\vartheta ^{-1}[{\cal B}'] \cap {\cal H} = \emptyset $,
 and so ${\cal A}' \cap{\cal B}' = \emptyset $.
\\[0.7ex]
\indent 
Let $\tilde {{\cal A}} = \vartheta ^{-1}[{\cal A}']$ and 
$\tilde {{\cal B}} = \vartheta ^{-1}[{\cal B}']$. Then
\\[0.7ex]
$\tilde {{\cal H}} = \tilde {{\cal A}} \cup \tilde {{\cal B}}$,
$\tilde {{\cal A}} \cap \tilde {{\cal B}} = \emptyset $,
$|{\cal A}| = |{\cal A}'|$, and $|{\cal B}| = |{\cal B}'|$.
\\[1ex]
\indent
It is easy to see that the following holds. 
\\[1.5ex]
\indent
{\sc Claim 2.}
{\em Let $H'_1, H'_2 \in {\cal H}'$. Then
$H'_1 \ne H'_2 \Leftrightarrow 
 \tilde {\vartheta } (H'_1) \ne  \tilde {\vartheta } (H'_2)$.
}
\\[1ex]
\indent
Suppose first that ${\cal B}' = \emptyset $. Then 
by {\sc Claim}  2, $ \vartheta ^{-1}$ is an injection from 
${\cal H}'$ to 
${\cal H}$ and we are done.
\\[1ex]
\indent
Now suppose that  ${\cal B}' \ne \emptyset $. 
Let  
$ \tilde {\vartheta }(H')  =   \vartheta ^{-1} (H')$ if 
$H' \in {\cal A}'$ and 
$ \tilde {\vartheta }(H')  =  \bar{H}$ if  $H' \in {\cal B}'$.
Then  by {\sc Claims} 1 and 2, $\tilde {\vartheta }$ is an injection from 
${\cal H}'$ to ${\cal H}$.
Thus, $|{\cal H}| \ge |{\cal H}'|$.
\ep
\\[1.5ex]
\indent
Let $x \in  V(G)$ and $s$  a positive integer.
\\[1ex]
\indent
Let $h_0^x(G, s, k)$  denote 
the number of  subgraphs of $G$ containing $x$ and having $s$ vertices and at most $k$ components each being a cycle.
\\[1ex]
\indent
Let $h_1^x(G,s, k)$ denote 
the number of  subgraphs of $G$ containing $x$ and having $s$ vertices and exactly $k$ components each being a non-trivial path.
\\[1.5ex]
\indent
Using the arguments similar to those in the above proof
of {\bf \ref{HamCycle,diamond-operation}}, the following generalization of  {\bf \ref{HamCycle,diamond-operation}} can be proved.
 \bs 
\label{subgraphs,diamond-operation}
Let $G \in {\cal G}_n^m$, 
$G'$ be the graph obtained from $G$ by an $\diamondsuit _{xy}$-operation.
Let 
\\[1ex]
$\eta (G) \in \{h_0^x(G, s, k)$, $h_1^x(G,s, k) 
\}$. 
Then 
$\eta (G) \ge \eta (G')$. 
\es

It can be shown that if in the definition of $h_0^x(G, s, k)$
condition ``at most $k$ components'' is replaced by 
 ``exactly $k$ components'' (as in the definition of 
 $h_0^x(G, s, k)$), then theorem {\bf \ref {subgraphs,diamond-operation}} will be no longer true for  $h_0^x(G, s, k)$.
\\[1.5ex]
\indent
The results similar to 
{\bf \ref{subgraphs,diamond-operation}}
are also true for some other types of  subgraphs of 
a graph.

 \section{On weighted graphs and digraphs}
 \label{weighted-graphs}
 
 \indent
 
Many notions and results above can be naturally generalized to weighted graphs and digraphs of different type.
Here are some of them.
\\[1.5ex]
\indent
Let $D = (V, E)$ be a directed graph or simply, a digraph (and so $E \subseteq V\times V$) and $w$ a function from $E \to 
\mathbb{R}$.
We call the pair $(D, w)$ {\em a weighted digraph} also denoted by $D^w$.
If $(x,y) \in E \Leftrightarrow (y,x) \in E$  in $D$, then $D$ can be interpreted as an undirected graph $G$ and $(G^w)$ the corresponding  weighted graph.
\\[1.5ex]
\indent
For $v \in V$, 
let $d(v, G^w) = \sum \{w(vx): x \in V(G - v)\}$.
For $V = V(G) = 
\\[0.5ex]
\{v_1, \ldots , v_n\}$, let 
$A(G^w)$ be the $(V\times V)$-matrix $(a_{ij})$, 
where  each $a_{ii} = 0$ and each 
\\[0.5ex]
$a_{ij}=  w(v_iv_j)$ for $i \ne j$.
Let $R(D^w)$ be the diagonal  $(V\times V)$-matrix 
$(r_{ij})$, where 
\\[0.5ex]
each $r_{ii} = d(v_i, G^w)$ and each $r_{ij} = 0$ for $i \ne j$.
Let $L(D^w) = R(G^w) - A(G^w)$.
\\[1.5ex]
\indent
A  {\em directed tree $($or ditree$)$ $T$ rooted at $z$}  is a digraph  such that its  underlying graph is a tree, $z \in V(T)$,  $d_{out}(z) = 0$,
and $d_{out}(v) = 1$ for every $v \in V(T) \setminus z$.
\\[1ex]
\indent
Let
 ${\cal T}_z(D)$ denote the set of spanning ditrees of $D$ rooted at $z$, 
\\[1ex]
\indent 
 $w(T) = \prod \{w(e): e \in E(T)\}$ for $T \in {\cal T}_z(D)$, and 
$t_z(G^w) = 
\sum \{w(T)): T \in  {\cal T}_z(G)\}$. 
 \\[2ex]
\indent
Here is an obvious generalization of  the classical Matrix Tree Theorem for weighted digraphs. 
 \bs 
\label{MTT,digraph,w}
Let $D^w$ be a weighted digraph and $r \in V(D)$.
Then $t_r(G^w) = det (L_r(D^w))$.
\es 

Obviously,  
$det (\lambda I - L(D^w)) = \lambda P(\lambda , D^w)$, 
where $P(\lambda , D^w)$ is a polynomial of degree 
\\[1ex]
$n-1$.
Let $S(G^w) = (\lambda _1(D^w), \ldots , \lambda _{n-1}(D^w))$
be  the list of root roots of   $P(\lambda , D^w)$.
\\[1ex]
\indent
Let $C = (V, E)$ be a directed graph 
such that $E = V\times V$. Digraph 
$C$ is called the {\em complete digraph with the vertex set $V$}, and so every vertex of $C$ has a loop. Let $w$ be a function from  $V\times V$ to $\mathbb{R}$, and so 
$C^w$ is a weighted digraph. Given  a function 
$a  : V\times V \to \mathbb{R}$,
let $u = a  - w$.
Then $C^u$ is another digraph with the vertex set $V$. We call 
digraph $C^u$ $a $-{\em complement of} $C^w$ and write $C^u = [C^w]^a $.
\\[2ex]
\indent
Here is a generalization of the Reciprocity Theorem {\bf \ref{reciprocity}} for weighted digraphs. 
\bs 
{\em \cite{KelLatvia,Ktree1}}
\label{w-di-reciprocity}
Let $C^w$ be a weighted digraph with $n$ vertices.
\\[1ex]
Let 
$\sigma : (V\times V) \to \mathbb{R}$ be a  constant
 function: $\sigma (xy) = q  \in \mathbb{R}$ for every $xy \in V\times V$. 
Then 
\\[1.5ex]
$(a1)$ 
there is a bijection 
$\alpha : \{1, \ldots , n-1\} \to \{1, \ldots , n-1\}$ \\[1ex]
such that  
$\lambda _i(C^w) + 
\lambda _{\alpha (i)}([C^w]^q) = qn$
for every $i \in \{1, \ldots , n-1\}$ or, equivalently,
\\[1.5ex]
$(a2)$ $(qn - \lambda) L(\lambda , [C^w]^q) = (-1)^{n-1} \lambda L(qn - \lambda , C^w)$.
\es

A weighted digraph $C^w$ is called {\em $r$-out-regular} ({\em $r$-in-regular}) if $d_{out}(v, C^w) = r$ (respectively, $d_{in}(v, C^w) = r$) for every $v \in V(D)$.
Let  $A(\lambda , C^w) = 
det (\lambda I_n - A(D^w))$.
Obviously, if $C^w$ is  $r$-out-regular} or $r$-in-regular, then $A(r ,C^w) = 0$.
From {\ref{w-di-reciprocity} we have the corresponding Reciprocity Theorem on $A( \lambda, C^w)$ for  an out-regular (respectively, in-regular) weighted digraph $D^w$.
\bs {\em \cite{KelLatvia}}
\label{w-di-reciprocity-adj}
Let $C^w$ be a weighted $r$-out-regular or
$r$-in-regular digraph with $n$ vertices.
Then 
\\[1ex]
\indent
$(\lambda + r) A(\lambda , [C^w]^q) = 
(-1)^n (\lambda - qn + r)A(- \lambda , C^w)$.
\es

Let $p(x,y)$ be a polynomial of two variables $x$ and $y$.
Assuming that $xy \ne yx$, we call 
$p(x,y)$ a {\em $xy$-polynomial} ({\em $yx$-polynomial}) if $p$ does not contain ``the product $yx$'' (respectively, ``the product $xy$''). Therefore if $xy= yx$, then $p(x,y)$ is both $xy$ and $yx$-polynomial.
\\[1.5ex]
\indent
Using {\bf \ref{w-di-reciprocity-adj}} it is also easy to prove the following useful fact.
\bs
Let $C^w$ be a weighted $r$-out-regular
($r$-in-regular) digraph with $n$ vertices.
Let $J_{nn}$ be the $n\times n$-matrix with all entries equal 1 and 
$\{ \alpha _1,  \ldots, \alpha_n\}$ be the set of eigenvalues of $A(C^w)= A$, where 
$\alpha_n = r$.
 Let $p(x, y)$ be a polynomial with two variables $x$ and  $y$, and with real coefficients. 
 Suppose that $p(x, y)$ is an $xy$-polynomial if $C^w$ is $r$-out-regular and is an $yx$-polynomial if $C^w$ is $r$-in-regular.
 Then matrix $p(A, J_{nn})$ has the eigenvalues $p(r, n)$ and $p(\alpha _i, 0 )$ for $i = 1,\cdots, n-1.$
\es 

Let $G^w$ be a weighted graph and let
 ${\cal T}(G)$ denote the set of spanning trees of $G$.
Given $X \subseteq V$, let ${\cal F}(G, X)$ denote 
the set of spanning forests $F$ of $G$ such that every component of $F$ has exactly one vertex 
in $X$, and so
${\cal T}(G) = {\cal F}(G, x)$ for every $x \in V$.
For $T \in {\cal T}(G)$ and $F \in {\cal F}(G, X)$, let 
\\[1.5ex]
\indent
$w(T) = \prod \{w(e): e \in E(T)\}$ and 
$w(F) = \prod \{w(e): e \in E(F)\}$ 
\\[1.5ex]
and let
\\[1.5ex]
\indent
$t(G^w) = 
\sum \{w(T)): T \in  {\cal T}(G)\}$ and
$f(G, X) =  \sum \{w(F)): F \in {\cal F}(G, X)\}$. 
\\[1.5ex]
Then
$t(G^w) = f(G, x)$ for every $x \in V$, and so
$v(G)~t(G^w) = \sum \{ f(G, x): x \in V\}$.
\\[1.5ex]
\indent
For $x, y \in V(G)$, $x \ne y$, let $G^w/\{x,y\}$ denote the weighted graph $F^u$ such that 
\\[0.5ex]
$V(F) = V(G) \setminus \{x,y\} \cup t$ and $u(tz) = w(xz) + w(yz)$ for every $z \in V(G - \{x,y\})$ and 
\\[0.5ex]
$u(ab) = w(ab)$ for $a,b \in V(G - \{x,y\})$.
We say that $G^w/\{x,y\}$ {\em is obtained from $G^w$ 
\\[0.5ex]
by identifying its vertices $x$ and $y$}. For  $h = xy \in E(G)$, let
$G^w / xy = G^w/\{x,y\}$
 and 
 \\[0.5ex]
 $G^w - h = G^{w'}$, where
$w'(h) = 0$ and $w'(e) = w(e)$ for every $e \in E(G -h)$.
For a graph 
\\[0.5ex]
$G^w$ and $X \subseteq V(G)$, let 
$L_X(G^w)$ denote the matrix obtained from $L(G^w)$ by removing 
\\[0.5ex]
the rows and columns corresponding to every 
vertex $x \in X$ and
let $G^w_X = G^w /X$ denote 
\\[0.5ex]
the graph obtained from $G$ by identifying all vertices in $X$. 
\\[1.5ex]
\indent
Now it is easy to see that the following generalization of Theorems {\bf \ref{MTT}} and {\bf \ref{generMTT}} holds.
\bs {\em \cite{Kcourse}}
\label{MTTw}
Let $G^w$ be a weighted graph, $V = V(G)$, and 
$X \subseteq V(G)$.
Then
\\[0.5ex]
$t(G^w_{X}) = f(G^w, X) = det (L_X(G^w))$, and so  
$t(G^w) = det (L_v(G^w))$ for every $v$ in $V(G)$.
\es

As before,  
$det (\lambda I - L(G^w)) = \lambda P(\lambda , G^w)$, 
\\[1ex]
where $P(\lambda , G^w)$ is a polynomial of degree $n-1$ with the 
root sequence  
\\[1ex]
\indent
 $S(G^w) = (\lambda _1(G^w) \le \ldots , \le \lambda _{n-1}(K^w))$ and 
 \\[1.5ex]
\indent
$P(\lambda , G^w) = \sum \{(-1)^s c_s(G^w) \lambda ^{n- 1 - s}:  s \in \{0, \ldots , n-1\}\}$.
\\[1.5ex]
\indent
The following generalizations of {\bf \ref{coefficients}} and {\bf \ref{recursion}} 
 are true.
\bs {\em \cite{Kcourse}}
\label{w-coefficients}
Let $G^w$ be a  weighted graph with $n$ vertices and $m$ edges, $s \in \{0, \ldots , n-1\}$, and
$V = V(G)$.
Then
\\[1.3ex]
$(a1)$
$c_s(G) = 
\sum\{f(G^w, V - X) : X \subseteq V, |X| = s\} = 
\\[1.3ex]
\sum\{\gamma (F)\prod\{w(e): e \in E(F)\}: 
F \in {\cal F}(K), e(F) = s\}$ and
\\[1.7ex]
$(a2)$
$(m - s) c_s(K^w) = \sum \{c_s(K^w - e): e \in E(K)\}$
for $s \in \{0, \ldots , m\}$.
\es

Let, us above, $K = (V, E$ be a complete graph, and so $E =  {V \choose 2})$, and 
$K^w$ a weighed graph with the vertex set $V$.
Given a function $\sigma : {V \choose 2} \to \mathbb{R}$, let $u = \sigma n - w$.
We call $(K^ u)$ 
the $\sigma $-{\em complement of} $K^w$ and write $K^u = [K^w]^\sigma $.
 \\[1.5ex]
\indent
Here is a generalization of the Reciprocity Theorem 
{\bf \ref{reciprocity}}  for weighted graphs.
\bs 
{\em \cite{Kcourse,KeigvMtr,Ktree1}}
\label{w-reciprocity}
Let $K^w$ be a weighted  graph with $n$ vertices.
Let 
$\sigma : {V \choose 2} \to \mathbb{R}$ be a  constant
 function, namely, $\sigma (uv) = q  \in \mathbb{R}_+$ for every  $uv \in {V \choose 2}$.
 %
Then
\\[1.5ex]
$(a1)$ 
$\lambda _i(K,w) + \lambda _{n-i}(K,w)^q) = qn$ for every $i \in \{1, \ldots , n-1\}$, and so 
\\[1.5ex]
$(a2)$ $(qn - \lambda) L(\lambda , (K, w)^q) = (-1)^{n-1} \lambda L(qn - \lambda , (K, w))$
and 
\\[1.5ex]
$(a3)$ if $q \ge max\{w(e): e \in E(K)\}$ and 
$w(e) \ge 0$ for every $e \in E(K)$, then 
\\[1ex]
\indent
~~$0 \le \lambda _i(K,w) \le qn$ 
and $0 \le \lambda _i(K,w)^q) \le qn$ 
for every 
$i \in \{1, \ldots , n-1\}$.
\es

Theorem {\bf \ref{w-reciprocity}} was used in \cite{KeigvMtr} to give a  simple 
procedure for finding the spectrum of 
so-called totally decomposable symmetric matrices.
\\[1ex] 
\indent
Let $K = (V, E)$ be a complete graph
and
$w: E \to \mathbb{R}$ be a function. Let
\\[0.5ex]
$\{v_1, v_2\} \subseteq V$,  $S \subseteq V \setminus \{v_1, v_2\}$,
and $\varepsilon : S \to \mathbb{R}$ be a function.
For $i \in \{1,2\}$, let
\\[0.5ex]
$w_i: E \to  \mathbb{R}$ be a function such that
$w_i(e) = w(e)$ for every $e \in E \setminus [S, v_i]$ and 
\\[0.5ex]
$w_i(sv_i) = w(sv_i) + \varepsilon (s)$  for every $s \in S$.
The following result establishes an inequality 
\\[0.5ex]
between 
$t(K,w_1)$ and $ t(K,w_2)$ under certain conditions on functions  $w$ and $ \varepsilon $.
\bs {\em \cite{Kcourse,KoperTrees}}
\label{weightedDiamondTree}
Let $K = (V, E)$ be a complete graph.
Suppose that
\\[1ex]
$(h1)$
$w: E \to \mathbb{R}$ and $\varepsilon : S \to \mathbb{R}$ are non-negative functions and
\\[1ex]
$(h2)$
$w(sv_2) \ge w(sv_1)$ for every $s \in S$.
\\[1ex]
\indent
Then 
$t(K,w_1) \ge t(K,w_2)$.
\es

Obviously, {\bf \ref{t(G),H-operation}} follows from 
{\bf \ref{weightedDiamondTree}} when an $H$-operation is 
$\diamondsuit  $-operation. 
\\[1ex]
\indent
Here is an analog of {\bf \ref{weightedDiamondTree}} for graphs with randomly deleted  edges.
Let
$K = (V, E)$ be a complete graph
and
$q: E \to [0,1]$ be a function.
Let $\{v_1, v_2\} \subseteq V$,  $S \subseteq V \setminus \{v_1, v_2\}$,
and $\varepsilon : S \to [0,1]$ be a function.
We call $(K, q)$ a {\em random graph},  interpret $q(e)$ as the probability that edge $e$ of the complete graph $K$ does not exist  and assume that all edge events are mutually independent. Let $R^k(K,q)$ denote the probability that the random graph $(K,q)$ has at most $k$ components.
For $i \in \{1,2\}$, let $q_i: E \to  [0,1]$ be a function such that
$q_i(e) = q(e)$ for every $e \in E \setminus [S, v_i]$ and 
$q_i(sv_i) = q(sv_i)  \varepsilon (s)$  for every $s \in S$.
The following result establishes an inequality between 
$R^k(K,q_1)$ and $ R^k(K,q_2)$ under certain conditions on function  $q$.
\bs {\em \cite{Kcourse,KoperProb}}
\label{weightedDiamondRandom}
Suppose that
$q(sv_2) \ge q(sv_1)$ for every $s \in S$.
Then 
$R^k(K,q_2) \ge R^k(K,q_1)$.
\es

Obviously, {\bf \ref{R(G),H-operation}} follows from 
{\bf \ref{weightedDiamondRandom}} when an $H$-operation is 
$\diamondsuit  $-operation.
\\[1.5ex]
\indent
Now we will define an analogue of the ${\cal K}_{xy}$-operation for weighted graphs.
\\
Let $K = (V, E)$ be a complete graph, 
$\{v_1,v_2\} \subseteq V$,  $S \subseteq V \setminus \{v_1,v_2\}$ and let
 $\varepsilon $ be a function from $S$ to $\mathbb{R}$.
For $i \in \{1,2\}$,
let $r_i: S \to \mathbb{R}$ be a function such that
$r_i(s)  = w(sv_i)$ for $s \in S$.
Let $r'_1(s)  = r_1(s) + \varepsilon (s)$ and
$r'_2(s)  = r_2(s) - \varepsilon (s)$ for $s \in S$.
Put $u(e) = w(e)$ for $e \in E \setminus [\{v_1,v_2\}, S]$ and
$u(sv_i) = r'_i(s)$  for $s \in S$, and so $u$ is a function from $E$ to $\mathbb{R}$.
We say that $K^u$ {\em is obtained from $K^w$ by 
the $(v_1,v_2, \varepsilon)$-operation} and denote it 
$K^w(v_1,v_2, \varepsilon)$.

Let
$a:  E \to \mathbb{R}$ be the function such that    
$a(e) = w(e)$ for every $e \in E \setminus 
[\{v_1,v_2\}, S]$ and $a(e) = 0$, otherwise.
We call the $(v_1,v_2, \varepsilon)$-operation
{\em symmetric on $K^w$} if 
there exists an authomorphism 
$\sigma: V(K) \to V(K) $ of $K^a$ such that 
$\sigma (v_1) = v_2$  and $\sigma (v_2) = v_1$.
Let $B(K^w) = D'(K^w) + A(K^w)$, where $D'(K^w)$ is a diagonal $(V\times V)$-matrix.
\\[1ex]
\indent
It is easy to prove the following claim.
\bs
\label{delta(x^tAx)}  
Let $x \in \mathbb{R}^n$ and $x_S = x|_S$. 
Then
$x^tB(K^u)x - x^tB(K^w) x = 2(x_1 - x_2) 
(\varepsilon  \cdot x_S)$.
\es

Since $B(K^w)$ is a symmetric matrix, all eigenvalues of 
$B(K^w)$ are real numbers. Let $\beta (K^w)$ denote the maximum eigenvalue of  $B(K^w)$.
\bs {\em \cite{Kcourse}}
\label{alpha(As)<alpha(Az)} 
Let $K = (V, E)$ be a complete graph, $\{v_1, v_2\} \subseteq V$,  $S \subseteq V \setminus 
\{v_1, v_2\}$, and
 $\varepsilon : S \to \mathbb{R}$ and
$w: E \to \mathbb{R}$ be  functions.
Let $K^u = K^w(v_1,v_2, \varepsilon )$ be obtained from $K^w$ by the $(v_1,v_2, \varepsilon )$-operation.
Suppose that
\\[1ex]
$(h1)$
$w: E \to \mathbb{R}$ is a non-negative function,
\\[1ex]
$(h2)$
a weighted graph $K^w$ is  connected,
\\[1ex]
$(h3)$
$r_1 \ge r_2 \ge \varepsilon \ge \bar{0}$, i.e.,
$w(v_1s) \ge w(v_2s) \ge \varepsilon (s) \ge 0$ for every 
$s \in S$, and 
\\[1ex]
$(h4)$
the $(v_1,v_2, \varepsilon)$-operation is 
symmetric on $K^w$.
\\[1ex]
\indent
Then
$\beta (K^w(v_1,v_2, \varepsilon )) \ge  \beta (K^w) > 0$.
\es

{\bf Proof.}
Since  $B(K^w)$ is symmetric, all its eigenvalues are real numbers.
Since by $(h1)$, $w$ is a non-negative function, all entries of $B(K^w)$ are non-negative.
Since by $(h2)$, $K^w$ is connected,
matrix $B(K^w)$ is irreducible.
Let $x \in \mathbb{R}^n$ be an eigenvector of matrix 
$B(K^w)$ corresponding to its maximum eigenvalue 
$\beta (K^w)$. We can assume that $\parallel x \parallel = 1$.
Now by the Perron-Frobenius theorem \cite{G},
 $ \beta (K^w) > 0$ and all coordinates of $x$ are positive. 
By the Courant-Swartz theorem, 
\bqn
\label{alpha^w}
\beta (K^w) =  
\max \{z^tB(K^w)z: z \in \mathbb{R}^n,  \parallel z \parallel = 1\} = x^tA(K^w)x 
\eqn
and
\bqn
\label{alpha^u}
\beta (K^u) = 
\max \{z^tA(B^u ) z: z \in \mathbb{R}^n,  \parallel z \parallel = 1\} \ge x^tB(K^u) x.
\eqn

By {\bf \ref{delta(x^tAx)}},
\bqn
\label{x^tA(K^u)x=}
x^tB(K^w(a,b, \varepsilon ) x = x^tB(K^w)x + 
2(x_1 - x_2) (\varepsilon \cdot x_S) = 
\beta (K^w) + 2(x_1 - x_2) (\varepsilon \cdot x_S).
\eqn

By  $(h4)$, $K^u = K^w(v_1,v_2, \varepsilon) $ is isomorphic to
$K^w(v_2,v_1, \varepsilon ')$, where
 $\varepsilon '  = r_1 - r_2 + \varepsilon $.
 Since by $(h3)$, $r_1 \ge r_2 \ge \varepsilon \ge 0$, 
 we have $\varepsilon ' \ge 0$.
 Therefore we can assume that $x_1 \ge x_2$.
By $(h3)$, $u$ and $\varepsilon $  are non-negative functions. Since $x$ is positive, obviously $x|_S$ is also positive.  Therefore 
$(x_1 - x_2) (\varepsilon \cdot x_S) \ge 0$.
Now from (\ref{alpha^w}),  (\ref{alpha^u}), and 
(\ref{x^tA(K^u)x=}) we have
$\beta (K^u) \ge \beta (K^w)$.
\ep
\\[1.5ex]
\indent
Notice that if $w$ is a function from $E(K)$ to $\{0,1\}$, then $K^w$ is an ordinary graph $G$,
an $(x, y, \varepsilon)$-operation on $K^w$ is 
a ${\cal K}_{xy}$-operation on $G$, and 
$(v_1,v_2, \varepsilon)$-operation
is symmetric on $K^w$ if and only if
the ${\cal K}_{xy}$-operation on $G$ is symmetric.
The above results on a ${\cal K}_{xy}$-operation for a graph $G$ can be generalized to a symmetric  
$(x, y, \varepsilon)$-operation for a weighted graph 
$K^w$.
\\[1ex]
\indent
Let $Sp(K^w) = (V(K), sup(w))$.
It turns out that {\bf \ref{alpha(As)<alpha(Az)}} is also true for $\lambda (K^w)$ of a bipartite weighted graph $K^w$. 
Let 
$L(K^w) = R(K^w) - A(K^w)$ and 
$L^+(K^w) =R(K^w) + A(K^w)$.  Accordingly, let
$L(\lambda ,K^w) = det(\lambda I - L(K^w))$ and
$L^+(\lambda ,K^w) = det(\lambda I - L^+(K^w))$.
Let $n$ be the number of vertices of $K^w$.
Let $bip(K^w) = n + 1 - s$, where  $s - 1$ is  the degree of the polynomial  $L(\lambda ,K^w) - L^+(\lambda ,K^w)$, and so $L(\lambda ,K^w) = L^+(\lambda ,K^w)$ if and only if $bip(K^w) = n + 1$.
Let
$odc (K^w)$ denote the length of a shortest odd cycle in 
$Sp(K^w)$ if any exists, and put $odc (K^w) = n + 1$, otherwise, and so $K^w$ is bipartite if and only if 
$odc (K^w) = n + 1$. Thus, $K^w$ is bipartite if and only if $bip(K^w) = odc (K^w) = n + 1$.
\\[1ex]
\indent
Let $\lambda ^+(K^w)$ denote the maximum eigenvalue of $L^+(K^w)$ and, as before, $\lambda ((K^w)$ the maximum eigenvalue of $L(K^w)$.
It turns out that  the following holds.
\bs {\em \cite{Kcourse}}
\label{bip=odd} 
Let $K^w$ be a weighted graph. Then
$bip(K^w) = odc(K^w)$, and so $K^w$ is bipartite if and ony if $L(\lambda ,K^w) = L^+(\lambda ,K^w)$. In particular, if $K^w$ is bipartite, then  
$\lambda ((K^w) = \lambda ^+(K^w)$.
\es

From 
{\bf \ref{alpha(As)<alpha(Az)}}  and  {\bf \ref{bip=odd}}, applied to  the non-negative matrix $L^+(K^w)$, we have: 
\bs  {\em \cite{Kcourse}}
\label{biplambda(K^w)>}
Let $K^u$ and $K^w$ be weighted graphs.
Suppose that $K^w$ is obtained from $K^u$ by  
a symmetric $(v_1,v_2, \varepsilon )$-operation and that
$K^w$ is bipartite.
Then 
$\lambda (K^w) \ge \lambda (K^u)$. 
\es

Now using {\bf \ref{biplambda(K^w)>}}, we can obtain the following generalization of {\bf \ref{bipartite,diamondsuit}} for a symmetric $K_{xy}$-operation.
\bs {\em \cite{Kcourse}}
\label{bipartite,K_xy}
Let $G \in {\cal G}$ and $G'$ be obtained from $G$ by 
a symmetric $K_{xy}$-operation.
If $G'$ is a bipartite graph, then $G \succeq ^\phi   G'$.
\es

The notion of a vertex comparable (or threshold) graph can also be naturally generalized to weighted graphs and digraphs.
Let, as above,   $K^w$ be a weighted graph with the vertex set $V$, and so $w$ is a function from ${V \choose 2}$ to  $\mathbb{R}$. For $x, y \in V$,
$x \ne y$, we write $x \trianglelefteq y$ if $w(xv) \le w(yv)$ for every $v \in V \setminus \{x,y\}$. 
We call  $K^w$  a {\em vertex comparable $($or threshold$)$ weighted graph} if for every two distinct vertices $x, y$ of $K^w$ either 
$x  \trianglelefteq y$ or $y  \trianglelefteq x$.
Many results on threshold graphs (in particular, many results described above) can be generalized to weighted threshold graphs and digraphs.


\newpage


\begin{thebibliography}{99}


\bibitem{BM} J.A. Bondy and U.S.R. Murty, 
{\em Graph 
Theory}, Springer, 2007.

\bibitem{B} F.T. Boesch, On unreliability polynomials and graph connectivity in reliable network synthesis,
{\em J. Graph Theory} {bf 10} No. 3 (1986) 339-352.

\bibitem{CH} V. Chvatal and P. Hammer,  Aggregation of inequalities in integer programming,
{\em Ann. Discrete Math.} {\bf 1} (1977) 145-162.

\bibitem{CLBS} D. G. Corneil, H. Lerchs, L. Burlingham, L. Stewart, "Complement reducible graphs", Discrete Applied Mathematics 3 (3) (1981) 163Ð174.

\bibitem{CDS} D.M. Cvetkovi\'c, M. Doob, and H. Sacks,
{\em Spectra of Graphs},
Academic Press, 3rd ed., Johann Ambrosius Barth, 1995.

\bibitem{D} R. Deistel, {\em Graph Theory}, Springer-Verlag, New York, 2005.

\bibitem{DKZ} E. Dinic, A. Kelmans, M. Zaitcev, 
Non-isomorphic trees with the same T-polynomials. 
 {\em  Information Processing 
Letters} {\bf 6} (3), (1977) 73-76.


\bibitem{G} F. R. Gantmacher, {\em The Theory of Matrices}, Chelsea, New York, 1959.



\bibitem{HK} P. Hammer and A. Kelmans,
Laplacian spectra and spanning  trees of threshold
graphs, 
DIMACS Technical Report 92-49, Rutgers University,1992, 1-18.
%


\bibitem{H} F. Harary, {\em Graph Theory}, 
Addison-Wesley,  Reading  MA (1969).

\bibitem{HKS}  N. Hartsfield, A. Kelmans, and Yun-Qiu Shen, On the Laplacian polynomoal of a $k$-cube extension. {\em Congressus Numerantium}, 119 (1996) 73-77.

\bibitem{HJ} R.A. Horn, C.R. Johnson, {\em Matrix Theory}, Cambridge University Press, New York, 1992.

\bibitem{JK} D. Jakobson and A. Kelmans, 
The Reciprocity theorem for the Laplacian spectra of regular bipartite graphs, Manuscript,  1997.

\bibitem{Jung} H. A. Jung, On a class of posets and the corresponding comparability graphs, {\em  Journal of Combinatorial Theory}, B 24 (2) (1978) 125Ð133

\bibitem{Kar} M. Karonski, Random graphs,
{\em J. Graph Theory} {\bf 6} (1982) 349-389.

\bibitem{Kcross} A. Kelmans, Crossing properties of reliability functions of a graph.
{\em J. Graph Theory} {\bf 35} (2000) 206-221
(submitted in September, 1995).

\bibitem{KgraphTrnsf} A. Kelmans, Transformations of a
graph increasing its Laplacian polynomials and the
number of trees.  
{\em Europian Journal of Combinatorics}, {\bf 18} (1997) 35-48
(see also RUTCOR report RRR 18-95, May 1995). 

\bibitem{Kmxtr} A. Kelmans, On graphs with the
maximum number of spanning trees. 
 {\em Random
Structures and Algorithms}, {\bf 9} (1996) 177-192
[Correction on page 189: replace in $(s1)$  
$D_n^m \in {\cal S}_0$  by $D_n^m \in {\cal S}_1$ and $(0,1)$ by $(0,2)$ and replace in $(s2)$ $(0,2)$ by $(0,1)$].

\bibitem{Kcourse}  A. Kelmans, Lectures on Algebra and Combinatorics, 
Rutgers University, 1992 - 1993 and
University of Puerto Rico, 1995 - 2009.


\bibitem{Kproceed} A. Kelmans, On network reliability and related topics,
{\em 
%
5-th All--Union School--Seminar 
on Distributed Automatic Systems}, IPU, CNII ASU, Moscow, 1988.

\bibitem{KK} A. Kelmans and B. Kimelfeld,
Multiplicative submodularity of a matrix's principal
minor as a function of the set of the rows and some
combinatorial applications. 
{\em  Discrete Mathematics} {\bf 44} (1983) 113-116.

\bibitem{Krandgr} A. Kelmans, On graphs with
randomly deleted edges. 
{\em  Acta Math.Acad. Sci. Hung.} {\bf 37} (1-3), (1981) 259-267.

\bibitem{Kextr} A. Kelmans, Graphs with an
extremal number of spanning trees. 
{\em  J. Graph Theory} {\bf 4} (1980) 119-122
[Correction on page 121: in line (-1) replace ``Ref. 5'' by ``Ref. 1'' and in line  (-3) remove ``have the fewest vertices, and''].

\bibitem{Knotuniform} A. Kelmans, The graph with the
maximum probability of connectivity depends on the
edge removal probability. 
{\em  Graph Theory
Newsletters} {\bf 9}  (1), (1979) 2-3.


\bibitem{KoperProb}  A. Kelmans, Comparison of
graphs by their probability of connectivity.
\\
In:
{\em  Kombinator. Asympt. Analiz}, Krasnoyarsk
(1977) 69 - 81.


\bibitem{Kcmpr} A. Kelmans, 
Comparison of graphs by their number of spanning trees. 
{\em  Discrete Mathematics} {\bf 16} (1976) 241-261.



\bibitem{KoperTrees}  A. Kelmans, Operations on graphs
increasing the number of their spanning trees. 
In:
{\em  Issledovanie po Discretnoy Optimizacii},
Nauka, Moscow (1976) 406-424.
 

\bibitem{Ktree-forest} A. Kelmans, The number of graph
spanning trees containing a given forest.
{\em  Acta Math. Acad. Sci. Hungar.} {\bf 27} (1-2) (1976)
89-95.


\bibitem{KCh} A. Kelmans and V. Chelnokov, 
A certain polynomials of a graph and graphs with an extremal numbers of trees. 
{\em  J. Combinatorial Theory $(B)$}
{\bf 16} (1974) 197-214.

\bibitem{Kcmpr2} A. Kelmans, Comparison of graphs by their Laplacian polynomials. Manuscript 1972.

 \bibitem{KeigvMtr}  A. Kelmans, On the determination
of the eigenvalues of some special  matrices. 
{\em 
Ekonomika i Matematicheskie Metodi}  {\bf 8} (2),
(1972) 266-272.


 
\bibitem{Krandgraph} A. Kelmans, Studies on random graphs, In: {\em Trudy 4  Vsesoys. Sov. po Avtomat. Uprav.} Nauka, Moscow (1972).




\bibitem{KelLatvia}  A. Kelmans, Laplacian polynomials, spectra,  the reliability, and the number of spanning trees of graphs,
Lectures at the Winter Workshop on Discrete Mathematics, Vaivary, Latvia, 1970.

\bibitem{Kpolinom}  A. Kelmans, On properties of the
characteristic polynomial of a graph. 
In {\em Kiber. na Sl. Kom.} 
{\bf4}, Energiya,
Moscow--Leningrad (1967).
%

\bibitem{Kprob-network}  A. Kelmans, On the
connectivity of probabilistic networks. 
{\em 
Avtomat. i Telemeh. $($Automat. Remote Control$)$} {\bf
28} (3) (1967) 98-116.

\bibitem{Kasympt} A. Kelmans, On analysis and synthesis of probabilistic networks, {\em Selected Transl. in Math. Statist. and Probability}, Volume 16, (1985) 127-138
[Thanslation of Adaptive Systems, Large Systems (Proc. First All-Union Symp. Statist. Problems in Engrg. Cybernetics, Moscow, 1967; N.A. Andreev et al, editors), ``Nauka'', Moscow, 1971. pp. 264-273 ({\bf MR 52} \#7729)].

\bibitem{Ktree2} 
A. Kelmans, The number of trees in a graph  II. 
{\em  Avtomat. i Telemeh.}  
 {\bf 2} (1966) 56-65 
(English translation in   {\em Automat. Remote Control}, 
{\bf 27} (1966)).


\bibitem{Ktree1} 
A. Kelmans, The number of trees in a graph I. 
{\em  Avtomat. i Telemeh.}  
{\bf 12} (1965) 2194-2204  
(English translation in   {\em Automat. Remote Control}, {\bf 26} (1965)).



\bibitem{KuniGr} A. Kelmans, The 
distinguishing capability of the characteristic polynomial of
the conductive matrix of a graph, manuscript,
1963 (described by Cvetkovi\'c D.M., Graphs and Their Spectra. 
{\em   Elektrohn. Fak., Ser. Mat. Fiz.}, No. 354 - No. 356. (1971) 1-50).

\bibitem{Lrch}
H. Lerchs, 
On cliques and kernels,
Tech. Report {\em Dept. of Comp. Sci., Univ. of Toronto} 1971


\bibitem{LP} M. Lomonosov and V. Polessky, On the maximum of the probability of connectivity of a graph,
{\em Problemy Peredachi Informacii} {\bf 8} , (1972) 68-73.

\bibitem{MP} N.V.R. Mahadev and U.N. Peled,
{\em Threshold graphs and Related Topics}, Annals of Discrete Mathematics 56, North-Holland, Amsterdam 1995.

\bibitem{RW} R. C. Read and R. J. Wilson,
{\em An Atlas of Graphs}, Oxford University Publications,
1998.

\bibitem{Sei} D. Seinsche, On a property of the class of n-colorable graphs", {\em Journal of Combinatorial Theory},  B 16 (2) (1974) 191Ð193.

\bibitem{Sum} D. P. Sumner, Dacey graphs, {\em J. Austral. Math. Soc. }18 (04)(1974) 492Ð502.

\bibitem{W} H. Weiner, Structural determination of paraffin boiling poins. {\em J. Am. Chem. Soc.} {\bf 69(1)} (1947), 
17-20.






\end{thebibliography}
\end{document}